# GENERALIZED DEGREES AND DENSITIES FOR FAMILIES OF SETS

by

EMANUEL KNILL

A thesis submitted to the
Faculty of the Graduate School of the
University of Colorado in partial fulfillment
of the requirements for the degree of
Doctor of Philosophy
Department of Mathematics
1991

This thesis for the Doctor of Philosophy degree by
Emanuel Knill
has been approved for the
Department of
Mathematics
by

_______________________________

Richard Laver

_______________________________

Andrzej Ehrenfeucht

Date _______________________________

Knill, Emanuel (Ph. D., Mathematics)
Generalized Degrees and Densities for Families of Sets
Thesis directed by  Richard Laver

Let $\mathcal{F}$ be a family of subsets of $\{1, 2, \ldots, n\}$. The *width-degree* of an element $x \in \bigcup \mathcal{F}$ is the width of the family $\{U \in \mathcal{F} | x \in U\}$. If $\mathcal{F}$ has maximum width-degree at most $k$, then $\mathcal{F}$ is *locally k-wide*.

Bounds on the size of locally $k$-wide families of sets are established. If $\mathcal{F}$ is locally $k$-wide and centered (every $U \in \mathcal{F}$ has an element which does not belong to any member of $\mathcal{F}$ incomparable to $U$), then $|\mathcal{F}| \leq (k+1)(n - \frac{k}{2})$; this bound is best possible. Nearly exact bounds, linear in $n$ and $k$, on the size of locally $k$-wide families of arcs or segments are determined. If $\mathcal{F}$ is any locally $k$-wide family of sets, then $|\mathcal{F}|$ is linearly bounded in $n$. The proof of this result involves an analysis of the combinatorics of antichains.

Let $P$ be a poset and $L$ a semilattice (or an intersection-closed family of sets). The *P-size* of $L$ is $|L^P|$. For $u \in L$, the *P-density* of $u$ is the ratio $\frac{|(u)^P|}{|L^P|}$. The *density* of $u$ is given by the [1]-density of $u$. Let $p$ be the number of filters of $P$. $L$ has the *P-density property* iff there is a join-irreducible $a \in L$ such that the $P$-density of $a$ is at most $\frac{1}{p}$.

Which non-trivial semilattices have the $P$-density property? For $P = [1]$, it has been conjectured that the answer is: "all" (the union-closed sets conjecture). Certain subdirect products of lower-semimodular lattices and, for $P = [n]$, of geometric lattices have the $P$-density property in a strong sense. This generalizes some previously known results. A fixed lattice has the $[n]$-density property if $n$ is large enough. The density of a generator $U$ of a union-closed family of sets $L$ with $\emptyset \in L$ is estimated. The estimate depends only on the local properties of $L$ at $U$. If $L$ is generated by sets of size at most two, then there is a generator $U$ of $L$ with estimated density at most $\frac{1}{2}$.

# CONTENTS

# CHAPTER 1

# INTRODUCTION

## 1.1 To the Reader

Section 1.2 is an overview of the contents of this work.

Chapter 2 is a summary of the background material, and contains definitions and some fundamental results from the combinatorics of finite sets and the theory of posets. Most of the definitions and terminology are standard. Table 2.2 in Section 2.1 contains a list of notation and terminology.

Except in Chapters 1 and 2, most definitions appear after the paragraph heading *Definition*. Terms are italicized when first defined. *Observations* are statements which are true by definition or follow from the discussion and do not require further proof.

Chapters 3 and 4 are independent and can be read in either order.

The *Notes* at the end of each chapter contain a discussion of the literature relevant to each section, attributions of results and references to related work.

## 1.2 Overview

The study of extremal properties of families of sets involves determining the relationships between their various size measures. Two fundamental size measures of a family of sets $\mathcal{F}$ are: $|\mathcal{F}|$, the cardinality of the family; and $|\bigcup \mathcal{F}|$, the number of elements in the union of its members. Many size measures are obtained by considering $\mathcal{F}$ as a poset ordered by inclusion. One such measure is the width of the family $\mathcal{F}$. Another measure is the *P-size* of $\mathcal{F}$ given by the number of order-preserving maps from the poset $P$ into $\mathcal{F}$. For example if $P$ is the $n$-element chain, then the $P$-size of $\mathcal{F}$ is the number of multichains of length $n-1$ of $\mathcal{F}$.

Let $w$ be a size measure for families of sets. Let $\mathcal{F}_x = \{U \in \mathcal{F} \mid x \in U\}$. The *domain* of $\mathcal{F}$ is the set $\bigcup \mathcal{F}$. The measure $w$ induces a generalized degree measure as follows: If $x$ is an element of the domain of $\mathcal{F}$, then the $w$-degree of $x$ is $w(\mathcal{F}_x)$. If $\mathcal{F}$ is a graph (i.e. every member of $\mathcal{F}$ is a pair) and $w$ is the cardinality function, then the $w$-degree coincides with the usual notion of degree for graphs.

The size measure $w$ also induces a generalized notion of density. If $\mathcal{G}$ is a subfamily of $\mathcal{F}$, then the *w-density* in $\mathcal{F}$ of $\mathcal{G}$ is the ratio $\frac{w(\mathcal{G})}{w(\mathcal{F})}$. If $x$ is an element of the domain of $\mathcal{F}$, then the $w$-density in $\mathcal{F}$ of $x$ is the $w$-density in $\mathcal{F}$ of $\mathcal{F}_x$.

Much research has been devoted to the extremal theory of sets. However, many generalized degrees and densities have only recently been introduced. There are some well known open problems concerning these notions, in particular concerning the notion of density. Perhaps the most famous is the union-closed sets conjecture due to P. Frankl (Conjecture 4.3.7). A list of others can be found in Section 4.12.

This work consists of two parts. The first introduces the study of the width-degree. The second examines $P$-densities in intersection- and union-closed families of sets (i.e. semilattices).

**The width-degree.** A family of sets of maximum width-degree at most $k$ is called *locally $k$-wide*. Except for Knill et al. [27], there has been no published work on locally $k$-wide families of sets. The concept was first introduced by Ehrenfeucht and Haussler as part of the definition of pseudotrees. Their aim was to find structures more general than trees for classifying textual information.

A *tree* on the set $X$ is a locally 1-wide family $\mathcal{F}$ of subsets of $X$ such that $\mathcal{F}$ contains $X$ and the singletons of $X$. Since $\mathcal{F}$ is locally 1-wide iff every pair of members of $\mathcal{F}$ is either comparable or disjoint, a tree on $X$ has the familiar tree-like structure where $X$ is the root and the singletons of $X$ are the leaves. A tree $\mathcal{F}$ also has the property that every $U \in \mathcal{F}$ has an element not contained in any member of $\mathcal{F}$ incomparable to $U$. Families of sets with this property are called *centered*. Thus centered, locally $k$-wide families of sets are a natural generalization of trees. Such families are called *$k$-pseudotrees*.

In order for $k$-pseudotrees to be useful from a data structure perspective, they should not be too large relative to the size of the domain. Ehrenfeucht and Haussler asked the following question: Is the size of a $k$-pseudotree linearly bounded in the size of the domain? The special case for 2-pseudotrees of segments was solved by Knill in 1990. Knill and Ehrenfeucht subsequently showed that if the domain of a $k$-pseudotree $\mathcal{F}$ has $n$ elements where $n > k$, then $|\mathcal{F}| \leq (k+1)(n-\frac{k}{2})$ (Theorem 3.4.1). This bound is best possible. Linear bounds have now been obtained for locally $k$-wide families of arcs (Theorem 3.6.1), locally $k$-wide families of segments (Theorem 3.6.6) and arbitrary locally $k$-wide families of sets (Theorem 3.8.1). The results are summarized in Table 1.1, where $f(k, n)$ is the maximum size of a locally $k$-wide family of subsets of $\{1, 2, \ldots, n\}$ in the given class.

The proof of the linear bound for locally $k$-wide families of sets requires an analysis of the combinatorics of certain antichains of posets. An antichain $A$ of the poset $P$ is a *maximal $r$-antichain* iff the width of the filter $[A]$ generated by $A$ is $r$ and $A$ is the only antichain of size $r$ in $[A]$. A *maximal $*$-antichain* is a maximal $r$-antichain for some $r$. The family of maximal $*$-antichains is a meet-subsemilattice of the family of antichains of $P$ with the filter order (defined by $A \leq B$ iff $B \subseteq [A)$ ) (Theorem 3.7.3). If $A$ and $B$ are incomparable $*$-antichains, then $|A \wedge B| > \max(|A|, |B|)$ (Corollary 3.7.4). If $w$ is the width of $P$, then the union of the $*$-antichains consists of at most $\frac{w(w+1)}{2}$ elements (Theorem 3.7.8).

**$P$-densities.** Let $P$ be a poset and let $p$ be the number of filters of $P$. An intersection-closed family of sets $L$ has the *$P$-density property* iff there is an element $x$ in the domain of $L$ such that the $P$-density of $x$ is at most $\frac{1}{p}$. The *density property* is the $[1]$-density property (note that for $P = [1]$, $p = 2$).

The problem discussed in the second part of this work is the following:

**Problem 1.2.1** *What intersection-closed families of sets have the $P$-density property?*

One version of the union-closed sets conjecture asserts that for every intersection-closed family $L$ with at least two members, there is an element in the domain of $L$ contained in at most one half of the members of $L$ (see Duffus [14]). If true, this implies that all non-trivial intersection-closed families have the density property.

If the intersection-closed family $L$ is trivial (i.e. $|L| = 1$), then $L$ does not have the $P$-density property for any poset $P$. No other examples without the $P$-density property are known.

Every intersection-closed family is a meet-semilattice, and conversely every meet-semilattice can be represented by an intersection-closed family (e.g. the family of principal ideals of the semilattice). Therefore one can use lattice theoretical techniques to analyze Problem 1.2.1. The notion of an element in the domain is replaced by that of a (proper) join-irreducible. For the remainder of this section, all (semi-)lattices, intersection- and union-closed families are assumed to have at least two members.

Known results related to the union-closed sets conjecture are:

  (i) If the intersection-closed family $L$ contains $(\bigcup L) \setminus \{x, y\}$ for some $x, y \in \bigcup L$, then $L$ has the density property.

 (ii) If the average size of the members of the intersection-closed family $L$ is at most $\frac{1}{2} |\bigcup L|$, then $L$ has the density property.

(iii) If the meet-semilattice $L$ has the density property, then so does $L \times M$ for any meet-semilattice $M$.

 (iv) Modular and geometric lattices have the density property.

  (v) If $L$ is a meet-semilattice and $n$ is the maximum cardinality of $L \setminus [a)$ for join-irreducibles $a$, then $|L| \leq Cn \log n$ for some constant $C$.

Results (i) and (ii) are obtained by double counting. They can be generalized and used as a starting point to establish lower bounds on counter-examples to the union-closed sets conjecture (see Sarvate and Renauld [38]). Result (iii) follows from the observation that if $a$ is a join-irreducible of $L$, then $\langle a, \hat{0}_M \rangle$ is a join-irreducible of $M$. Results (iv) and (v) are of unknown origins. They are generalized in Sections 4.7 and 4.10 respectively.

To generalize (iii) and (iv), the idea of decreasing matchings among different types of order-preserving maps in $L^P$ is introduced. Let $F$ be a filter of $P$ and $a$ a join-irreducible of $L$. Define $\mathrm{T}(L^P, F, a) = \{f \in L^P \mid f(x) \geq a \text{ iff } x \in F\}$. The order-preserving map $f : P \to L$ is of *type* $(F, a)$ iff $f \in \mathrm{T}(L^P, F, a)$. The pair $(L, a)$ has the *top matching property* for $P$ iff for every filter $F$ of $P$ there is a one-to-one map $\sigma : \mathrm{T}(L^P, P, a) \to \mathrm{T}(L^P, F, a)$ such that $\sigma(f) \leq f$. The semilattice $L$ has the *top matching property* iff there is a join-irreducible $a \in L$ such that $(L, a)$ has the top matching property. The matching property implies the $P$-density property. There are non-trivial semilattices which do not have the matching property for any poset (Example 4.2.8).

In Section 4.7 it is shown that lattices with a lower-semimodular coatom have the matching property for any poset (Theorem 4.7.4). This includes all lower-semimodular and in particular all modular lattices. It is also shown that (non-trivial) geometric lattices have the matching property for linearly ordered sets (Theorem 4.7.9).

The matching property is preserved not only by direct products, but also by certain subdirect products which preserve local properties of lattices (Section 4.5). To make the notion of locality precise, lattices are considered as union-closed families of sets (Section 4.6). The join-irreducibles of the lattice $L$ form a family of sets $\mathcal{J}$ which generates $L$ (by union). The *lattice neighborhood* $\mathcal{N}_L(U)$ of $U \subseteq \bigcup \mathcal{J}$ is the union-closed family generated by $\{\emptyset\} \cup \{A \in \mathcal{J} \mid A \cap U \neq \emptyset\}$. If $U \in \mathcal{J}$ and $\mathcal{N}_L(U)$ is a geometric lattice, then $(L, U)$ has the matching property for the $n$-element chain. This and similar results appear in Section 4.7.

In Section 4.11 a technique is developed for estimating the density of a join-irreducible $a$ of $L$ which depends only on $\mathcal{N}_L(\bigcup \mathcal{N}_L(a))$. This estimate is used to show that if the lattice $L$ (considered as a union-closed family of sets) is generated by a graph, then $L$ has the density property (Theorem 4.11.2). This implies that for every graph $G$ there is an edge contained in at most half of the unions of edges of $G$.

Result (v) listed above is generalized in Section 4.10: If for every join-irreducible $a$ the number of order-preserving maps not of type $(P, a)$ is at most $n$, then $\left| L^P \right| \leq M(n)$, where $M(n) \sim n \log_p(n)$ asymptotically. The bound is exact for a generalization of the notion of lattices to multi-semilattices, where every non-maximal member of the lattice is assigned a multiplicity.

In Section 4.9 an expression for the Zeta polynomial is used to show that if $L$ is a fixed non-trivial lattice and $n$ is sufficiently large, then $L$ has the $[n]$-density property (where $[n]$ is the $n$-element chain).

Table 1.1

| Class | $f(k, n)$ | Location |
|---|---|---|
| $U \in \mathcal{F} \Rightarrow |U| = r$ | $\leq \lfloor \frac{kn}{r} \rfloor$ | Theorem 3.3.5. |
| All families | $\leq 2 + n + kn \ln(n)$ | Theorem 3.3.6. |
| | $\leq (2k)^{k-1} n$ | Theorem 3.8.1 |
| | $\geq 3kn - o(kn)$ if | |
| | $1 = o(k)$ and $k = o(n)$ | Example 3.8.2. |
| Centered families | $(k + 1)(n - \frac{k}{2})$ if | |
| | $n > k$ | Theorem 3.4.1. |
| Families of segments | $2kn - 2k^2 + k + 1$ if | |
| | $n \geq 2k$ | Theorem 3.6.5. |
| Families of arcs | $\leq 2kn - k + 1$ | Theorem 3.6.1. |
| | $\geq 2kn - k^2 - k + 2$ if | |
| | $n \geq k$ | Example 3.6.2. |

# CHAPTER 2

# PRELIMINARIES

## 2.1 Conventions and Notation

All sets and structures are finite by default. On the few occasions where an infinite set is used, this will be stated explicitly.

Subscripts and function arguments will be omitted when it is possible to do so without loss of clarity.

The guidelines for symbol use are given in Table 2.1.

Table 2.1: Guidelines for symbol use.

| | |
|---|---|
| $\mathcal{C}, \mathcal{L}$ | classes of structures |
| $\mathcal{A}, \mathcal{B}$ | antichains of sets |
| $\mathcal{F}, \mathcal{G}, \mathcal{H}$ | families of sets |
| $A, B, C \ldots$ | sets, antichains |
| $F, G, H$ | filters |
| $L, M$ | semilattices |
| $P, Q$ | posets |
| $U, V, W$ | sets, members of a family of sets |
| $X, Y, Z$ | sets, domains |
| $a, b, c \ldots$ | elements, join-irreducibles, atoms etc. |
| $i, j, k \ldots$ | integers |
| $u, v, w$ | elements of a semilattice |
| $x, y, z$ | arbitrary elements, real numbers |

A list of notation and terminology is given in Table 2.2.

## 2.2 Sets, Arcs and Segments

An *n-set* is a set with exactly $n$ elements. The sets $A$ and $B$ *intersect* iff $A \cap B$ is non-empty. $A$ and $B$ are *comparable* iff $A \subseteq B$ or $B \subseteq A$. $A$ and $B$ are *incomparable* iff they are not comparable. $A$ and $B$ *overlap* iff they intersect and are incomparable.

If $i$ and $j$ are integers, then $[i, j]$ denotes the *segment* of integers $k$ such that $i \le k \le j$. If $j < i$, then $[i, j]$ is the empty segment. For $i \le j$, the *left* and *right endpoints* of $[i, j]$ are $i$ and $j$ respectively. A segment of $[n]$ is a segment $[i, j]$ with $[i, j] \subseteq [n]$. Let $[i, j]$ and $[k, l]$ be segments. Note that if $k = i$ or $j = l$, then the two segments are comparable. The segment $[i, j]$ is *to the left* of $[k, l]$ iff $[i, j]$ and $[k, l]$ are incomparable and $i < k$. The intersection of $[i, j]$ and $[k, l]$ is also a segment. If $[i, j]$ is to the left of $[k, l]$, then their intersection is the segment $[k, j]$ and is called the *left overlap* of $[i, j]$ with $[k, l]$. If this intersection is non-empty, then $[i, j]$ *overlaps* $[k, l]$ *from the left*. The *right overlap* of two segments is defined similarly.

Table 2.2: Notation.

| | |
|---|---|
| **Z** | the integers |
| **N** | the non-negative integers |
| **P** | the strictly positive integers |
| **R** | the real numbers |
| $[n]$ | the set $\{1, 2, \ldots, n\}$ with the linear order |
| $[i, j]$ | the set of integers $k$ such that $i \leq k \leq j$ |
| $\mathbf{Z}_n$ | the integers modulo $n$, represented by $\{0, 1, \ldots, n-1\}$ |
| $\emptyset$ | the empty set |
| $|X|$ | the number of elements (cardinality, size) of the set $X$ |
| $\lfloor x \rfloor$ | the greatest integer $\leq x$ |
| $\lceil x \rceil$ | the least integer $\geq x$ |
| $\log(x)$ | the logarithm base two of $x$ |
| $\ln(x)$ | the logarithm base $e$ of $x$ |
| $n!$ | the factorial of $n$: $n! = 1 \cdot 2 \cdot \ldots \cdot n$ |
| $(n)_k$ | the $k$'th falling factorial of $n$: $(n)_k = n(n-1)\ldots(n-k+1)$ |
| $\binom{n}{k}$ | the binomial coefficient: $\binom{n}{k} = \frac{(n)_k}{k!}$ |
| $f \circ g$ | the composition of the maps $f$ and $g$ |
| $f \restriction A$ | the restriction of the map $f$ to $A$ |
| $f = o(g)$ | means $\lim_{x \to \infty} \frac{f(x)}{g(x)} = 0$ |
| $f = O(g)$ | means that for some $c > 0$, $f(x) \leq C \cdot g(x)$ |
| $f \sim g$ | means $\lim_{x \to \infty} \frac{f(x)}{g(x)} = 1$ |
| $A \subset B$ | the set $A$ is a proper subset of $B$ |
| $A \subseteq B$ | $A = B$ or $A \subset B$ |
| $A \cap B$ | the intersection of the sets $A$ and $B$ |
| $A \cup B$ | the union of the sets $A$ and $B$ |
| $A \setminus B$ | set difference: the set of elements of $A$ not in $B$ |
| $\langle a, b \rangle$ | the ordered pair $a$ and $b$ |
| $A \times B$ | the cartesian product of the sets $A$ and $B$ |
| $\prod_{i=1}^{n} A_i$ | the cartesian product of the sets $A_1 \ldots A_n$ |
| $\pi_i$ | the projection onto the $i$'th component of a cartesian product |
| $\mathrm{d}_{\mathcal{G}}(x)$ | the degree in the graph $\mathcal{G}$ of $x$ |

Table 2.2: continued.

| | |
|---|---|
| $2^X$ | the family of subsets of $X$ |
| $\bigcup \mathcal{F}$ | the union of the members of the family of sets $\mathcal{F}$ |
| $\bigcap \mathcal{F}$ | the intersection of the members of $\mathcal{F}$ |
| $\mathcal{F} \vee \mathcal{G}$ | for families of sets $\mathcal{F}$ and $\mathcal{G}$: $\mathcal{F} \vee \mathcal{G} = \{U \cup V \mid U \in \mathcal{F}, V \in \mathcal{G}\}$ |
| $\mathcal{F} \wedge \mathcal{G}$ | for families of sets $\mathcal{F}$ and $\mathcal{G}$: $\mathcal{F} \wedge \mathcal{G} = \{U \cap V \mid U \in \mathcal{F}, V \in \mathcal{G}\}$ |
| $\mathcal{F}_{\subseteq X}$ | the family induced by $\mathcal{F}$ in $X$: $\mathcal{F}_{\subseteq X} = \{U \subseteq X \mid U \in \mathcal{F}\}$ |
| $\mathcal{F}_{\supseteq X}$ | the family induced by $\mathcal{F}$ above $X$: |
| | $\mathcal{F}_{\supseteq X} = \{U \supseteq X \mid U \in \mathcal{F}\}$ |
| $\mathcal{F}_{\cap X}$ | the restriction of $\mathcal{F}$ to $X$: $\mathcal{F}_{\cap X} = \{U \cap X \mid U \in \mathcal{F}\}$ |
| $\mathcal{F}_{\setminus X}$ | the restriction of $\mathcal{F}$ to the complement of $X$: |
| | $\mathcal{F}_{\setminus X} = \{U \setminus X \mid U \in \mathcal{F}\}$ |
| $[\mathcal{F}]_r$ | the family of members of $\mathcal{F}$ of size (or rank) $r$: |
| | $[\mathcal{F}]_r = \{U \in \mathcal{F} \mid |U| = r\}$ |
| $\mathrm{C}(U)$ | the center of $U$ in a family of sets |
| $\hat{0}_P$ | the least element of the poset $P$ (if it exists) |
| $\hat{1}_P$ | the greatest element of the poset $P$ (if it exists) |
| $[x, y]$ | for $x \leq y$ in the poset $P$, the interval of elements |
| | between $x$ and $y$: $[x, y] = \{z \in P \mid x \leq z \leq y\}$ |
| $(A]$ | the ideal generated by $A$ in the poset $P$: |
| | $(A] = \{x \in P \mid x \leq y$ for some $y \in A\}$ |
| $[A)$ | the filter generated by $A$ in the poset $P$: |
| | $[A) = \{x \in P \mid x \geq y$ for some $y \in A\}$ |
| $\bigvee A$ | the least upper bound of $A$ in a poset (if it exists) |
| $\bigwedge A$ | the greatest lower bound of $A$ in a poset (if it exists) |
| $J(L)$ | the set of proper join-irreducibles of the meet-semilattice $L$ |
| $M(L)$ | the set of proper meet-irreducibles of the join-semilattice $L$ |
| $\mathcal{P}(L)$ | the family of principal ideals of the semilattice $L$ |
| $\mathrm{Sub}(L)$ | the family of subsemilattices of $L$ |
| $P + Q$ | the disjoint union of the posets $P$ and $Q$ |
| $Q^P$ | the set of order-preserving maps from $P$ to $Q$ |
| $\mathrm{w}(P)$ | the width of the poset $P$ |
| $\mathrm{T}(L^P, F, a)$ | for the semilattice $L$, poset $P$, filter $F$ of $P$ |
| | and join-irreducible $a \in L$: |
| | $\mathrm{T}(L^P, F, a) = \{f \in L^P \mid f(x) \geq a$ iff $x \in F\}$ |
| $P \cong Q$ | means that $P$ is isomorphic to $Q$ |

The number $i$ *precedes* $j$ in $\mathbf{Z}_n$ iff $j = i + 1 \bmod(n)$. This defines the *clockwise cyclic order* on $\mathbf{Z}_n$. The number $j$ occurs *before* $k$ *clockwise from* $i$ iff $j = k$ or $j$ appears before $k$ in the sequence

$$i, i + 1 \bmod(n), i + 2 \bmod(n), \ldots, i - 2 \bmod(n), i - 1 \bmod(n).$$

In *clockwise cyclic order* $i \leq j \leq k$ means that $j$ occurs before $k$ clockwise from $i$. The expression $[i, k]$ denotes the *arc* of elements $j \in \mathbf{Z}_n$ such that in clockwise cyclic order $i \leq j \leq k$. The *left* and *right endpoints* of the arc $[i, k]$ are $i$ and $k$ respectively. By default $[i, k] \neq \mathbf{Z}_n$, so that $k + 1 \bmod(n) \neq i$.

The intersection of two incomparable arcs is not necessarily an arc, but consists of the union of at most two disjoint arcs. The arc $[i, j]$ *overlaps* $[k, l]$ *from the left* iff $k$ occurs before $j$ clockwise from $i$, and the arcs $[i, j]$ and $[k, l]$ are incomparable. The arc $[k, j]$ is called the *left overlap* of $[i, j]$ with $[k, l]$.

The *cyclic closure* of $[n]$ is obtained by letting the integer $n$ precede $1$ in $[n]$. The cyclic closure of $[n]$ is identical to $\mathbf{Z}_n$ with the clockwise cyclic order if $k \in [n]$ is identified with $k - 1 \bmod(n) \in \mathbf{Z}_n$.

## 2.3 Families of Sets

A *family of sets* $\mathcal{F}$ is a set of sets. The *domain* of $\mathcal{F}$ is given by $\bigcup \mathcal{F}$. $\mathcal{F}$ is a *family of arcs (segments)* of $\mathbf{Z}_n$ ($[n]$) iff every member of $\mathcal{F}$ is an arc (segment) of $\mathbf{Z}_n$ ($[n]$). A *hypergraph* is a pair $(\mathcal{F}, X)$, where $\mathcal{F}$ is a family of sets and $X$ includes the domain of $\mathcal{F}$. The members of $\mathcal{F}$ are the *edges* of the hypergraph. The elements of $X$ not in the domain are *isolated* elements of the hypergraph. The hypergraph $(\mathcal{F}, X)$ is *$r$-uniform* iff every member of $\mathcal{F}$ has exactly $r$ elements.

If $\mathcal{F}$ and $\mathcal{G}$ are families of sets, then $\mathcal{F} \vee \mathcal{G} = \{U \cup V \mid U \in \mathcal{F}, V \in \mathcal{G}\}$ and $\mathcal{F} \wedge \mathcal{G} = \{U \cap V \mid U \in \mathcal{F}, V \in \mathcal{G}\}$. The family of sets $\mathcal{F}$ is *union-closed* iff for every $U, V \in \mathcal{F}$, $U \cup V \in \mathcal{F}$. Equivalently, $\mathcal{F}$ is union-closed iff $\mathcal{F} \vee \mathcal{F} \subseteq \mathcal{F}$. Similarly, $\mathcal{F}$ is *intersection-closed* iff $\mathcal{F} \wedge \mathcal{F} \subseteq \mathcal{F}$.

**Observation 2.3.1** *Let $X$ be a subset of the domain of the family of sets $\mathcal{F}$. If $\mathcal{F}$ is union-closed (intersection-closed) then so are the induced families of sets $\mathcal{F}_{\subseteq X}$ and $\mathcal{F}_{\supseteq X}$, as well as the restrictions $\mathcal{F}_{\cap X}$ and $\mathcal{F}_{\setminus X}$ of $\mathcal{F}$.*

Let $\mathcal{F}$ be a union-closed family of sets. The family $\mathcal{F}$ is *generated* by $\mathcal{G}$ iff $\mathcal{F}$ consists of the unions of non-empty subfamilies of $\mathcal{G}$, that is iff $\mathcal{F} = \{\bigcup \mathcal{H} \mid \mathcal{H} \subseteq \mathcal{G}, \mathcal{H} \neq \emptyset\}$. The *generators* of $\mathcal{F}$ are the members $U$ of $\mathcal{F}$ such that if $\mathcal{G} \subseteq \mathcal{F}$ and $\bigcup \mathcal{G} = U$, then $U \in \mathcal{G}$. Note that $\mathcal{F}$ is generated by the family of generators of $\mathcal{F}$.

A *forest (of sets)* on $X$ is a family $\mathcal{F}$ of non-empty subsets of $X$ such that for every $U, V \in \mathcal{F}$, $U$ and $V$ are either comparable or disjoint. If in addition $\mathcal{F}$ contains $X$ and the singletons of $X$, then $\mathcal{F}$ is a *tree (of sets)*.

A *partition* of $X$ is a family of sets $\mathcal{F}$ with domain $X$ such that the members of $\mathcal{F}$ are pairwise disjoint. The partition $\mathcal{G}$ of $X$ is coarser than $\mathcal{F}$ iff every member of $\mathcal{F}$ is included in a member of $\mathcal{G}$.

## 2.4 Partially Ordered Sets

A *partially ordered set* $P$ (*poset*, for short) is a non-empty set (also denoted by $P$) together with a reflexive, transitive and antisymmetric binary relation $\leq$ on $P$. Thus for every $x, y, z \in P$, $x \leq x$ (reflexivity); $x \leq y$ and $y \leq z$ imply that $x \leq z$ (transitivity); and $x \leq y$ and $y \leq x$ imply that $x = y$ (antisymmetry). The relation $\leq$ is the *partial order* on $P$. The expression $x < y$ means $x \neq y$ and $x \leq y$. The element $x$ is *below* $y$ iff $x \leq y$, and $x$ is *above* $y$ iff $y \leq x$.

Let $P$ be a poset. The *dual* partial order of $P$ is obtained by reversing the partial order of $P$. The set $P$ together with the dual partial order is denoted by $P^*$ (the *dual* of $P$). For $x, y \in P^*$, $x \leq y$ in $P^*$ iff $y \leq x$ in $P$. Whenever appropriate, the symbols $\leq$ and $<$ are reversed to denote the dual partial order.

If $\mathcal{F}$ is a family of sets, the subset relation induces the *inclusion order* on $\mathcal{F}$ defined by $U \leq V$ iff $U \subseteq V$. When terminology for posets is used in the context of a family of sets, the intended partial order is the inclusion order (unless otherwise specified).

Let $P$ be a poset. The elements $x$ and $y$ of $P$ are *comparable* iff $x \leq y$ or $y \leq x$; otherwise, $x$ and $y$ are *incomparable*. The element $y$ is *between* $x$ and $z$ iff $x \leq y \leq z$. The element $x$ *covers* $y$ in $P$ iff $y < x$ and the only elements between $x$ and $y$ are $x$ and $y$. The cover relation of $P$ defines a directed graph, the *Hasse diagram* of $P$, obtained by connecting $x$ to $y$ iff $x$ covers $y$. The Hasse diagram of $P$ completely determines the partial order of $P$ (assuming that $P$ is finite). Pictorial representations of the Hasse diagram, where a downward edge from $x$ to $y$ means that $x$ covers $y$, are frequently used to represent posets.

The poset $P$ has the *discrete order* iff for every $x, y \in P$ with $x \neq y$, $x$ and $y$ are incomparable. $P$ is a *linearly ordered set* (or a *chain*) iff for every $x, y \in P$, $x \leq y$ or $y \leq x$. Thus, $[n]$ is linearly ordered by the usual ordering of the integers. Figure 2.1 shows the Hasse diagrams of $\{1, 2, 3, 4\}$ with the discrete order, $2^{\{1,2\}}$ and $[4]$.

The partial order $<'$ is an extension of the partial order $<$ of $P$ iff for every $x, y \in P$, $x < y$ implies that $x <' y$. The following theorem is one of the fundamental results in the theory of posets (see *Notes* at the end of this chapter).

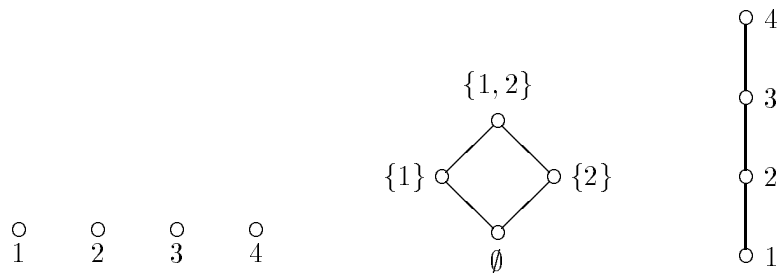

Figure 2.1: $\{1, 2, 3, 4\}$ (with the discrete order), $2^{\{1,2\}}$ and [4]

**Theorem 2.4.1** *For every poset $P$ there exists a linear order which extends the partial order of $P$.*

Let $P$ be a poset. If $x$ and $y$ are elements of $P$ and $x \le y$, then the *interval* from $x$ to $y$, denoted by $[x, y]$, consists of the elements of $P$ between $x$ and $y$. The subset $Q$ of $P$ is *convex* iff for every $x$ and $y$ in $Q$ such that $x \le y$, $[x, y] \subseteq Q$. The order *induced* by $P$ on the subset $Q$ is the restriction of the partial order of $P$ to $Q$. Thus, if $x, y \in Q$, then $x \le y$ in $Q$ iff $x \le y$ in $P$. A *maximal* element of $Q$ is an element $x \in Q$ such that for every $y \in Q$, $x \not< y$. Dually, a *minimal* element of $Q$ is an element $x \in Q$ such that for every $y \in Q$, $y \not< x$. Observe that $x$ is minimal in $Q$ iff $x$ is maximal in $Q^*$. If $x$ is the only maximal element of $Q$, then $x$ is called the *greatest* element of $Q$. If $x$ is the only minimal element of $Q$, then $x$ is called the *least* element of $Q$. If $P$ has a greatest element, this element is denoted by $\hat{1}_P$. Similarly, the least element of $P$, if it exists, is denoted by $\hat{0}_P$. An *atom* is an element which covers $\hat{0}$; a *coatom* is an element which is covered by $\hat{1}$.

The element $x$ of $P$ is an *upper bound* of $Q$ iff for every element $y$ of $Q$, $y < x$. The element $x$ is the *least upper bound* of $Q$ iff $x$ is the least element of the set of upper bounds of $Q$. If it exists, the least upper bound of $Q$ is denoted by $\bigvee Q$. *Lower bounds* and *greatest lower bounds* of $Q$ are defined dually. If it exists, the greatest lower bound of $Q$ is denoted by $\bigwedge Q$. If $Q = \{x, y\}$, then $\bigvee Q$ and $\bigwedge Q$ are also denoted by $x \vee y$ (the *join* of $x$ and $y$) and $x \wedge y$ (the *meet* of $x$ and $y$) respectively. To avoid ambiguity, $\vee_P$ and $\wedge_P$ will be used to denote the join and the meet operations in $P$, when necessary.

The basic properties of the (in general partial) join and meet operations are as follows:

- $\vee$ and $\wedge$ are idempotent, i.e. $x \vee x = x$ and $x \wedge x = x$.

- $\vee$ and $\wedge$ are symmetric, i.e. $x \vee y = y \vee x$ if either side exists; and similarly for $\wedge$.

- $\vee$ and $\wedge$ are associative, i.e. if $x \vee y$, $y \vee z$, $(x \vee y) \vee z$ and $x \vee (y \vee z)$

all exist, then $(x \vee y) \vee z = x \vee (y \vee z)$; and similarly for $\wedge$.

- $\vee$ and $\wedge$ satisfy the absorption identities, i.e. $(x \vee y) \wedge x = x$ (if $x \vee y$ exists) and $(x \wedge y) \vee x = x$ (if $x \wedge y$ exists).

If $P$ has a greatest element $\hat{1}$ then $\bigwedge \emptyset = \hat{1}$ and if $P$ has a least element $\hat{0}$ then $\bigvee \emptyset = \hat{0}$.

A *(proper) chain* of $P$ of length $n$ is a strictly increasing sequence $x_0 < x_1 < \ldots < x_n$ of elements of $P$. The maximum length of a chain of $P$ is called the *height* of $P$. A *multichain* of *length* $n$ is an increasing sequence $x_0 \leq x_1 \leq \ldots \leq x_n$ of elements of $P$. Note that a (multi)chain of length $n$ has $n + 1$ elements.

An *antichain* of $P$ is a subset $A$ of $P$ such that every pair of distinct elements of $A$ are incomparable. The maximum cardinality of an antichain of $P$ is called the *width* of $P$ and is denoted by $\mathrm{w}(P)$. An antichain of $P$ with $\mathrm{w}(P)$ elements is called a *Sperner antichain*.

If $P$ is the disjoint union of $n$ chains, then the width of $P$ is at most $n$. The converse is given by the following result due to Dilworth [10]:

**Theorem 2.4.2** *If $P$ is a poset of width $w$, then $P$ is the disjoint union of $w$ chains.*

An *(order) ideal* of the poset $P$ is a subset $I$ of $P$ such that if $x \in I$ and $y \leq x$ in $P$, then $y \in I$. An *(order) filter* of $P$ is a subset $F$ of $P$ such that if $x \in F$ and $y \geq x$ in $P$, then $y \in F$. The family of ideals of $P$ is both intersection- and union-closed, as is the family of filters of $P$.

If $A$ is an antichain of $P$ and $x \in P$, then $x$ is *below* $A$ iff there exists $y \in A$ such that $x \leq y$. The set of elements of $P$ below $A$ is the ideal *generated* by $A$ and is denoted by $(A]$. The element $x$ of $P$ is *above* $A$ iff there exists $y \in A$ such that $x \geq y$. The set of elements of $P$ above $A$ is the filter *generated* by $A$ and is denoted by $[A)$. For $x \in P$, the *principal* ideal generated by $x$ is $(x] = (\{x\}]$ and the *principal* filter generated by $x$ is $[x) = [\{x\})$.

Let $\mathcal{A}$ be the family of antichains of $P$. Identifying each antichain $A$ with $(A]$ induces a partial order on $\mathcal{A}$, where $A \leq B$ iff $(A] \subseteq (B]$. This is called the *ideal order* of $\mathcal{A}$. The *filter order* is the partial order of $\mathcal{A}$ defined by $A \leq B$ iff $[A) \supseteq [B)$. Note that this is the dual of the partial order induced on $\mathcal{A}$ by the inclusion order of the family of filters $[A)$ for $A \in \mathcal{A}$.

Let $I$ be an ideal of $P$. Let $A$ be the set of minimal elements of $I$. Then $A$ is the (unique) antichain which generates $I$. Similarly, if $F$ is a filter of $P$, the set of maximal elements of $F$ is the unique antichain which generates $F$. This establishes bijective correspondences between antichains and ideals, and between antichains and filters.

Let $\mathcal{F}$ be a family of sets partially ordered by inclusion. The subfamilies of $\mathcal{F}$ can be ordered by defining $\mathcal{G} \leq \mathcal{H}$ iff $\bigcup \mathcal{G} \subseteq \bigcup \mathcal{H}$ for $\mathcal{G}, \mathcal{H} \subseteq \mathcal{F}$. This is an extension of the inclusion order of the subfamilies of $\mathcal{F}$. Since for an antichain $\mathcal{A}$ of $\mathcal{F}$, $\bigcup \mathcal{A} = \bigcup [\mathcal{A}]$, the restriction of this order to the antichains of $\mathcal{F}$ extends the ideal order. Alternatively, the subfamilies of $\mathcal{F}$ can be ordered by defining $\mathcal{G} \leq \mathcal{H}$ iff $\bigcap \mathcal{G} \subseteq \bigcap \mathcal{H}$ for $\mathcal{G}, \mathcal{H} \subseteq \mathcal{F}$. This extends the dual of the inclusion order of the subfamilies of $\mathcal{F}$. Since for an antichain $\mathcal{A}$ of $\mathcal{F}$, $\bigcap \mathcal{A} = \bigcap [\mathcal{A}]$, the restriction of this order to the antichains of $\mathcal{F}$ extends the filter order. Thus, deciding which partial order of the antichains of $\mathcal{F}$ to use often depends on whether we are interested in intersections or in unions of antichains.

Let $P$ and $Q$ be posets. A map $f : P \to Q$ is *order-preserving* iff $x \leq y$ in $P$ implies that $f(x) \leq f(y)$. The composition of order-preserving maps is order-preserving. $P$ and $Q$ are *isomorphic* iff there is an order-preserving bijection $f : P \to Q$ with an order-preserving inverse. The set of order-preserving maps from $P$ to $Q$ is denoted by $Q^P$. The set $Q^P$ is partially ordered by defining $f \leq g$ iff for every $x \in P$, $f(x) \leq g(x)$ in $Q$.

Let $f$ be an order-preserving map from $P$ into the two element chain [2]. Then the set of elements $x$ of $P$ such that $f(x) = \hat{1}$ is a filter of $P$. Conversely if $F$ is a filter of $P$, then the map $f : P \to [2]$ defined by $f(x) = \hat{1}$ iff $x \in F$ is order-preserving. This correspondence between $[2]^P$ and the set of filters of $P$ is order preserving:

**Theorem 2.4.3** *The family of filters of $P$ is isomorphic to $[2]^P$.*

Two other useful constructions on posets are the cartesian product and the disjoint union. The *cartesian product* $P \times Q$ of $P$ and $Q$ has the partial order defined by $\langle x, y \rangle \leq \langle x', y' \rangle$ iff $x \leq x'$ and $y \leq y'$. The *disjoint union* $P + Q$ of $P$ and $Q$ has the partial order defined by $x \leq y$ iff either $x, y \in P$ and $x < y$ in $P$ or $x, y \in Q$ and $x < y$ in $Q$.

## 2.5 Lattices and Semilattices

A *meet-semilattice* is a poset $L$ such that for every $u, v \in L$, the greatest lower bound $u \wedge v$ of $u$ and $v$ exists. The dual of a meet-semilattice is a *join-semilattice*. A *semilattice* is by default a meet-semilattice. If all pairwise meets in a poset exist, then every non-empty (finite) subset has a greatest lower bound. Thus, every meet-semilattice has a least element (its greatest lower bound), and every join-semilattice has a greatest element.

A non-empty intersection-closed family of sets is a meet-semilattice with the meet operation given by intersection; a non-empty union-closed family of sets is a join-semilattice with the join operation given by union.

A *lattice* is a poset which is both a meet- and a join-semilattice. Lattices can be defined algebraically using the properties of the meet and join operations listed in Section 2.4. If $L$ is a meet-semilattice and $U \subseteq L$ has an upper bound, then $U$ has a least upper bound $\bigvee U$, given by the greatest lower bound of the set of upper bounds of $U$. This implies that if $L$ has a greatest member, then $L$ is a lattice. If the semilattice $L$ does not have a greatest member, it suffices to adjoin a new greatest element $\hat{1}$ to $L$, obtaining the lattice $\hat{L} = L \cup \{\hat{1}\}$; this is the *completion* of $L$.

The subset $M$ of the meet-semilattice $L$ is a *meet-subsemilattice* of $L$ iff $M \neq \emptyset$ and for every $u, v \in M$, $u \wedge v$ is in $M$. *Join-subsemilattices* are defined dually. $M$ is a *sublattice* of the lattice $L$ iff $M$ is both a meet-subsemilattice and a join-subsemilattice of $L$.

Let $L$ and $M$ be meet-semilattices. The map $f : L \rightarrow M$ is a *meet-homomorphism* iff for every $u, v \in L$, $f(u \wedge v) = f(u) \wedge f(v)$. If $f$ is a meet-homomorphism, $f$ is said to *preserve* meets. Note that a meet-homomorphism is an order-preserving map. *Join-homomorphisms* are defined dually. If $L$ and $M$ are lattices and $f$ is both meet- and join-preserving, then $f$ is a *lattice homomorphism*. A *semilattice isomorphism* is a meet-preserving bijection. A *lattice isomorphism* is a meet- and join-preserving bijection.

The element $u$ of $L$ is *join-irreducible* iff whenever $v \vee w = u$, either $u = v$ or $u = w$. If $u$ covers two distinct elements $v$ and $w$, then $u = v \vee w$. This implies that $u$ is join-irreducible iff $u$ covers at most one element of $L$. The only join-irreducible of $L$ which covers no other element of $L$ is $\hat{0}$. The *proper join-irreducibles* of $L$ are the join-irreducibles of $L$ other than $\hat{0}$. Note that every element $u$ of $L$ is the join of the join-irreducibles below $u$. *Meet-irreducibles* are defined dually.

The semilattice $L$ is *atomic* iff every proper join-irreducible of $L$ is an atom of $L$. The join-semilattice $L$ is *coatomic* iff every proper meet-irreducible of $L$ is a coatom of $L$.

The *Boolean lattice* $B_n$ generated by $n$ atoms is (isomorphic to) the family of subsets of an $n$-set. The atoms of $B_n$ correspond to the one-element subsets.

The *distributive laws* are the identities

$$(u \wedge v) \vee w = (u \vee w) \wedge (v \vee w), \qquad (u \vee v) \wedge w = (u \wedge w) \vee (v \wedge w).$$

A *distributive* lattice is a lattice satisfying the distributive laws. Note that the dual of a distributive lattice is distributive. Since union distributes over intersection and intersection distributes over union, a family of sets which is both union- and intersection-closed is a distributive lattice; in particular, families of ideals and families of filters of posets are distributive lattices. The converse is part of the Fundamental Theorem of Distributive Lattices (see Grätzer [17] for a proof):

**Theorem 2.5.1** *For every distributive lattice $L$ with $|L| \geq 2$, there is a unique (up to isomorphism) poset $P$ such that $L$ is isomorphic to the family of ideals of $P$.*

The poset $P$ in the theorem is given by the set of proper join-irreducibles of $L$. The isomorphism is obtained by assigning to each element $u \in L$ the ideal of all proper join-irreducibles below $u$. Since the family of ideals of $P$ is the family of filters of $P^*$, Theorems 2.4.3 and 2.5.1 imply that the class of distributive lattices $L$ with $|L| \geq 2$ is equivalent to the class of posets of the form $[2]^P$.

Let $P$ be a poset. Consider the family $\mathcal{A}$ of antichains of $P$. The correspondence between antichains and ideals discussed in Section 2.4 shows that $\mathcal{A}$ with the ideal order is a distributive lattice, where $A \vee B$ is given by the set of maximal elements of $A \cup B$. The correspondence between antichains and filters shows that $\mathcal{A}$ with the filter order is a distributive lattice, where $A \wedge B$ is given by the set of minimal elements of $A \cup B$. These two lattice structures on $\mathcal{A}$ are in general not the same. An important subfamily of the family of antichains of $P$ for which the two lattice structures do coincide is the family of Sperner antichains of $P$:

**Theorem 2.5.2** (Dilworth [12]) *Let $A$ and $B$ be Sperner antichains of the poset $P$. Then for both the filter and the ideal order of antichains, $A \wedge B$ and $A \vee B$ are Sperner antichains of $P$. In either case, $A \wedge B$ is the set of minimal elements of $A \cup B$ and $A \vee B$ is the set of maximal elements of $A \cup B$.*

Let $L$ be a lattice. The *upper covering condition* for $L$ asserts that if $u$ covers $v$ in $L$, then for every $w \in L$, $u \vee w$ covers $v \vee w$ or $u \vee w = v \vee w$. The lattice $L$ is *upper semimodular* iff $L$ satisfies the upper covering condition. The *lower covering condition* for $L$ asserts that if $u$ covers $v$ in $L$, then for every $w \in L$, $u \wedge w$ covers $v \wedge w$ or $u \wedge w = v \wedge w$. The lattice $L$ is *lower semimodular* iff $L$ satisfies the lower covering condition. The lattice $L$ is *modular* iff $L$ is both upper and lower semimodular.

The lattice $L$ is *geometric* iff $L$ is atomic and upper semimodular. A subset $B$ of the atoms of the geometric lattice $L$ is *independent* iff for every proper subset $A$ of $B$, $\bigvee A < \bigvee B$. Independent sets of atoms have many of the properties of independent sets of vectors of vector spaces:

**Theorem 2.5.3** *If $x$ is an element of the geometric lattice $L$, then there is an independent set of atoms $A$ such that $\bigvee A = x$. If $x \leq y$, and $A$ is an independent set of atoms such that $\bigvee A = x$, then there is an independent set of atoms $B$ such that $B \supseteq A$ and $\bigvee B = y$. If $A$ is an independent set of atoms and $a$ is an atom with $a \not\leq \bigvee A$, then $A \cup \{a\}$ is independent.*

## 2.6 Two Results from the Combinatorics of Sets

**Theorem 2.6.1** *Let $\mathcal{F}$ be a tree on an $n$-set. Let $\{U_1, \ldots, U_b\}$ be the family of members $U$ of $\mathcal{F}$ with $|U| \geq 2$. Suppose that for $1 \leq i \leq b$, $U_i$ covers $2 + r_i$ sets in $\mathcal{F}$. Then*

$$|\mathcal{F}| = 2n - 1 - \sum_{i=1}^{b} r_i.$$

**Proof.**    Let $X$ be the domain of $\mathcal{F}$. We can assume that $U_i$ is maximal in $\{U_1, \ldots, U_i\}$. For $i \geq 0$, let $\mathcal{F}_i$ consist of the singletons of $X$ and the sets $U_j$ for $1 \leq j \leq i$. Let $\mathcal{M}_i$ be the set of maximal members of $\mathcal{F}_i$. Since $\mathcal{F}$ is a tree, for $U, V \in \mathcal{F}$, $U$ and $V$ are either comparable or disjoint. It follows that $\mathcal{M}_i$ is a partition of $X$ for each $i$. If $i < j$, then $U_i \not\supseteq U_j$, so if $U_i$ covers $U$ in $\mathcal{F}$, then $U$ is maximal in $\mathcal{F}_{i-1}$. Thus $U_i$ covers $2 + r_i$ members of $\mathcal{M}_{i-1}$, which gives

$$\left| \mathcal{M}_i \right| = \left| \mathcal{M}_{i-1} \right| - 1 - r_i$$

for each $i \geq 1$. This implies that

$$\left| \mathcal{M}_b \right| = \left| \mathcal{M}_0 \right| - b - \sum_{i=1}^{b} r_i.$$

Since $\mathcal{M}_b = \{X\}$ and $\mathcal{M}_0 = \{\{x\} \mid x \in X\}$, this gives

$$1 = n - b - \sum_{i=1}^{b} r_i.$$

Solving for $b$ and using $\mathcal{F} = n + b$ yields

$$\left| \mathcal{F} \right| = 2n - 1 - \sum_{i=1}^{b} r_i,$$

as required.    ∎

**Corollary 2.6.2** *If $\mathcal{F}$ is a tree on an $n$-set, then $\left| \mathcal{F} \right| \leq 2n - 1$.*

**Proof.**    If $U \in \mathcal{F}$ and $\left| U \right| \geq 2$, then $U$ covers at least two members of $\mathcal{F}$. The result follows by Theorem 2.6.1.    ∎

The next result, known as Kleitman's lemma, has many applications. Its proof can be found in Anderson [3].

**Theorem 2.6.3** (Kleitman [25]) *Let $X$ be an $n$-set. If $\mathcal{F}$ and $\mathcal{G}$ are filters in $2^X$, then*

$$\frac{\left| \mathcal{F} \cap \mathcal{G} \right|}{2^n} \geq \frac{\left| \mathcal{F} \right|}{2^n} \frac{\left| \mathcal{G} \right|}{2^n}.$$

This says that the density (in $2^X$) of the intersection of two filters is at least the product of the densities of each.

## 2.7  Notes

**Section 2.1.**    The terminology for induced families of sets and restrictions of families of sets given in Table 2.2 loosely follows that in Lovász [31].

**Section 2.3.**    Families of sets are a fundamental research area of combinatorics. An excellent text on the subject is Anderson [3].

Families of sets are often studied as hypergraphs. Hypergraphs are so-called because they are a natural generalization of graphs. An up-to-date overview of research on hypergraphs can be found in Füredi [15].

The notion of a tree is usually defined for graphs. In graph theory, a *tree* is a connected graph without circuits. An *arborescence* is a directed graph with root $a$ such that for every vertex $x$ of the graph, there is a unique path from $a$ to $x$. Arborescences can be obtained from trees by selecting a root and directing every edge of the graph away from the root. A tree of sets $\mathcal{F}$ can be made into an arborescence by connecting $U$ to $V$ iff $U$ covers $V$ in $\mathcal{F}$. Conversely, an arborescence in which every vertex has either out-degree zero (the *leaves*) or out-degree at least two, can be made into a tree of sets by associating to each vertex $v$ the set of leaves which can be reached from $v$.

Trees and arborescences are well researched structures and have many applications in combinatorics (see Lovász [31]), as well as in computer science (see Cormen et al. [6]) and classification theory.

**Section 2.4.**    Posets are used in all areas of mathematics. The combinatorial properties of posets have attracted a lot of interest, particularly since Rota's 1964 article on the theory of Möbius functions [36]. Yet, in available texts which cover posets, the focus is usually on lattice theory (Crawley and Dilworth [7], Grätzer [17]). For an account of the theory of posets from a combinatorial perspective, including a collection of interesting exercises with solutions and references, see Stanley [41]. Good sources of information on current research in this area are Proceedings [34] and [35].

Altwegg [2] shows that a poset can be defined in terms of its betweenness relation (up to duality). Using Hasse diagrams to represent posets is standard practice. An interesting open problem is to characterize the undirected graphs which arise from Hasse diagrams of posets (see [35]).

The fact that every poset has a linear extension is mentioned as the first of three fundamental results in the theory of ordered sets in Rival [32]. He attributes this result to Szpilrajn [42].

Every subset of a poset can be considered as a sub-poset with the induced order. However, this notion of sub-poset is too general for many purposes. A good alternative is the retract. The subset $Q$ of $P$ with the induced partial order is a *retract* of $P$ iff there is an order-preserving map $f$ from $P$ to $Q$ which is the identity on $Q$. The retract construction is discussed in Rival [33].

Chains and antichains figure prominently in the theory of posets. The height and width of a poset are two of the most important invariants.

Sperner antichains are named after E. Sperner. In 1928 he proved that the width of the family of subsets of an $n$-set is $\binom{n}{\lfloor \frac{n}{2} \rfloor}$. This result, known as Sperner's theorem, is the starting point for the study of the combinatorics of sets (Anderson [3]).

Dilworth's chain decomposition result is the second fundamental result mentioned by Rival in [32]; the third is the fixed point theorem for lattices due to Knaster [26] and Tarski [43].

Posets and order-preserving maps constitute the category of posets. Exponentiation $(Q^P)$ of posets has many of the properties of exponentiation of sets. The operations of exponentiation, cartesian product and disjoint sum can be used to define an arithmetic of posets. This is discussed in Jónsson [24]. See also Birkhoff [4], pp. 66-69.

**Section 2.5.** Lattices are the best studied class of posets. These structures are studied from both an algebraic (Crawley and Dilworth [7], Grätzer [17]) and a combinatorial (Stanley [41]) perspective. The Fundamental Theorem of Distributive Lattices is sometimes used to argue that the study of posets essentially reduces to the study of distributive lattices.

The structure of the lattice of antichains has attracted substantial attention. For a generalization of Theorem 2.5.2 and some surprising duality results, see Greene and Kleitman [18] and Greene [19].

The set of bases of a geometric lattice is often called a matroid. Matroids occur naturally not only as the sets of bases of vector spaces, but also in graph theory and the study of matchings in hypergraphs. Matroids and geometric lattices are discussed in Crapo [8], Grätzer [17] and Welsh [44].

**2.6.** The expression for the size of trees in Theorem 2.6.1 is usually stated as an identity involving the number of leaves and the number of nodes for $n$-ary trees (see Grimaldi [20]).

Kleitman's lemma is considerably strengthened by an inequality due to Ahlswede and Daykin [1] (see Anderson [3]).

## FAMILIES OF SETS WITH WIDTH RESTRICTIONS.

### 3.1 Definitions

Let $\mathcal{F}$ be a family of subsets of the set $X$.

**Definition.** The *center* of $A \subseteq X$ is given by

$$C(A) = \{x \in A \mid \text{if } x \in U \in \mathcal{F}, \text{ then } U \subseteq A \text{ or } A \subseteq U\}.$$

Thus $C(A)$ consists of the elements of $A$ not contained in any member of $\mathcal{F}$ incomparable to $A$. The family $\mathcal{F}$ is *centered* iff for every $U \in \mathcal{F}$, $C(U) \neq \emptyset$.

**Definition.** The family $\mathcal{F}$ is *locally k-wide* iff for every $x \in X$, the width of $\mathcal{F}_{\supseteq \{x\}}$ is at most $k$. Equivalently, $\mathcal{F}$ is locally k-wide iff every antichain $\mathcal{A} \subseteq \mathcal{F}$ with $|\mathcal{A}| = k + 1$ has empty intersection.

Since $X$ and the singletons of $X$ are either comparable to or disjoint from every subset of $X$:

**Observation 3.1.1** *If $\mathcal{F}$ is locally k-wide, then so is the family $\mathcal{F}' = \mathcal{F} \cup \{X\} \cup \{\{x\} \mid x \in X\}$.*

**Observation 3.1.2** $C(\{x\}) = \{x\}$ *and* $C(X) = X$.

**Observation 3.1.3** *If $\mathcal{F}$ is centered, then so is the family $\mathcal{F}' = \mathcal{F} \cup \{X\} \cup \{\{x\} \mid x \in X\}$.*

**Definition.** The family $\mathcal{F}$ is a *pseudotree* on $X$ iff $\mathcal{F}$ is centered, $X \in \mathcal{F}$ and for every $x \in X$, $\{x\} \in \mathcal{F}$. The family $\mathcal{F}$ is a *k-pseudotree* on $X$ iff $\mathcal{F}$ is a locally k-wide pseudotree.

Since $\mathcal{F}$ is locally 1-wide iff for every $U, V \in \mathcal{F}$, $U$ and $V$ are either comparable or disjoint:

**Observation 3.1.4** $\mathcal{F}$ *is a 1-pseudotree iff $\mathcal{F}$ is a tree.*

### 3.2 The Structure of Centered Families of Sets

Let $\mathcal{F}$ be a family of subsets of $X$. Since the center of a set $U$ does not intersect any member of $\mathcal{F}$ incomparable to $U$:

**Observation 3.2.1** *For every $A \subseteq X$, the family $\{U \in \mathcal{F} \mid C(U) = A\}$ is a chain.*

**Observation 3.2.2** *Let $\mathcal{A}$ be an antichain of $\mathcal{F}$ with $|\mathcal{A}| \geq 2$. Then the centers of the members of $\mathcal{A}$ are pairwise disjoint and disjoint from $\bigcap \mathcal{A}$.*

It follows that if $\mathcal{F}$ is centered and $\mathcal{A}$ is an antichain of $\mathcal{F}$ with $|\mathcal{A}| \geq 2$ and $\bigcap \mathcal{A} \neq \emptyset$, then $|\bigcup \mathcal{A}| \geq k + 1$. This gives the following observation:

**Observation 3.2.3** *If $n \geq 2$ and $\mathcal{F}$ is a centered family of subsets of an $n$-set, then $\mathcal{F}$ is an $(n-1)$-pseudotree.*

**Theorem 3.2.4** *Let $U, V \in \mathcal{F}$. If $U \subseteq V$ and $U \cap \mathrm{C}(V) \neq \emptyset$, then $\mathrm{C}(U) \subseteq \mathrm{C}(V)$.*

**Proof.** We show that $U \setminus \mathrm{C}(V) \subseteq U \setminus \mathrm{C}(U)$. Let $x \in U \setminus \mathrm{C}(V)$. Then $x \in V \setminus \mathrm{C}(V)$, so there is a set $W \in \mathcal{F}$ such that $x \in W$ and $W$ is incomparable to $V$. The assumptions on $U$ and $W \cap \mathrm{C}(V) = \emptyset$ imply that $W$ is incomparable to $U$. Hence $x \notin \mathrm{C}(U)$, as required. ∎

**Corollary 3.2.5** *The family $\{\mathrm{C}(U) \mid U \in \mathcal{F}\}$ of centers of members of $\mathcal{F}$ is a forest.*

**Proof.** Let $U, V \in \mathcal{F}$ such that $\mathrm{C}(U) \cap \mathrm{C}(V) \neq \emptyset$. Then $U$ and $V$ are comparable. By Theorem 3.2.4, $\mathrm{C}(U)$ and $\mathrm{C}(V)$ are comparable. ∎

**Theorem 3.2.6** *If $\mathcal{F}$ is a centered family of sets and $U \in \mathcal{F}$, then $\mathrm{w}(\mathcal{F}_{\supseteq \mathrm{C}(U)}) = \mathrm{w}(\mathcal{F}_{\supseteq U})$.*

**Proof.** Let $w = \mathrm{w}(\mathcal{F}_{\supseteq \mathrm{C}(U)})$ and $w' = \mathrm{w}(\mathcal{F}_{\supseteq U})$. Since $\mathcal{F}_{\supseteq \mathrm{C}(U)} \supseteq \mathcal{F}_{\supseteq U} \supseteq \{U\}$, $w \geq w' \geq 1$. To show that $w \leq w'$, let $\{U_1, \ldots, U_w\}$ be an antichain of $\mathcal{F}_{\supseteq \mathrm{C}(U)}$. Since $U_i$ intersects $\mathrm{C}(U)$, the sets $U_i$ and $U$ are comparable for each $i$. If $U_i \supseteq U$ for each $i$, then $\{U_1, \ldots, U_w\} \subseteq \mathcal{F}_{\supseteq U}$ which implies that $w \leq w'$. Suppose that for some $i$, $U_i \subseteq U$. Then by Theorem 3.2.4, $\mathrm{C}(U_i) \subseteq \mathrm{C}(U)$. Observation 3.2.2 implies that $w = 1 \leq w'$. ∎

**Theorem 3.2.7** *Let $\mathcal{F}$ be a centered family of arcs of $\mathbf{Z}_n$ which contains $\mathbf{Z}_n$ and the singletons of $\mathbf{Z}_n$. Then $\mathcal{F}$ is a 2-pseudotree.*

**Proof.** Suppose that $\mathcal{F}$ has three pairwise incomparable arcs $U = [x_\mathrm{l}, x_\mathrm{r}]$, $V = [y_\mathrm{l}, y_\mathrm{r}]$ and $W = [z_\mathrm{l}, z_\mathrm{r}]$ such that $U \cap V \cap W \neq \emptyset$. Let $x \in U \cap V \cap W$. We can assume that in clockwise cyclic order

$$x_\mathrm{l} < y_\mathrm{l} < z_\mathrm{l} \leq x \leq x_\mathrm{r} < y_\mathrm{r} < z_\mathrm{r}.$$

This implies that $V \subseteq U \cup W$, so the center of $V$ is empty, contradicting the assumption that $\mathcal{F}$ is centered. ∎

### 3.3 Basic Extremal Results

**Theorem 3.3.1** *If $\mathcal{F}$ is a centered family of subsets of an $n$-set, then $|\mathcal{F}| \leq \frac{(n+1)n}{2}$. This bound is best possible.*

**Proof.**    To show that the bound is attained:

**Example 3.3.2** Let

$$\mathcal{F} = \{[1, i-1] \cup \{j\} \mid 1 \leq i \leq n, \ i \leq j \leq n\},$$

so that $\mathcal{F}$ consists of the initial segments of $[n]$ with a singleton adjoined. Then $\mathcal{F}$ is a centered family of subsets of $[n]$ where the center of the member $[1, i-1] \cup \{j\}$ of $\mathcal{F}$ is $\{j\}$. The cardinality of $\mathcal{F}$ is $\frac{(n+1)n}{2}$.

**Proof of the bound.**    Let $1 \leq r \leq n$ and consider the family $[\mathcal{F}]_r$ of $r$-element members of $\mathcal{F}$. Let $U \in [\mathcal{F}]_r$. Since the members of $[\mathcal{F}]_r$ are pairwise incomparable, if $V \in [\mathcal{F}]_r$ and $V \neq U$, then $C(V) \cap U = \emptyset$. This implies that $|U| + \left|[\mathcal{F}]_r \setminus \{U\}\right| \leq n$, hence

$$\left|[\mathcal{F}]_r\right| \leq n + 1 - r.$$

The bound is obtained by summing this inequality for $r = 1, \ldots, n$. ∎

**Theorem 3.3.3** *If $\mathcal{F}$ is a centered family of segments of $[n]$, then $|\mathcal{F}| \leq 3n - 3$. This bound is best possible.*

**Proof.**    To show that the bound is attained:

**Example 3.3.4** Let

$$\mathcal{F} = \{\{i\} \mid i \in [n]\} \cup \{[1, i] \mid i \in [n]\} \cup \{[i, n] \mid i \in [n]\},$$

so that $\mathcal{F}$ consists of the singletons of $[n]$ and the segments of $[n]$ containing either 1 or $n$. In any family of segments of $[n]$, $1 \in C([1, i])$ and $n \in C([j, n])$. Thus $\mathcal{F}$ is centered. The cardinality of $\mathcal{F}$ is $3n - 3$, as desired.

**Proof of the bound.**    By Observation 3.1.3, we can assume that $\mathcal{F}$ contains $[n]$ and the singletons of $[n]$. By cyclic closure, families of segments of $[n]$ can be considered as families of arcs. By Theorem 3.2.7, any centered family of arcs of $\mathbf{Z}_n$ which contains $\mathbf{Z}_n$ and the singletons of $\mathbf{Z}_n$ is a 2-pseudotree. By Theorem 3.4.1 below, $3n - 3$ is the bound on the size of 2-pseudotrees on an $n$-set, so the result follows. ∎

**Theorem 3.3.5** *If $\mathcal{F}$ is a locally $k$-wide family of subsets of an $n$-set $X$ and $r \geq 2$, then $\left|[\mathcal{F}]_r\right| \leq \frac{nk}{r}$.*

**Proof.**    Every element of $X$ is contained in at most $k$ members of $[\mathcal{F}]_r$. Hence, $kn = k\big|X\big| \geq \big|[\mathcal{F}]_r\big|r$. ∎

From this we can deduce a bound on the maximum size of a locally k-wide family of sets.

**Theorem 3.3.6** *If $\mathcal{F}$ is a locally $k$-wide family of subsets of an $n$-set, then $\big|\mathcal{F}\big| \leq 2 + n + kn\ln(n)$.*

**Proof.**    For $n \geq 2$,

$$
\begin{aligned}
\big|\mathcal{F}\big| &= \big|[\mathcal{F}]_0\big| + \big|[\mathcal{F}]_1\big| + \ldots + \big|[\mathcal{F}]_n\big| \\
&\leq 1 + n + \frac{kn}{2} + \ldots + \frac{kn}{n-1} + 1 \\
&\leq 2 + n + kn\ln(n-1).
\end{aligned}
$$

∎

A bound on the size of locally $k$-wide families of sets linear in the size of the domain will be obtained in Section 3.8.

## 3.4    The Size of $k$-pseudotrees

**Theorem 3.4.1** *Let $n > k \geq 1$. If $\mathcal{F}$ is a $k$-pseudotree on an $n$-set $X$, then $\big|\mathcal{F}\big| \leq (k+1)n - \frac{(k+1)k}{2}$. This bound is best possible.*

By Observation 3.2.3, if $n \geq 2$, a pseudotree on an $n$-set is an $(n-1)$-pseudotree. Therefore the assumption that $n > k$ does not restrict the generality of the theorem.

**Proof.**    To show that the bound is attained:

**Example 3.4.2** Let

$$
\mathcal{F} = \big\{\{i\} \mid i \leq n\big\} \cup \big\{[1,i] \cup \{j\} \mid 1 \leq i < j \leq \min(i+k,n)\big\},
$$

so that $\mathcal{F}$ consists of the singletons of $[n]$ and the initial segments of $[n]$ with one of the next $k$ singletons adjoined.

The cardinality of $\mathcal{F}$ is $n + (n-1) + \ldots + (n-k) = (k+1)n - \frac{(k+1)k}{2}$.

To see that $\mathcal{F}$ is centered, let $U \in \mathcal{F}$. If $V \in \mathcal{F}$ and $V$ contains the greatest integer in $U$, then either $V \subset U$ or $U \subset V$. Therefore $j \in \mathrm{C}(U)$.

To show that $\mathcal{F}$ is locally $k$-wide, suppose that $\mathcal{U} = \{U_1, \ldots, U_{k+1}\}$ is an antichain of $\mathcal{F}$. For each $i$, let $x_i$ be the greatest integer in $U_i$. Since $x_i \in \mathrm{C}(U_i)$ for each $i$, the $x_i$ are distinct. Without loss of generality, assume that $x_1 < x_2 < \ldots < x_{k+1}$. Then $x_{k+1} - x_1 \geq k$, so if $U_{k+1} \neq \{x_{k+1}\}$ then $U_{k+1} \supseteq U_1$, contrary to assumption. Therefore $U_{k+1} = \{x_{k+1}\}$ which implies that $\bigcap \mathcal{U} = \emptyset$, as desired. Note that this argument shows that the family of non-singleton members of $\mathcal{F}$ has width $k$.

**Proof of the bound.**     By associating some of the sets in $\mathcal{F}$ with their centers in a one-to-one way, and then associating the remaining sets in $\mathcal{F}$ with elements of $X$ in a $k$-to-one way, we will first establish an upper bound of $(k+2)n-1$. The bound will then be reduced to $(k+1)n-1$ by an analysis of these associations. To further reduce this to the bound given in the theorem, a notion of *deficiency* will be defined for arbitrary members of $k$-pseudotrees. The deficiency of $X$ will be the minimum value of $((k+1)n-1) - \lfloor \mathcal{F} \rfloor$.

Let $\mathcal{C}$ be the set of centers of members of $\mathcal{F}$. By Corollary 3.2.5 and since $\mathcal{F}$ contains $X$ and the singletons of $X$, $\mathcal{C}$ is a tree. By Observation 3.2.1, for every $C \in \mathcal{C}$, the family $\{U \in \mathcal{F} \mid \mathrm{C}(U) = C\}$ is a chain. Associate the minimal member of this chain with $C$. Since the maximum size of a tree on an $n$-set is $2n-1$ (Theorem 2.6.1), this association (which is one-to-one) accounts for $|\mathcal{C}| \leq 2n-1$ members of $\mathcal{F}$.

Let $\mathcal{F}'$ consist of the members of $\mathcal{F}$ not associated with their center. For every $U \in \mathcal{F}'$, we will select an element $p(U) \in U \setminus \mathrm{C}(U)$. This will be done in such a way that if $p(U) = p(V)$ and $U \neq V$, then $U$ and $V$ are incomparable.

**Lemma 3.4.3** *Let $U \in \mathcal{F}'$. Then there is an element $x \in U \setminus \mathrm{C}(U)$ such that for every $W \in \mathcal{F}$ with $x \in W \subset U$, $x \in \mathrm{C}(W)$.*

**Proof.**     Let $V$ be the maximal member of the chain of sets

$$\{W \in \mathcal{F} \mid W \subset U \text{ and } \mathrm{C}(W) = \mathrm{C}(U)\}.$$

Let $V_1$ be a maximal member of

$$\{W \in \mathcal{F} \mid W \subset U \text{ and } W \cap (U \setminus V) \neq \emptyset\}.$$

Since the singletons of $X$ are in $\mathcal{F}$, such sets exist and $V_1$ is well-defined. If $V \subset V_1$, then by Theorem 3.2.4, $\mathrm{C}(V) \subseteq \mathrm{C}(V_1) \subseteq \mathrm{C}(U)$, contradicting maximality of $V$. Therefore $V_1$ and $V$ are incomparable and $\mathrm{C}(V_1) \subseteq U \setminus V$. Let $V_1 \supset V_2 \supset \ldots \supset V_r$ be a (strictly) decreasing chain of members of $\mathcal{F}$ such that

(i) $V_r = \{x\}$ for some $x \in X$,

(ii) $V_{j+1}$ is a maximal member of

$$\{W \in \mathcal{F} \mid W \subset V_j \text{ and } W \cap \mathrm{C}(V_j) \neq \emptyset\}.$$

Such a chain exists. To show that $x$ is as desired, let $W$ be a member of $\mathcal{F}$ with $x \in W \subset U$. By Theorem 3.2.4, the centers of the $V_j$ form a decreasing chain. By the maximality condition on $V_1$, $W \not\supseteq V_1$. Since $W$ intersects the center of $V_j$, $W$ is comparable to $V_j$ for each $j$. In particular, $W \subseteq V_1$. Let $j$ be the greatest index such that $W \subseteq V_j$. Since $V_r = \{x\}$, if $j = r$, then $W = V_j$. If $j < r$, then the maximality condition on $V_{j+1}$ gives $W = V_j$. Hence $x \in \mathrm{C}(W)$. ∎

For each $U \in \mathcal{F}'$, let $p(U)$ be an element of $U \setminus \mathrm{C}(U)$ which satisfies the conclusion of Lemma 3.4.3. If $p(U) = p(V)$ and $U \neq V$, then $U$ and $V$ are incomparable. For every $x \in X$, define

$$s(x) = \left| \{ U \in \mathcal{F}' \mid p(U) = x \} \right|$$

and

$$l(x) = \mathrm{w}(\mathcal{F}_{\supseteq \{x\}}).$$

Since $\mathcal{F}$ is locally $k$-wide, $s(x) \leq l(x) \leq k$. Therefore $\left| \mathcal{F}' \right| \leq kn$. Since $\left| \mathcal{F} \setminus \mathcal{F}' \right| = \left| \mathcal{C} \right|$,

$$
\begin{aligned}
\left| \mathcal{F} \right| &= \left| \mathcal{C} \right| + \left| \mathcal{F}' \right| \\
&= \left| \mathcal{C} \right| + \sum_{x \in X} s(x) \\
&\leq (k+2)n - 1.
\end{aligned}
$$

Let $X' = \{ x \in X \mid s(x) = l(x) \}$.

**Lemma 3.4.4** $\left| \mathcal{C} \right| \leq 2n - 1 - \left| X' \right|$.

**Proof.**     For $x \in X'$, let $C_x$ be the minimal member of

$$\{ C \in \mathcal{C} \mid x \in C \text{ and for every } U \in \mathcal{F}' \text{ with } p(U) = x, \, \mathrm{C}(U) \subseteq C \}.$$

Since $\mathrm{C}(X) = X$, such sets exist and $C_x$ is well-defined. Let $\mathcal{C}' = \{ C_x \mid x \in X' \}$. For $C \in \mathcal{C}'$, let $r_C = \left| \{ x \in X' \mid C_x = C \} \right|$. We will show that if $C \in \mathcal{C}'$, then $C$ covers at least $2 + r_C$ members of $\mathcal{C}$. Since $\sum_{C \in \mathcal{C}'} r_C = \left| X' \right|$ and by the bound on the size of trees in Theorem 2.6.1, it will follow that $\left| \mathcal{C} \right| \leq 2n - 1 - \left| X' \right|$.

Let $x \in X'$. Let $\mathcal{U} = \{ U_1, \ldots, U_{l(x)} \}$ be the family of members $U$ of $\mathcal{F}'$ with $p(U) = x$. Let $C_i$ be the center of $U_i$. Then $x \notin C_i$ and since $\mathcal{U}$ is an antichain, the $C_i$ are disjoint. Observe that $l(x) \geq 2$. Proof: If $l(x) = 1$, then the sets $U \in \mathcal{F}$ with $x \in U$ form a chain. This implies that if $x \in U \in \mathcal{F}$, then $x \in \mathrm{C}(U)$ so that $x \neq p(U)$. Hence $s(x) = 0$, contradicting $s(x) = l(x)$.

**Lemma 3.4.4.1** *If $C \in \mathcal{C}$ contains $x$ and $C \cap C_i \neq \emptyset$ for some $i$, then $C \supseteq C_x$.*

**Proof.**     Without loss of generality, assume that $C$ intersects $C_1$. Let $W$ be a set in $\mathcal{F}$ with center $C$. Since $x \in C$, the $U_i$ are comparable to $W$. Since $W$ intersects $C_1$ and by Theorem 3.2.4, if $W \subseteq U_1$, then $C \subseteq C_1$ which contradicts $x \in C$. Therefore $W \supseteq U_1$. For $i \geq 2$, $U_i$ and $U_1$ are incomparable, so $W \supset U_i$ for each $i$. Since $x \in U_i \cap C$, we have $C_i \subseteq C$ for each $i$. By minimality of $C_x$ it follows that $C_x \subseteq C$, as required. ∎

**Lemma 3.4.4.2** *If $C \in \mathcal{C}$ and $C \supseteq C_i$ for each $i$, then $C \supseteq C_x$.*

**Proof.** It suffices to show that $x \in C$. Let $W$ be a member of $\mathcal{F}$ with center $C$. Since the $C_i$ are disjoint, $C \supset C_i$ and therefore $W \supset U_i$ for each $i$. Suppose that $x \notin C$. Let $V$ be a set in $\mathcal{F}$ incomparable to $W$ such that $x \in V$. Since $C$ and $V$ are disjoint, $V$ and $U_i$ are incomparable for each $i$. It follows that $\{V\} \cup \mathcal{U}$ is an antichain with $x \in \bigcap(\{V\} \cup \mathcal{U})$, contradicting $l(x) = |\mathcal{U}|$. ∎

**Lemma 3.4.4.3** *Suppose that $y \in X'$ with $y \neq x$ and $C_y = C_x$. If $C$ is a member of $\mathcal{C}$ with $x \in C$ and $y \in C$, then $C \supseteq C_x$.*

**Proof.** Let $\mathcal{V} = \{V_1, \ldots, V_{l(y)}\}$ be the family of members $V$ of $\mathcal{F}$ such that $p(V) = y$. Let $W$ be a member of $\mathcal{F}$ with center $C$.

Suppose that $y \notin U_i$ for some $i$. Then $W \nsubseteq U_i$. Since $U_i$ intersects $C$, we have $U_i \subset W$ and $C_i \subseteq C$, so that by Lemma 3.4.4.1, $C \supseteq C_x$. By symmetry, if $x \notin V_i$ for some $i$, then $C \supseteq C_y = C_x$.

Suppose that $y \in \bigcap \mathcal{U}$ and $x \in \bigcap \mathcal{V}$. Since $l(y)$ is the maximum size of an antichain of $\mathcal{F}_{\supseteq \{y\}}$, $U_1$ is comparable to at least one of the $V_i$. We can assume that $U_1$ is comparable to $V_1$. If $U_1 \subset V_1$, then by the property expressed in Lemma 3.4.3, $y \in C_1$, so that $C$ intersects $C_1$ and by Lemma 3.4.4.1, $C \supseteq C_x$. The case $V_1 \subset U_1$ is symmetric. This completes the proof. ∎

Let $C \in \mathcal{C}'$ and let $\{x_1, \ldots, x_r\}$ be the set of elements $x \in X'$ such that $C_x = C$ (so that $r = r_C$). Suppose that $C$ covers $C'$ in $\mathcal{C}$. Then $C'$ contains at most one of the $x_i$ (Lemma 3.4.4.3). If $x_i \in C'$, then $C'$ does not include the center of any set $U$ with $p(U) = x_i$ (Lemma 3.4.4.1). Consider $x_1$. If $\mathcal{U} = \{U_1, \ldots, U_{l(x_1)}\}$ is the family of $U \in \mathcal{F}$ with $p(U) = x_1$, then there is an $i$ such that $C'$ does not include $\mathrm{C}(U_i)$ (Lemma 3.4.4.2). Since $\mathrm{C}(U_i) \subseteq C$ for each $i$ and $l(x_1) \geq 2$, this implies that $C$ covers at least $r + 2$ members of $\mathcal{C}$, as required. ∎

Lemma 3.4.4 yields

$$\begin{aligned}
|\mathcal{F}| &= |\mathcal{C}| + \sum_{x \in X} s(x) \\
&\leq 2n - 1 - |X'| + \sum_{x \in X} s(x).
\end{aligned}$$

For $x \in X'$, $s(x) = l(x)$ and for $x \in X \setminus X'$, $s(x) \leq l(x) - 1$. This gives

$$\begin{aligned}
|\mathcal{F}| &\leq 2n - 1 + \sum_{x \in X} (l(x) - 1) \\
&= n - 1 + \sum_{x \in X} l(x) \\
&= (k+1)n - 1 - \sum_{x \in X} (k - l(x)).
\end{aligned}$$

Let $d(x) = k - l(x)$. Then

$$|\mathcal{F}| \leq (k+1)n - 1 - \sum_{x \in X} d(x).$$

The number $d(x)$ is called the *deficiency* (in $\mathcal{F}$) of $x$. For $U \in \mathcal{F}$, the *deficiency* (in $\mathcal{F}$) of $U$ is defined by

$$d(U) = \sum_{x \in \mathrm{C}(U)} d(x).$$

Let $d(m, l)$ be the minimum value of the deficiency in $\mathcal{G}$ of $U$ for arbitrary $k$-pseudotrees $\mathcal{G}$ and members $U$ of $\mathcal{G}$ such that $l = \mathrm{w}(\mathcal{G}_{\supseteq U})$ and $m$ is the cardinality of the center (in $\mathcal{G}$) of $U$. Define

$$b(m, l) = \sum_{i=0}^{\min(m-1, k-l)} (k+1-l-i).$$

**Lemma 3.4.5** $d(m, l) + 1 \geq b(m, l)$.

In particular, this shows that $d(n, 1) \geq \frac{k(k+1)}{2} - 1$. Since $d(n, 1)$ is the deficiency of $X$ (the center of which is $X$), we obtain $|\mathcal{F}| \leq (k+1)n - 1 - d(n, 1) \leq (k+1)n - \frac{k(k+1)}{2}$, thus completing the proof of the theorem.

**Proof.** If $l = k$, then the lemma asserts that $d(m, k) + 1 \geq (k + 1 - k) = 1$. Since $d(x) \geq 0$ for each $x$, the inequality $d(m, k) + 1 \geq b(m, k)$ holds.

Suppose that $l < k$. The remainder of the proof proceeds by induction on $m$. Let $U$ be a set in the $k$-pseudotree $\mathcal{G}$ with center $C$ such that $\mathrm{w}(\mathcal{G}_{\supseteq U}) = l$ and $|C| = m$. We compute $l(x)$, $s(x)$, $d(x)$ and $\mathrm{C}(V)$ relative to $\mathcal{G}$. Let $\mathcal{C}$ be the tree of centers of members of $\mathcal{G}$.

Let $m = 1$. Then $C = \{x\}$ for some element $x$ and $d(U) = d(x) = k - l(x)$. To show that $d(U) \geq k - l$, observe that by Theorem 3.2.6, $l(x) = l$.

Let $m > 1$. Since by Theorem 3.2.6, the width of $\mathcal{G}_{\supseteq V}$ depends only on $\mathrm{C}(V)$, we can assume that $U$ is minimal among the sets $V$ of $\mathcal{G}$ with $\mathrm{C}(V) = C$.

Let $\mathcal{U}$ consist of the sets $V$ of $\mathcal{G}$ such that for some $C' \in \mathcal{C}$ covered by $C$, $V$ is the maximal member of the chain $\{W \in \mathcal{G} \mid \mathrm{C}(W) = C'\}$. Since $|C| \geq 2$ and the singletons are in $\mathcal{C}$, $C$ covers at least two members of $\mathcal{C}$. Therefore $|\mathcal{U}| \geq 2$. If $C$ covers $C'$ in $\mathcal{C}$ and $C'$ is the center of $V \in \mathcal{G}$, then $V \subset U$. This implies that every member of $\mathcal{U}$ is included in $U$.

We show that $\mathcal{U}$ has at least two maximal members. Let $U_1$ be maximal in $\mathcal{U}$. Let $x \in C \setminus \mathrm{C}(U_1)$. Let $V$ be maximal among the members $W$ of $\mathcal{G}$ such that $W$ is incomparable to $U_1$ and $x \in W$. Since $V$ is disjoint from $\mathrm{C}(U_1) \subset C$ but intersects $C$, it follows that $V \subset U$ and $\mathrm{C}(V) \subset C$. Let $C'$ be the member of $\mathcal{C}$ which contains $\mathrm{C}(V)$ and is covered by $C$. Let $U_2$ be the member of $\mathcal{U}$ with center $C'$. Then $U_2 \supseteq V$ so that $U_2 \not\subseteq U_1$. Since $U_1$ is maximal in $\mathcal{U}$, $U_2$ and $U_1$ are incomparable. By the maximality condition on $V$, $U_2 = V$ and $U_2$ is maximal in $\mathcal{U}$, as desired.

Write $\mathcal{U} = \{U_1, U_2, \ldots, U_t\}$ where $U_i$ is maximal in $\mathcal{U}_i = \{U_i, \ldots, U_t\}$ for each $i$, and $U_1$ and $U_2$ are maximal in $\mathcal{U}$. Let $c_i = \left| \mathrm{C}(U_i) \right|$. Let $s_i = c_1 + \ldots + c_i$. Note that $s_t = m$. Let $r$ be the least index $i$ such that $s_i \geq k$ or $i = t$.

**Lemma 3.4.5.1** *For* $2 < i \leq r$, $\mathrm{w}(\mathcal{G}_{\supseteq U_i}) \leq \max(s_{i-1}, l)$. *For* $i = 1, 2$, $\mathrm{w}(\mathcal{G}_{\supseteq U_i}) \leq l$.

**Proof.** Let $w = \mathrm{w}(\mathcal{G}_{\supseteq U_i})$. Let $\mathcal{V} = \{V_1, \ldots, V_w\}$ be an antichain of $\mathcal{G}_{\supseteq U_i}$. Since $V_j$ intersects $C$, $V_j$ is comparable to $U$ for each $j$. Suppose that for some $j$, $V_j \supseteq U$. Then $\mathcal{V} \subseteq \mathcal{G}_{\supseteq U}$ and therefore $w \leq l$. Suppose that $V_j \subset U$ for each $j$. Then $\mathrm{C}(V_j) \subset C$ for each $j$. Hence for every $j$ there is a $U_{k(j)}$ such that $V_j \subseteq U_{k(j)}$ and $\mathrm{C}(V_j) \subseteq \mathrm{C}(U_{k(j)})$. We have $U_i \subseteq V_j$ for each $j$, so by maximality of $U_1$ and $U_2$ in $\mathcal{U}$, if $i = 1$ or $i = 2$, then $w = 1$. Suppose that $i > 2$. Since $U_i$ is maximal in $\mathcal{U}_i$, either $w = 1$ or the $V_j$ are not included in any member of $\mathcal{U}_i$. If $w = 1$ we are done. If $w > 1$, then $k(j) < i$ for each $j$. This implies that

$$\bigcup_{j=1}^{w} \mathrm{C}(V_j) \subseteq \bigcup_{j=1}^{w} \mathrm{C}(U_{k(j)}) \subseteq \bigcup_{j=1}^{i-1} \mathrm{C}(U_j).$$

Since $\left| \bigcup_{j=1}^{i-1} \mathrm{C}(U_j) \right| = s_{i-1}$ and the centers of the $V_j$ are disjoint, $w \leq s_{i-1}$, as required. ∎

Since

$$\begin{aligned} d(U) &= d(U_1) + \ldots + d(U_t) \\ &\geq d(c_1, \mathrm{w}(\mathcal{G}_{\supseteq U_1})) + \ldots + d(c_r, \mathrm{w}(\mathcal{G}_{\supseteq U_r})) \end{aligned}$$

and $c_i < m$ for each $i$, we can apply the induction hypothesis to bound $d(U)$ from below. For $3 \leq i \leq r$, let $l_i = \max(s_{i-1}, l)$. For $i = 1$ and $i = 2$, let $l_i = l$. Let $r_i = \min(c_i - 1, k - l_i)$. By Lemma 3.4.5.1, $\mathrm{w}(\mathcal{G}_{\supseteq U_i}) \leq l_i$. Since the inductively obtained lower bound $b(c_i, l')$ for $d(c_i, l') + 1$ is decreasing in $l'$, it follows that

$$\begin{aligned} d(U) + r &\geq b(c_1, l_1) &&+ b(c_2, l_2) + \ldots + b(c_r, l_r) \\ &= (k + 1 - l_1) &&+ \sum_{i=1}^{r_1}(k + 1 - l_1 - i) \\ &\quad + 1 + (k - l_2) &&+ \sum_{i=1}^{r_2}(k + 1 - l_2 - i) \\ &\quad \vdots && \\ &\quad + 1 + (k - l_r) &&+ \sum_{i=1}^{r_r}(k + 1 - l_r - i). \end{aligned}$$

Since the sum of the leading 1's is $r - 1$,

$$
\begin{aligned}
d(U) + 1 \;\geq\;\; & (k + 1 - l_1) & + \;\; & \textstyle\sum_{i=1}^{r_1}(k + 1 - l_1 - i) \\
& + (k - l_2) & + \;\; & \textstyle\sum_{i=1}^{r_2}(k + 1 - l_2 - i) \\
& \;\;\vdots & & \\
& + (k - l_r) & + \;\; & \textstyle\sum_{i=1}^{r_r}(k + 1 - l_r - i).
\end{aligned}
$$

By using $l_1 = l_2 = l$ and reducing some of the terms, we obtain

$$
\begin{aligned}
d(U) + 1 \;\geq\;\; & (k + 1 - l) & + \;\; & (k + 1 - l_1 - 1) + \ldots + (k + 1 - l - r_1) \\
& + (k - l) & + \;\; & (k - l - 1) + \ldots + (k - l - r_2) \\
& \;\;\vdots & & \\
& + (k - l_r) & + \;\; & (k - l_r - 1) + \ldots + (k - l_r - r_r).
\end{aligned}
$$

Let $q_2 = \min(r_1 + r_2 + 1, k - l)$ and for $i \geq 3$, let $q_i = r_i$. We have $(k + 1 - l - r_1) \leq (k - l) + 1$, so that we can combine the first two lines of the above inequality to obtain

$$
\begin{aligned}
d(U) + 1 \;\; \geq \;\; & \\
(2) \qquad & (k + 1 - l) + (k + 1 - l - 1) + \ldots + (k + 1 - l - q_2) \\
(3) \qquad & + (k - l_3) + (k - l_3 - 1) + \ldots + (k - l_3 - q_3) \\
\vdots \qquad & \\
(r) \qquad & + (k - l_r) + (k - l_r - 1) + \ldots + (k - l_r - q_r).
\end{aligned}
$$

Note that each line is a sum of successive non-negative integers. Let $T_{i,j}$ be the $(j + 1)$'th term of line $(i)$, and let $T$ be the sum of the $T_{i,j}$. To complete the proof, we show that $T \geq b(m, l)$.

If $q_2 = k - l$ or $q_2 = m - 1$, then

$$
T \geq \sum_{j=0}^{q_2} T_{2,j} \geq \sum_{j=0}^{\min(m-1, k-l)} (k + 1 - l - j) = b(m, l).
$$

Suppose that neither $q_2 = k - l$ nor $q_2 = m - 1$. Then

   (i) $r_1 = c_1 - 1 < k - l$,

   (ii) $r_2 = c_2 - 1 < k - l$,

   (iii) $q_2 = c_1 + c_2 - 1 = s_2 - 1$,

   (iv) $s_2 < m$,

   (v) $s_2 < k - l + 1 \leq k$.

Inequalities (iv) and (v) imply that $r \geq 3$. Let $r'$ be the least index $i$ such that $q_i = k - l_i$ or $i = r$. By (i) and (ii), $r' \geq 3$. Let $2 \leq i < r'$. We show that

$T_{i,q_i} \leq T_{i+1,0} + 1$. If $i = 2$, we have

$$
\begin{aligned}
T_{2,q_i} &= (k + 1 - l) - q_2 \\
&= k - (l + s_2 - 1) + 1 \\
&\leq k - \max(s_2, l) + 1 \\
&= k - l_3 + 1 \\
&= T_{3,0} + 1.
\end{aligned}
$$

Suppose that $i > 2$. Since $q_i = r_i \neq k - l_i$, we have $q_i = c_i - 1$. Using $l_{i+1} = \max(s_i, l) \leq \max(s_{i-1}, l) + c_i = l_i + c_i$, we obtain

$$
\begin{aligned}
T_{i,q_i} &= k - l_i - r_i \\
&= k - l_i - (c_i - 1) \\
&= k - (l_i + c_i) + 1 \\
&\leq k - l_{i+1} + 1 \\
&= T_{i+1,0} + 1.
\end{aligned}
$$

Observe that

$$
\sum_{i=2}^{r'-1} (q_i + 1) = s_{r'-1}.
$$

Let $q = s_{r'-1} + (q_{r'} + 1)$. We can now bound

$$
\begin{aligned}
T &= \sum_{i=2}^{r} \sum_{j=0}^{q_i} T_{i,j} \\
&\geq \sum_{i=2}^{r'} \sum_{j=0}^{q_i} T_{i,j} \\
&\geq \sum_{j=0}^{\min(q-1, k-l)} (k + 1 - l - j).
\end{aligned}
$$

It remains to show that $\min(q - 1, k - l) \geq \min(m - 1, k - l)$. Either $q_{r'} = k - l_{r'}$ or $r' = r$ and $q_{r'} = q_r = c_r - 1$. Suppose that $q_{r'} = k - l_{r'}$. Then either $l_{r'} = l$ or $l_{r'} = s_{r'-1}$. If $l_{r'} = l$, then $q - 1 \geq q_{r'} = k - l$ and the result follows. If $l_{r'} = s_{r'-1}$, then $q = s_{r'-1} + (k - s_{r'-1} + 1) = k + 1$ and we are done. Suppose that $r' = r$ and $q_{r'} = c_r - 1$. Then $q = s_r$. By definition of $r$, $s_r \geq \min(m, k)$. This gives $\min(q - 1, k - l) \geq \min(m - 1, k - l)$, as desired. ∎

The proof of the theorem is complete. ∎

### 3.5  Intersecting Antichains of Arcs

**Theorem 3.5.1** *If $\mathcal{A}$ is an antichain of arcs of $\mathbf{Z}_n$ such that $|\mathcal{A}| = k$ and $\bigcap \mathcal{A} \neq \emptyset$, then every $A \in \mathcal{A}$ has at least $k$ elements.*

**Proof.** Suppose that $\mathcal{A} = \{A\}$. Then $\bigcap \mathcal{A} \neq \emptyset$ implies that $A$ contains at least one element. Suppose that $|\mathcal{A}| \geq 2$. Let $A \in \mathcal{A}$ and $x \in \bigcap \mathcal{A}$. Then $A \neq \mathbf{Z}_n$, so that $A = [a_l, a_r]$ for some $a_l, a_r \in \mathbf{Z}_n$. By rotation we can assume that $a_l = 0 \leq x \leq a_r$. If $B \in \mathcal{A}$ and $B \neq A$, then $B$ has either a left endpoint $b_l$ with $a_l < b_l \leq x$ or a right endpoint $b_r$ with $x \leq b_r < a_r$, but not both. Since $\mathcal{A}$ is an antichain, no two members of $\mathcal{A}$ have the same left or the same right endpoints. Therefore

$$\big|[a_l, x] \setminus \{a_l\}\big| + \big|[x, a_r] \setminus \{a_r\}\big| \geq k - 1,$$

which implies that $|A| \geq k$. ∎

**Corollary 3.5.2** *If $\mathcal{F}$ is a maximal locally $k$-wide family of arcs, then $\mathcal{F}$ contains the $i$-element arcs for $0 \leq i \leq k$.*

### 3.6 The Size of Locally $k$-wide Families of Arcs and Segments

**Theorem 3.6.1** *If $n \geq 1$ and $\mathcal{F}$ is a locally $k$-wide family of arcs of $\mathbf{Z}_n$, then $|\mathcal{F}| \leq 2kn - k + 1$.*

Let $\mathrm{A}(k, n)$ be the maximum cardinality of a locally $k$-wide family of arcs of $\mathbf{Z}_n$.

The bound on the size of locally $k$-wide families of arcs is not optimal. Since for $n \leq k + 1$, the family of arcs of $\mathbf{Z}_n$ is locally $k$-wide, we have $\mathrm{A}(k, k) = k(k - 1) + 2$ and $\mathrm{A}(k, k + 1) = (k + 1)k + 2$. For $n = k + 2$ there are at most $k$ many $(k + 1)$-element arcs in $\mathcal{F}$, so that $\mathrm{A}(k, k + 2) = k(k + 2) + k + 2$. This shows that $\mathrm{A}(k, n) = 2kn - k^2 + k + 2$ for $n = k, k + 1, k + 2$. It is conjectured that $\mathrm{A}(k, n) = 2kn - k^2 - k + 2$ for all $n \geq k$.

**Example 3.6.2** Let

$$\mathcal{F} = \{[i, j] \mid \big|[i, j]\big| \leq k \text{ or } 0 \leq i \leq k - 1\} \cup \{\mathbf{Z}_n\}.$$

so that $\mathcal{F}$ consists of the $i$-element arcs for $i \leq k$, $\mathbf{Z}_n$ and the arcs with left endpoint $j$, where $0 \leq j \leq k - 1$.

By Theorem 3.5.1 and its corollary, to show that $\mathcal{F}$ is locally $k$-wide it suffices to consider the family $[\mathcal{F}]_{>k}$ of $i$-element arcs in $\mathcal{F}$ with $i \geq k + 1$. Since $[\mathcal{F}]_{>k}$ is included in the disjoint union of the $k$ chains $\{[i, j] \mid j \in \mathbf{Z}_n\}$ for $0 \leq i \leq k - 1$, $\mathrm{w}([\mathcal{F}]_{>k}) \leq k$, which implies that $[\mathcal{F}]_{>k}$ is locally $k$-wide, as required.

The number of $i$-element arcs of $\mathbf{Z}_n$ with $i \leq k$ is $1 + kn$ (including the empty arc). There are $n - k - 1$ many $j$-element arcs with left endpoint $i$ and $k < j < n$. Therefore

$$|\mathcal{F}| = 1 + kn + k(n - k - 1) + 1 = 2kn - k^2 - k + 2.$$

**Proof of Theorem 3.6.1.** Let $\mathcal{F}$ be a locally $k$-wide family of arcs of $\mathbf{Z}_n$. To prove the bound, we can assume that $\mathcal{F}$ contains $\emptyset$, $\mathbf{Z}_n$ and the singletons of $\mathbf{Z}_n$. For $k = 1$, $\mathcal{F} \setminus \{\emptyset\}$ is a tree, so by Theorem 2.6.1 $|\mathcal{F}| \leq 2n$ and the bound follows.

Let $k \geq 2$. We assume that the bound holds for locally $(k-1)$-wide families of sets. We will reconstruct $\mathcal{F}$ from a $(k-1)$-wide family of arcs and a tree derived from $\mathcal{F}$.

In this proof, ordering relations involving more than two terms are assumed to be in clockwise cyclic order.

Recall that if $A = [a_l, a_r]$ and $B = [b_l, b_r]$ are incomparable arcs and $a_r \in B$, then the left overlap of $A$ with $B$ is the arc $[b_l, a_r]$. For every arc $U \in \mathcal{F}$, let $L(U)$ be an arc in $\mathcal{F}$ which is incomparable to $U$ and has maximal left overlap with $U$. If there is no such arc, let $L(U) = \emptyset$. Let $\overline{U}$ be the left overlap of $L(U)$ with $U$. Let $\overline{\mathcal{F}} = \{\overline{U} \mid U \in \mathcal{F}\}$. Let $\mathcal{F}_r = \{U \setminus \overline{U} \mid U \in \mathcal{F}\}$. Note that if $U \in \mathcal{F}$ is non-empty, then $U \setminus \overline{U}$ is non-empty.

**Lemma 3.6.3** $\mathcal{F}_r$ *is a locally 1-wide family of arcs.*

**Proof.** Let $U, V \in \mathcal{F}$ and suppose that $U \setminus \overline{U}$ intersects $V \setminus \overline{V}$. We show that $U \setminus \overline{U}$ and $V \setminus \overline{V}$ are comparable. Since $\mathbf{Z}_n \setminus \overline{\mathbf{Z}_n} = \mathbf{Z}_n$ and $\overline{\emptyset} = \emptyset$, we can assume that $U = [u_l, u_r]$ and $V = [v_l, v_r]$ for some $u_l, u_r, v_l, v_r \in \mathbf{Z}_n$. The right endpoints of $U \setminus \overline{U}$ and $V \setminus \overline{V}$ are $u_r$ and $v_r$ respectively. By interchanging $U$ and $V$ if necessary, we can assume that $u_r \in V \setminus \overline{V}$. If $u_r = v_r$, the $U \setminus \overline{U}$ and $V \setminus \overline{V}$ are comparable and we are done. Suppose that $u_r \neq v_r$. If $U$ and $V$ are incomparable, then $U$ overlaps $V$ from the left and the left overlap of $U$ with $V$ strictly includes $\overline{V}$, contradicting the definition of $\overline{V}$. Therefore $U \subset V$. This implies that $L(V)$ overlaps $U$ from the left. Since the left overlap of $L(V)$ with $U$ is included in $\overline{U}$, it follows that $U \setminus \overline{U} \subseteq V \setminus \overline{V}$, as desired. ∎

**Lemma 3.6.4** $\overline{\mathcal{F}}$ *is a locally $(k-1)$-wide family of arcs.*

**Proof.** Let $\overline{\mathcal{U}} = \{\overline{U_1}, \ldots, \overline{U_l}\}$ be an antichain of $\overline{\mathcal{F}}$ with $\bigcap \overline{\mathcal{U}} \neq \emptyset$. Let $x \in \bigcap \overline{\mathcal{U}}$. Then $x \in U_i \cap L(U_i)$ for each $i$. Let $\mathcal{U}_i = \{U_1, \ldots, U_i\}$ and $L(\mathcal{U}_i) = \{L(U_1), \ldots, L(U_i)\}$. To show that $l \leq k - 1$, we will construct an antichain $\mathcal{V} \subseteq \mathcal{U}_l \cup L(\mathcal{U}_l)$ such that $|\mathcal{V}| = l + 1$. Since $\mathcal{U}_l \cup L(\mathcal{U}_l) \subseteq \mathcal{F}$ and $\mathcal{F}$ is locally $k$-wide, we can then deduce that $l + 1 \leq k$ as required.

For $1 \leq i \leq l$, let $[u_{il}, u_{ir}] = U_i$ and $[v_{il}, v_{ir}] = L(U_i)$. Then $\overline{U_i} = [u_{il}, v_{ir}]$. Since the $\overline{U_i}$ are pairwise incomparable, $u_{il} \neq u_{jl}$ and $v_{ir} \neq v_{jr}$ for $i \neq j$. By reordering we can assume that

$$u_{11} < u_{21} < \ldots < u_{l1} \leq x.$$

This implies that

$$x \leq v_{1r} < v_{2r} < \ldots < v_{lr}.$$

We also have

$$v_{il} < u_{il} \le x \le v_{ir} < u_{ir}$$

for each $i$.

We collect some facts about the relationships between the $U_i$ and the $L(U_i)$:

   (i) For $i \le j$, $L(U_i)$ overlaps $U_j$ from the left (since $v_{il} < u_{il} \le u_{jl} \le x \le v_{ir} \le v_{jr} < u_{jr}$).

   (ii) For $i < j$, $U_i \not\subseteq U_j$ (since $u_{il} < u_{jl} \le x$).

   (iii) For $i < j$, $L(U_i) \not\supseteq L(U_j)$ (since $x \le v_{ir} < v_{jr}$).

For $i = 1, 2, \ldots, l$ we recursively construct antichains $\mathcal{V}_i \subseteq \mathcal{U}_i \cup L(\mathcal{U}_i)$ such that $|\mathcal{V}_i| = i + 1$ and $U_i \in \mathcal{V}_i$ ($\mathcal{V}_{i+1}$ need not include $\mathcal{V}_i$).

Let $\mathcal{V}_1 = \{L(U_1), U_1\}$. Assume that $i \ge 2$ and the antichain $\mathcal{V}_{i-1}$ has been constructed. By fact (ii), either $U_{i-1}$ and $U_i$ are incomparable or $U_{i-1} \supset U_i$. We consider these cases in turn.

Suppose that $U_{i-1}$ and $U_i$ are incomparable. Let $\mathcal{V}_i = \mathcal{V}_{i-1} \cup \{U_i\}$. To see that $\mathcal{V}_i$ is an antichain, first note that by fact (i), $U_i$ is incomparable to each $L(U_j) \in \mathcal{V}_{i-1}$. To show that $U_i$ is incomparable to each $U_j \in \mathcal{V}_{i-1}$, suppose that $j < i-1$ and $U_j \in \mathcal{V}_{i-1}$. Since $U_{i-1} \in \mathcal{V}_{i-1}$, $U_j$ and $U_{i-1}$ are incomparable. Since $u_{jl} < u_{(i-1)l} < u_{il}$, incomparability of $U_{i-1}$ and $U_i$ implies that $x < u_{(i-1)r} < u_{ir}$ and incomparability of $U_j$ and $U_{i-1}$ implies that $x < u_{jr} < u_{(i-1)r}$. Therefore

$$u_{jl} < u_{il} \le x < u_{jr} < u_{ir},$$

so that $U_i$ is incomparable to $U_j$ as required.

Suppose that $U_{i-1} \supset U_i$. Then

$$u_{(i-1)l} < u_{il} \le x \le v_{ir} < u_{ir} \le u_{(i-1)r}.$$

Let $\mathcal{V}' = \mathcal{V}_{i-1} \setminus \{U_{i-1}\}$ and $\mathcal{V}_i = \mathcal{V}' \cup \{U_i, L(U_i)\}$. To show that $U_i$ and $L(U_i)$ are incomparable to the members of $\mathcal{V}'$, we first show that

$$(*) \qquad\qquad u_{(i-1)l} \le v_{il} < u_{il} \le x.$$

To that end, assume that $u_{(i-1)l} - 1 \bmod(n) \in L(U_i)$. Then $L(U_i)$ is incomparable to $U_{i-1}$ and has left overlap $[u_{(i-1)l}, v_{ir}]$ with $U_{i-1}$. Since $\overline{U_{i-1}} = [u_{(i-1)l}, v_{(i-1)r}] \subset [u_{(i-1)l}, v_{ir}]$, this contradicts the maximality condition in the definition of $\overline{U_{i-1}}$. Hence $u_{(i-1)l} - 1 \bmod(n) \notin L(U_i)$, and $(*)$ follows.

Using $(*)$ and the other known relationships between the endpoints, for $j \le i-1$ we obtain

$$v_{jl} < u_{jl} \le u_{(i-1)l} \le v_{il} < u_{il} \le x \le v_{jr} < v_{ir} < u_{ir} \le u_{(i-1)r}.$$

This implies that $U_i$ and $L(U_i)$ are incomparable to $L(U_j)$.

If $U_j \in \mathcal{V}'$, then $U_j$ and $U_{i-1}$ are incomparable, so that $U_j$ overlaps $U_{i-1}$ from the left. By the maximality condition in the definition of $\overline{U_{i-1}}$, $u_{j\mathrm{r}} \leq v_{(i-1)\mathrm{r}} < u_{(i-1)\mathrm{l}}$. Therefore, using $(*)$,

$$u_{j\mathrm{l}} < u_{(i-1)\mathrm{l}} \leq v_{i\mathrm{l}} < u_{i\mathrm{l}} \leq x < u_{j\mathrm{r}} \leq v_{(i-1)\mathrm{r}} < v_{i\mathrm{r}} < u_{i\mathrm{r}},$$

which implies that $U_i$ and $L(U_i)$ are incomparable to $U_j$ as required.

This completes the construction of the antichains $\mathcal{V}_i$. ∎

We construct a surjective map $\sigma$ from the disjoint union of $\mathcal{F}_\mathrm{r}$ and $\overline{\mathcal{F}} \setminus \{\emptyset\}$ onto $\mathcal{F}$. For $A \in \mathcal{F}_\mathrm{r}$, let $\sigma(A)$ be the minimal member of the chain of arcs $U \in \mathcal{F}$ with $U \setminus \overline{U} = A$. For $A \in \overline{\mathcal{F}} \setminus \{\emptyset\}$, let $\sigma(A)$ be the minimal member of the chain of arcs with $\overline{U} = A$. (Note that the arcs of $\mathbf{Z}_n$ with left (right) endpoint $i$ form a chain for each $i \in \mathbf{Z}_n$.)

To show that $\sigma$ is onto, let $U \in \mathcal{F}$. If $\overline{U} = \emptyset$, $U$ is the smallest arc with $U \setminus \overline{U} = U \in \mathcal{F}_r$, whence $\sigma(U \setminus \overline{U}) = U$. Suppose that $\overline{U} \neq \emptyset$ and $\sigma(\overline{U}) \neq U$. Then $U \neq \mathbf{Z}_n$ so $U = [u_\mathrm{l}, u_\mathrm{r}]$ for some $u_\mathrm{l}, u_\mathrm{r} \in \mathbf{Z}_n$ and there is an arc $V = [v_\mathrm{l}, v_\mathrm{r}] \in \mathcal{F}$ such that $\overline{V} = \overline{U}$ and $V \subset U$. We have $V \setminus \overline{V} \neq \emptyset$ and $U \setminus \overline{U} \neq \emptyset$. Let $W = [w_\mathrm{l}, w_\mathrm{r}]$ be an arc in $\mathcal{F}$ with $W \setminus \overline{W} = U \setminus \overline{U}$. Suppose that $W \subset U$ and let $x$ be the left endpoint of $U \setminus \overline{U}$. Then

$$u_\mathrm{l} = v_\mathrm{l} < w_\mathrm{l} \leq x \leq v_\mathrm{r} < u_\mathrm{r} = w_\mathrm{r}.$$

Thus $V$ overlaps $W$ from the left so that by definition of $\overline{W}$, $v_\mathrm{r} \in \overline{W}$. However, $\overline{V} = \overline{U}$ implies that $v_\mathrm{r} \in U \setminus \overline{U}$, contradicting the assumption that $U \setminus \overline{U} = W \setminus \overline{W}$. It follows that $W \supseteq U$, so that $U$ is minimal among the arcs $W \in \mathcal{F}$ with $W \setminus \overline{W} = U \setminus \overline{U}$ and thus $\sigma(U \setminus \overline{U}) = U$ as required.

Lemma 3.6.3 implies that $|\mathcal{F}_\mathrm{r}| \leq 2n$. Since $\emptyset \in \overline{\mathcal{F}}$, Lemma 3.6.4 implies that $\left|\overline{\mathcal{F}} \setminus \{\emptyset\}\right| \leq A(k-1, n) - 1$. Therefore

$$
\begin{aligned}
|\mathcal{F}| &\leq |\mathcal{F}_\mathrm{r}| + \left|\overline{\mathcal{F}} \setminus \{\emptyset\}\right| \\
&\leq 2n + A(k-1, n) - 1,
\end{aligned}
$$

so that by arbitrariness of $\mathcal{F}$,

$$A(k, n) \leq 2n - 1 + A(k-1, n).$$

Iterating this recursive inequality and using the fact that $\mathrm{A}(1, n) = 2n$, we get

$$\mathrm{A}(k, n) \leq 2kn - k + 1,$$

as required. ∎

**Theorem 3.6.5** *If $n \geq 2k$ and $\mathcal{F}$ is a locally $k$-wide family of segments of $[n]$, then $\mathcal{F} \leq 2kn - 2k^2 + k + 1$. This bound is best possible.*

This result follows from the next theorem, which removes the restriction on $n$. Let $S(k, n)$ be the maximum cardinality of a locally $k$-wide family of segments of $[n]$.

**Theorem 3.6.6** *The function* $S(k, n)$ *satisfies the recursive equations*

$$
\begin{array}{llll}
S(k, 0) & = & 1, & \\
S(0, n) & = & 1, & \\
S(k, 1) & = & 2 & \text{for } k \geq 1, \\
S(k, n) & = & 2n - 1 + S(k - 1, n - 2) & \text{for } k \geq 1 \text{ and } n \geq 2.
\end{array}
$$

To see that this generalizes the bound of $2kn - 2k^2 + k + 1$ for $n \geq 2k$ asserted in Theorem 3.6.5, let $S'(k, n) = 2kn - 2k^2 + k + 1$. We show that the function $S'$ solves the recursive equations in the statement of Theorem 3.6.6 for $2k \leq n$. Note that if $2k \leq n$, then $2(k - 1) \leq (n - 2)$. Compute $S'(k, n)$ as follows:

$$
\begin{array}{lll}
S'(k, n) & = & 2kn - 2k^2 + k + 1 \\
& = & 2n - 1 + 2(k - 1)(n - 2) - 2(k - 1)^2 + (k - 1) + 1 \\
& = & 2n - 1 + S'(k - 1, n - 2).
\end{array}
$$

The only relevant boundary condition is $S(0, n) = 1$ and since $S'(0, n) = 1$, this is satisfied.

**Proof of Theorem 3.6.6.** The bound is attained by the following example:

**Example 3.6.7** Let

$$
\mathcal{H}_{k,n} = \{[i, j] \subseteq [n] \mid j \leq i + k - 1\} \cup \{[i, j] \subseteq [n] \mid 1 \leq i \leq k\}.
$$

By cyclic closure, $\mathcal{H}_{k,n}$ is equivalent to the family of arcs of $\mathbf{Z}_n$ in Example 3.6.2 excluding the arcs $A$ with $A \neq \mathbf{Z}_n$, $0 \in A$ and $n - 1 \in A$. In particular, $\mathcal{H}_{k,n}$ is locally $k$-wide.

To show that $\left|\mathcal{H}_{k,n}\right| = S(k, n)$, observe that $\left|\mathcal{H}_{k,0}\right| = \left|\mathcal{H}_{0,n}\right| = 1$, $\left|\mathcal{H}_{k,1}\right| = 2$ (for $k \geq 1$) and $\left|\mathcal{H}_{1,n}\right| = 2n$. Let $k \geq 2$ and $n \geq 1$. Let $\mathcal{G}$ be the disjoint union of $\mathcal{H}_{1,n}$ and $\mathcal{H}_{k-1,n-2} \setminus \{\emptyset\}$. Define $\sigma : \mathcal{H}_{k,n} \to \mathcal{G}$ by

$$
\sigma([i, j]) = \begin{cases} [i, j] \in \mathcal{H}_{1,n} & \text{if } [i, j] = \emptyset, \ i = j \text{ or } i = 1, \\ [i - 1, j - 2] \in \mathcal{H}_{k-1,n-2} \setminus \{\emptyset\} & \text{otherwise.} \end{cases}
$$

Then $\sigma$ is a bijection. It follows that

$$
\begin{array}{lll}
\left|\mathcal{H}_{k,n}\right| & = & \left|\mathcal{H}_{1,n}\right| + \left|\mathcal{H}_{k-1,n-2}\right| - 1 \\
& = & 2n - 1 + \left|\mathcal{H}_{k-1,n-2}\right|.
\end{array}
$$

This shows that the function $S'$ defined by $S'(k, n) = \left|\mathcal{H}_{k,n}\right|$ satisfies the recursive identities of the theorem, so that $\left|\mathcal{H}_{k,n}\right| = S(k, n)$ as required.

**Proof of the bound.** Let $\mathcal{F}$ be a maximal locally $k$-wide family of segments of $[n]$. If $k = 0$ or $n = 0$, then $\mathcal{F} = \{\emptyset\}$. If $n = 1$ and $k \geq 1$, then $\mathcal{F} = \{\emptyset, \{1\}\}$, so the boundary conditions on the function $S$ are satisfied. If $k = 1$ and $n \geq 2$, $\mathcal{F} \setminus \{\emptyset\}$ is a tree, so that $|\mathcal{F}| \leq 2n$. Since $S(1, n) = 2n - 1 + S(0, n - 2) = 2n$, it remains to consider the case $k \geq 2$ and $n \geq 1$. We assume that the bound holds for locally $(k - 1)$-wide families of segments.

By considering $\mathcal{F}$ as a family of arcs on the cyclic closure of $[n]$, we can apply the reduction of the proof of Theorem 3.6.1 to $\mathcal{F}$. That is, we obtain a locally 1-wide family $\mathcal{F}_\mathrm{r}$ and a locally $(k - 1)$-wide family $\overline{\mathcal{F}}$ such that there is a surjective map from the disjoint union of $\mathcal{F}_\mathrm{r}$ and $\overline{\mathcal{F}} \setminus \{\emptyset\}$ to $\mathcal{F}$. The fact that $\mathcal{F}$ consists of segments implies that both $\mathcal{F}_\mathrm{r}$ and $\overline{\mathcal{F}}$ consist of segments. Furthermore (using notation from the proof of Theorem 3.6.1), if $U$ is a segment containing $n$, no segment which overlaps $U$ from the left contains $n$, hence $\overline{U}$ does not contain $n$. If $U$ is a segment containing 1, no segment overlaps $U$ from the left, so $\overline{U}$ is empty. This implies that $\overline{\mathcal{F}}$ consists entirely of segments included in $[2, n - 1]$. Therefore

$$\begin{aligned} |\mathcal{F}| &\leq |\mathcal{F}_\mathrm{r}| + \left|\overline{\mathcal{F}} \setminus \{\emptyset\}\right| \\ &\leq 2n + S(k - 1, n - 2) - 1. \end{aligned}$$

By arbitrariness of $\mathcal{F}$,

$$S(k, n) \leq 2n - 1 + S(k - 1, n - 2)$$

as required. ∎

## 3.7 Maximal $r$-antichains and Sperner Closure

The results in this section are used in Section 3.8 to obtain a bound on the size of locally $k$-wide families of sets.

Let $P$ be a poset. Let $\mathcal{A}$ be the family of antichains of $P$. We consider $\mathcal{A}$ ordered by $A \leq B$ iff $[A) \supseteq [B)$ (the filter order). Thus $\mathcal{A}$ is a (distributive) lattice where $A \wedge B$ consists of the minimal elements of $A \cup B$. If $\mathcal{B} \subseteq \mathcal{A}$, then the antichain $A \in \mathcal{B}$ is maximal in $\mathcal{B}$ iff for every $B \in \mathcal{B}$, $A \not< B$. Note that this notion of maximality is different from the one induced by the inclusion order on $\mathcal{A}$.

**Definition.** The antichain $A$ is a *maximal $r$-antichain* of $P$ iff $A$ is maximal among the $r$-element antichains of $P$. A *maximal $*$-antichain* is a maximal $r$-antichain for some $r$. Let $\mathcal{C}_{-1}(P)$ be the family of maximal $*$-antichains of $P$.

If $A$ is an antichain, then $A$ is the minimal antichain of the filter generated by $A$. Therefore:

**Observation 3.7.1** *$A$ is a maximal $r$-antichain iff for every antichain $B \subseteq [A)$ with $B \neq A$, $|B| \leq r - 1$.*

Suppose that $A$ and $B$ are antichains and $B \subseteq [A)$. Then $(A \setminus (B]) \cup B$ is an antichain. We have

$$\left|(A \setminus (B]) \cup B\right| = |A| - \left|A \cap (B]\right| + |B|.$$

Thus $\left|(A \setminus (B]) \cup B\right| < |A|$ iff $\left|A \cap (B]\right| > |B|$. This implies:

**Observation 3.7.2** *$A$ is a maximal $*$-antichain iff for every antichain $B \subseteq [A)$ with $B \not\subseteq A$, $\left|A \cap (B]\right| > |B|$.*

**Theorem 3.7.3** *If $A$ and $B$ are maximal $*$-antichains, then $A \wedge B$ is a maximal $*$-antichain.*

This theorem implies that $\mathcal{C}_{-1}$ is a subsemilattice of $\mathcal{A}$.

**Proof.**   Let $A$ and $B$ be maximal $*$-antichains and define

$$
\begin{array}{llll}
A_1 & = & A \setminus ((B] \cup [B)), & \qquad B_1 & = & B \setminus ((A] \cup [A)), \\
A_2 & = & A \cap ((B] \setminus B), & \qquad B_2 & = & B \cap ((A] \setminus A), \\
A_3 & = & A \cap ([B) \setminus B), & \qquad B_3 & = & B \cap ([A) \setminus A).
\end{array}
$$

Then $\{A \cap B, A_1, A_2, A_3, B_1, B_2, B_3\}$ partitions $A \cup B$ and

$$
\begin{array}{lll}
A & = & (A \cap B) \cup A_1 \cup A_2 \cup A_3, \\
B & = & (A \cap B) \cup B_1 \cup B_2 \cup B_3, \\
A \wedge B & = & (A \cap B) \cup A_1 \cup A_2 \cup B_1 \cup B_2.
\end{array}
$$

To show that $A \wedge B$ is a maximal $*$-antichain, suppose that $D$ is an antichain of $[A \wedge B)$. Let $m_1 = \left|((A \cap B) \cup A_1 \cup A_2) \cap (D]\right|$, $m_2 = \left|A_3 \cap (D]\right|$, $n_1 = \left|D \cap [A)\right|$ and $n_2 = \left|D \setminus [A)\right|$. Since $m_1 + m_2$ is the number of elements of $A$ below the antichain $D \cap [A)$, Observation 3.7.2 implies that

$$(*) \qquad
\begin{array}{lll}
m_1 + m_2 & = & n_1 \quad \text{if } D \cap [A) \subseteq A, \\
m_1 + m_2 & > & n_1 \quad \text{otherwise.}
\end{array}
$$

Let $m_3 = \left|(B_1 \cup B_2) \cap (D]\right| = \left|(A \wedge B) \cap (D]\right| - m_1$. Since $(B_1 \cup B_2) \cap (D]$ contains the elements of $B$ below the antichain $(D \setminus [A)) \cup (A_3 \cap (D])$, Observation 3.7.2 implies that

$$(**) \qquad
\begin{array}{lll}
m_3 & \geq & n_2 + m_2 \quad \text{if } (D \setminus [A)) \cup (A_3 \cap (D]) \subseteq B, \\
m_3 & > & n_2 + m_2 \quad \text{otherwise.}
\end{array}
$$

If $D \cap [A) \not\subseteq A$ or $(D \setminus [A)) \cup (A_3 \cap (D]) \not\subseteq B$, then we can combine the inequalities above to obtain

$$\left|(A \wedge B) \cap (D]\right| = m_1 + m_3 > n_1 + n_2 = |D|.$$

Suppose that $D \cap [A) \subseteq A$ and $(D \setminus [A)) \cup (A_3 \cap (D]) \subseteq B$. Then $D \subseteq A \cup B$. Since $B_3 \subseteq [A) \setminus A$, $D \cap B_3 = \emptyset$. Since $A_3 \cap B = \emptyset$, $D \cap A_3 = \emptyset$. Therefore $D \subseteq (A \cap B) \cup A_1 \cup A_2 \cup B_1 \cup B_2 = A \wedge B$.

By arbitrariness of $D$ and Observation 3.7.2, $A \wedge B$ is a maximal $*$-antichain. ∎

If $A$ and $B$ are incomparable antichains, then $A \wedge B$ is strictly below both $A$ and $B$. Therefore, Observation 3.7.1 implies:

**Corollary 3.7.4** *If $A$ is a maximal $r$-antichain and $B$ is a maximal $r'$-antichain where $A$ and $B$ are incomparable, then $A \wedge B$ is a maximal $q$-antichain for some $q > \max(r, r')$.*

Since the Sperner antichains of $P$ are the $\mathrm{w}(P)$-element antichains of $P$, the (unique) maximal Sperner antichain of $P$ is the least member of $\mathcal{C}_{-1}$.

**Corollary 3.7.5** *If $A$ is the maximal Sperner antichain of $P$, and $B$ is a maximal $*$-antichain of $P$, then $A \leq B$.*

**Definition.** Let $\mathcal{C}_i$ be the family of antichains $A \in \mathcal{A}$ such that if $B$ is an antichain of $[A)$ and $B \not\subseteq A$, then $|B| \leq |A| + i$. For $i = -1$, this agrees with the definition of $\mathcal{C}_{-1}$. If $i \geq 0$, the condition $B \not\subseteq A$ in the definition of $\mathcal{C}_i$ is unnecessary. Note that for $i < j$, $\mathcal{C}_i \subseteq \mathcal{C}_j$.

The next observation generalizes Observation 3.7.2.

**Observation 3.7.6** *The antichain $A$ is in $\mathcal{C}_i$ iff for every antichain $B \subseteq [A)$ with $B \not\subseteq A$, $\left| A \cap (B] \right| + i \geq |B|$.*

The proof of Theorem 3.7.3 can be adapted to show the following result:

**Theorem 3.7.7** *If $A \in \mathcal{C}_i$, $B \in \mathcal{C}_j$ and $i, j \geq 0$, then $A \wedge B \in \mathcal{C}_{i+j}$. If $i < 0$ or $j < 0$, then $A \wedge B \in \mathcal{C}_{\max(i,j)}$.*

This implies that $\mathcal{C}_0$ is also a meet-subsemilattice of $\mathcal{A}$.

**Proof.** Suppose that $i, j \geq 0$. If we follow the proof of Theorem 3.7.3, using the same notation and applying Observation 3.7.6 instead of Observation 3.7.2, then $(*)$ is replaced by

$(*)'$ $\qquad\qquad\qquad\qquad m_1 + m_2 + i \geq n_1,$

and $(**)$ is replaced by

$(**)'$ $\qquad\qquad\qquad\qquad m_3 + j \geq n_2 + m_2.$

The conclusion is

$$\left| (A \wedge B) \cap (D] \right| = m_1 + m_3 + i + j \geq n_1 + n_2 = |D|,$$

which proves the result for $i, j \geq 0$.

If $i < 0$ and $j \geq 0$, then $A \in \mathcal{C}_0$, so $A \wedge B \in \mathcal{C}_{0+j} = \mathcal{C}_{\max(i,j)}$ as required.

Suppose that $i < 0$ and $j < 0$. Then $(*)'$ holds if $D \cap [A) \nsubseteq A$ and $(**)'$ holds if $(D \setminus [A)) \cup (A_3 \cap (D]) \nsubseteq B$. Therefore, as in the proof of Theorem 3.7.3,

$$\big|(A \wedge B) \cap (D]\big| = m_1 + m_3 + \max(i,j) \geq n_1 + n_2 = \big|D\big|,$$

unless $D \subseteq A \wedge B$, as required. ∎

**Theorem 3.7.8** *Let* $w = \mathrm{w}(P)$. *Then* $\big|\bigcup \mathcal{C}_{-1}\big| \leq \frac{w(w+1)}{2}$. *This bound is best possible.*

**Proof.** The bound is attained by the following partial order:

**Example 3.7.9** Let $P$ be the union of $w$ disjoint antichains $A_1, A_2, \ldots, A_w$ where $A_i$ consists of $w - i + 1$ elements. Let $x < y$ iff $x \in A_i$ and $y \in A_j$ for $i < j$. Then $\big|P\big| = \frac{w(w+1)}{2}$ and for each $i$, $A(i)$ is the only maximal $(w - i + 1)$-antichain of $P$. Figure 3.1 shows $P$ for $w = 3$.

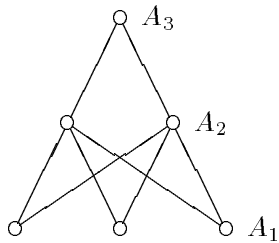

Figure 3.1

**Proof of the bound.** If $A$ is a maximal $r$-antichain of $P$, then $A$ is a maximal $r$-antichain of every subposet of $P$ which includes $A$. Hence we can assume that $P = \bigcup \mathcal{C}_{-1}$ and bound $\big|P\big|$.

For every $x \in P$, let $\mathcal{C}(x) = \{A \in \mathcal{C}_{-1} \mid x \in A\}$. By assumption, $\mathcal{C}(x) \neq \emptyset$. If $A, B \in \mathcal{C}(x)$, then $x \in A \wedge B$ and by Theorem 3.7.3, $A \wedge B \in \mathcal{C}_{-1}$, which implies that $A \wedge B \in \mathcal{C}(x)$. Thus $\mathcal{C}(x)$ is closed under meets. Let $A(x)$ be the minimal member of $\mathcal{C}(x)$. Then $A(x)$ is the maximum cardinality antichain in $\mathcal{C}(x)$.

**Lemma 3.7.10** *If* $x \notin [A(y))$, *then* $\big|A(x)\big| > \big|A(y)\big|$.

**Proof.** Suppose that $x \notin [A(y))$. Then $x \in (A(y) \wedge A(x)) \setminus A(y)$, which implies that $A(y) \wedge A(x) < A(y)$. Minimality of $A(x)$ in $\mathcal{C}(x)$ and $A(y) \wedge A(x) \in \mathcal{C}(x)$ imply that $A(x) = A(y) \wedge A(x)$. Hence $A(x) < A(y)$ and the result follows. ∎

For $1 \leq r \leq w$, let $N(r) = \{x \in P \mid |A(x)| = r\}$. The lemma implies that for $x, y \in N(r)$, $x \not\sim y$. Hence $N(r)$ is an antichain. Let $x \in N(r)$. By the lemma, if $y \in N(r)$, then $y \in [A(x))$, so that $N(r) \geq A(x)$. Since $A(x)$ is a maximal $r$-antichain, we have $|N(r)| \leq |A(x)| = r$. Using $P = \bigcup_{r=1}^{w} N(r)$, we obtain

$$
\begin{aligned}
|P| &\leq |N(1)| + |N(2)| + \ldots + |N(w)| \\
&\leq 1 + 2 + \ldots + w \\
&= \frac{w(w+1)}{2}
\end{aligned}
$$

as required. ∎

**Theorem 3.7.11** *Let $A_1, A_2, \ldots, A_l$ be a sequence of pairwise incomparable $r$-element antichains of a poset $P$ of width $w$. Suppose that for every $i < j$, $A_j$ is a maximal $r$-antichain of $A_i \cup A_j$ (with the induced order). Then*

$$
l \leq \binom{w}{r}.
$$

The bound is attained by the sequence (in any order) of all $r$-subsets of the discrete $w$-element poset.

**Proof.** Using Dilworth's Theorem 2.4.2, let $P = C_1 \cup C_2 \cup \ldots \cup C_w$ be a decomposition of $P$ into disjoint chains. For every antichain $A$ of $P$, let $c(A) = \{C_i \mid C_i \cap A \neq \emptyset\}$. Then $|c(A)| = |A| = r$. It remains to show that for $i \neq j$, $c(A_i) \neq c(A_j)$, as this implies that $l$ is at most the number of $r$-element subsets of a $w$-set.

Suppose that for some $i < j$, $c(A_i) = c(A_j)$. Then $A_i \cup A_j$ has width $r$, so $A_i$ and $A_j$ are Sperner antichains of $A_i \cup A_j$. By Theorem 2.5.2, $A_i \vee A_j$ is a Sperner antichain of $A_i \cup A_j$. Since $A_j$ is a maximal $r$-antichain of $A_i \cup A_j$, $A_i \vee A_j = A_j$, so that $A_i \leq A_j$, contradicting the incomparability assumption. ∎

**Corollary 3.7.12** *If $P$ is a poset of width $w$, then $P$ has at most $\binom{w}{r}$ maximal $r$-antichains.*

Let $\mathcal{F}$ be a family of sets ordered by inclusion.

**Definition.** If $A \subseteq \bigcup \mathcal{F}$ and $\mathcal{A}$ is the maximal Sperner antichain of $\mathcal{F}_{\supseteq A}$, then the *(Sperner) closure* in $\mathcal{F}$ of $A$ is the set $\mathrm{SC}_{\mathcal{F}}(A) = \bigcap \mathcal{A}$. If $\mathcal{F}_{\supseteq A} = \emptyset$, let $\mathrm{SC}_{\mathcal{F}}(A) = \bigcup \mathcal{F}$. Define

$$
\begin{aligned}
\mathrm{SC}(\mathcal{F}) &= \{\mathrm{SC}_{\mathcal{F}}(A) \mid A \subseteq \bigcup \mathcal{F}\}, \\
\mathrm{SC}_r(\mathcal{F}) &= \{\mathrm{SC}_{\mathcal{F}}(A) \mid \mathrm{w}(\mathcal{F}_{\supseteq A}) = r\}.
\end{aligned}
$$

The uniqueness of maximal Sperner antichains implies that $\mathrm{SC}_{\mathcal{F}}$ is a well-defined operation on the subsets of $\bigcup \mathcal{F}$.

**Observation 3.7.13** *For* $A \subseteq \bigcup \mathcal{F}$, $\mathrm{SC}_{\mathcal{F}}(A) \supseteq A$ *and* $\mathrm{SC}_{\mathcal{F}}(\mathrm{SC}_{\mathcal{F}}(A)) = \mathrm{SC}_{\mathcal{F}}(A)$.

Since $\mathcal{F}_{\supseteq A}$ is a filter of $\mathcal{F}$, we have:

**Observation 3.7.14** *If* $A \subseteq \bigcup \mathcal{F}$ *and* $\mathcal{A}$ *is the maximal Sperner antichain of* $\mathcal{F}_{\supseteq A}$, *then* $\mathcal{A}$ *is a maximal $*$-antichain of* $\mathcal{F}$.

**Theorem 3.7.15** *For* $A \subseteq B \subseteq \bigcup \mathcal{F}$, $\mathrm{SC}_{\mathcal{F}}(A) \subseteq \mathrm{SC}_{\mathcal{F}}(B)$.

**Proof.** Let $\mathcal{A}$ and $\mathcal{B}$ be the maximal Sperner antichains of $\mathcal{F}_{\supseteq A}$ and $\mathcal{F}_{\supseteq B}$ respectively. Then $\mathcal{B}$ is a maximal $*$-antichain of $\mathcal{F}_{\supseteq A}$, so by Corollary 3.7.5, $\mathcal{A} \leq \mathcal{B}$. This gives $\mathrm{SC}_{\mathcal{F}}(A) = \bigcap \mathcal{A} \subseteq \bigcap \mathcal{B} = \mathrm{SC}_{\mathcal{F}}(B)$. ∎

**Theorem 3.7.16** $\mathrm{SC}(\mathcal{F}) = \{\bigcap \mathcal{A} \mid \mathcal{A} \in \mathcal{C}_{-1}(\mathcal{F})\}$.

**Proof.** Observation 3.7.14 implies that $\mathrm{SC}(\mathcal{F}) \subseteq \{\bigcap \mathcal{A} \mid \mathcal{A} \in \mathcal{C}_{-1}(\mathcal{F})\}$. For the reverse inclusion, let $\mathcal{A}$ be a maximal $*$-antichain of $\mathcal{F}$ and let $\mathcal{B}$ be the maximal Sperner antichain of $\mathcal{F}_{\supseteq \bigcap \mathcal{A}}$. Since $\mathcal{A} \subseteq \mathcal{F}_{\supseteq \bigcap \mathcal{A}}$, $\mathcal{A}$ is a maximal $*$-antichain of $\mathcal{F}_{\supseteq \bigcap \mathcal{A}}$, hence $\mathcal{B} \leq \mathcal{A}$ (Corollary 3.7.5). This yields

$$
\bigcap \mathcal{A} \subseteq \mathrm{SC}_{\mathcal{F}}(\bigcap \mathcal{A}) = \bigcap \mathcal{B} \subseteq \bigcap \mathcal{A},
$$

so that $\bigcap \mathcal{A} = \mathrm{SC}_{\mathcal{F}}(\bigcap \mathcal{A}) \in \mathrm{SC}(\mathcal{F})$, as desired. ∎

If $\mathcal{A}$ and $\mathcal{B}$ are antichains of $\mathcal{F}$, then

$$
\bigcap(\mathcal{A} \wedge \mathcal{B}) = (\bigcap \mathcal{A}) \cap (\bigcap \mathcal{B}).
$$

Since $\mathcal{C}_{-1}(\mathcal{F})$ is a meet-subsemilattice of the family of antichains of $\mathcal{F}$, Theorem 3.7.16 implies that $\mathrm{SC}(\mathcal{F})$ is intersection-closed. This can also be deduced from the fact that $\mathrm{SC}_{\mathcal{F}}$ is a closure operation in the usual sense (Observation 3.7.13 and Theorem 3.7.15).

Corollary 3.7.4 gives the following observation:

**Observation 3.7.17** *If* $A \in \mathrm{SC}_i(\mathcal{F})$ *and* $B \in \mathrm{SC}_j(\mathcal{F})$, *then* $A \cap B \in \mathrm{SC}_k(\mathcal{F})$ *where* $k \geq \max(i, j)$. *The inequality is strict if* $A$ *and* $B$ *are incomparable.*

By Corollary 3.7.12, if $\mathrm{w}(\mathcal{F}_{\supseteq\{x\}}) = k$, then there are at most $\binom{k}{r}$ maximal $r$-antichains $\mathcal{A} \subseteq \mathcal{F}_{\supseteq\{x\}}$. This yields:

**Observation 3.7.18** *If $\mathcal{F}$ is locally $k$-wide, then $\mathrm{SC}_r(\mathcal{F})$ is locally $\binom{k}{r}$-wide.*

If $A \subseteq B \subseteq \bigcup \mathcal{F}$ and $\mathcal{A}$ and $\mathcal{B}$ are the maximal Sperner antichains of $\mathcal{F}_{\supseteq A}$ and $\mathcal{F}_{\supseteq B}$ respectively, then as in the proof of Theorem 3.7.15, $\mathcal{A} \leq \mathcal{B}$. If $A \neq B$, then $\mathcal{A} \neq \mathcal{B}$ and $|\mathcal{A}| > |\mathcal{B}|$. This implies

**Observation 3.7.19** *The families $\mathrm{SC}_r(\mathcal{F})$ are antichains of sets.*

## 3.8 The Size of Locally $k$-wide Families of Sets

**Theorem 3.8.1** *Let $k \geq 2$ and $n \geq 1$. If $\mathcal{F}$ is a locally $k$-wide family of subsets of an $n$-set, then $|\mathcal{F}| \leq (2k)^{k-1} n$.*

The bound $(2k)^{k-1} n$ is linear in $n$ but super-exponential in $k$. In general, unless $n$ is much larger than $k$ (asymptotically as $k \to \infty$, unless $\ln(\ln(n)) \gtrsim k \ln(k)$), the bound of Theorem 3.3.6 is smaller.

For $k = 2$, the theorem implies that a locally 2-wide family of sets $\mathcal{F}$ on an $n$-set satisfies $|\mathcal{F}| \leq 4n$. Since the maximum size of a locally 2-wide family of arcs on $\mathbf{Z}_n$ is at least $4n - 4$ (Example 3.6.2), this bound is optimal up to a constant (see Corollary 3.8.8 below).

Whether the bound on the size of locally $k$-wide families of sets on an $n$-set is asymptotically linear in both $k$ and $n$ is an open question. The largest known examples have $\sim 3kn$ members (asymptotically, for $1 = o(k)$ and $k = o(n)$). For comparison, note that the asymptotic bound for such families of arcs is $2kn$ (Theorem 3.6.1). These examples are obtained by restricting locally $k$-wide families of subsets of $\mathbf{Z}$ to $[1, n]$. They are constructed by enlarging the families of arcs in Example 3.6.2.

**Example 3.8.2** Let $\mathcal{S}_k$ consist of the $k$-element segments of $\mathbf{Z}$:

$$\mathcal{S}_k = \{[i, i + k - 1] \mid i \in \mathbf{Z}\}.$$

Let $\mathcal{T}_k$ consist of the segments of $\mathbf{Z}$ with left endpoint $k$:

$$\mathcal{T}_k = \{[k, l] \mid l \geq k\}.$$

For $k \geq 2$, let $\mathcal{G}_{k+1}$ be given by

$$\begin{aligned}
\mathcal{G}_{k+1} \;=\; & \{\, [x, i] \cup \{x + k\} \mid x \leq i \leq x + k - 2,\, x = 1 \bmod(2k)\,\} \\
& \cup \{\, \{x - k\} \cup [i, x] \mid x \geq i \geq x - k + 2,\, x = 0 \bmod(2k)\,\}.
\end{aligned}$$

Observe that if $A, B \in \mathcal{G}_{k+1}$, then $A$ and $B$ are either comparable or disjoint. Define $\mathcal{F}_k$ by

$$\begin{aligned}
\mathcal{F}_2 \;&=\; \{\emptyset\} \cup \{\mathbf{Z}\} \cup \mathcal{S}_1 \cup \mathcal{S}_2 \cup \mathcal{T}_1 \cup \mathcal{T}_2, \\
\mathcal{F}_{k+1} \;&=\; \mathcal{F}_k \cup \mathcal{S}_{k+1} \cup \mathcal{T}_{k+1} \cup \mathcal{G}_{k+1}.
\end{aligned}$$

**Theorem 3.8.3** *The family of sets $\mathcal{F}_k$ is locally $k$-wide for each $k \geq 2$.*

**Proof.**    The family $\mathcal{F}_2$ is an unbounded version of Example 3.6.2 with $k = 2$. Assume inductively that $\mathcal{F}_k$ is locally $k$-wide. Let $x \in \mathbf{Z}$. For each $l \geq 2$, let

$$\mathcal{H}_l = \{A \in \mathcal{F}_l \mid x \in A, \, |A| \leq l\}.$$

Since $\mathcal{F}_k$ is locally $k$-wide and by Dilworth's Theorem 2.4.2, $\mathcal{H}_k$ can be decomposed into $k$ disjoint chains $\mathcal{C}_1, \ldots, \mathcal{C}_k$. Since the maximal members of $\mathcal{H}_k$ are the $k$ $k$-element segments containing $x$, the maximal member $C_i$ of $\mathcal{C}_i$ is a $k$-element segment for each $i$. By definition of $\mathcal{F}_{k+1}$, $\mathcal{H}_{k+1}$ consists of the members of $\mathcal{H}_k$ together with the $(k+1)$-element segments (from $\mathcal{S}_{k+1}$) and the chain $\mathcal{C}_{k+1}$ of members $A$ of $\mathcal{G}_{k+1}$ with $x \in A$. The members of $\mathcal{C}_{k+1}$ are included in a unique $(k+1)$-element segment $D_{k+1}$ of $\mathbf{Z}$. Since every $k$-element segment is contained in two successive $(k+1)$-element segments, we can match each $C_i \in \{C_1, \ldots, C_k\}$ to a distinct $(k+1)$-element segment $D_i \neq D_{k+1}$ with $C_i \subset D_i$. The chains $\mathcal{C}_i'$, defined by $\mathcal{C}_i' = \mathcal{C}_i \cup \{D_i\}$ for $1 \leq i \leq k+1$, form a decomposition of $\mathcal{H}_{k+1}$ into $k+1$ disjoint chains.

Let

$$\mathcal{H}_{k+1}' = \{A \in \mathcal{F}_{k+1} \mid |A| \geq k+1\}.$$

Then $\mathcal{H}_{k+1}'$ consists of the $(k+1)$-element segments and the union of the $k+1$ disjoint chains $\{[i, j] \mid j \geq i + k\}$ for $1 \leq i \leq k+1$. The fact that $\mathcal{H}_{k+1}'$ is locally $(k+1)$-wide follows from Corollary 3.5.2. This implies that $(\mathcal{H}_{k+1}')_{\supseteq\{x\}}$ has a decomposition into disjoint chains $\mathcal{D}_1, \ldots, \mathcal{D}_{k+1}$, where the minimal member of each chain $\mathcal{D}_i$ is a $(k+1)$-element segment. Thus the minimal members of the chains $\mathcal{D}_i$ are the same as the maximal members of the chains $\mathcal{C}_i'$. This implies that the two chain decompositions can be combined to yield a chain decomposition of $(\mathcal{F}_{k+1})_{\supseteq\{x\}}$ into $k+1$ chains. By arbitrariness of $x$, $\mathcal{F}_{k+1}$ is locally $(k+1)$-wide. The induction is complete.    ∎

To obtain a locally $k$-wide family of sets on an $n$-set, we restrict $\mathcal{F}_k$ to $[n]$: let $\mathcal{F}_{k,n} = (\mathcal{F}_k)_{\cap[n]}$. Then $\mathcal{F}_{2,n}$ is the family of segments of Example 3.6.7 and has cardinality $4n - 5$ for $n \geq 4$.

Let $n \geq 2k + 2$. To estimate $\mathcal{F}_{k+1,n}$, observe that

$$\left|(\mathcal{S}_{k+1})_{\cap[n]} \setminus \mathcal{F}_{k,n}\right| = \left|\{[i, i+k] \mid k+1 \leq i \leq n-k\}\right| = n - 2k$$

and

$$\left|(\mathcal{T}_{k+1})_{\cap[n]} \setminus \left(\mathcal{F}_{k,n} \cup (\mathcal{S}_{k+1})_{\cap[n]}\right)\right| = \left|\{[k+1, l] \mid 2k+2 \leq l \leq n\}\right| = n - 2k - 1.$$

The cardinality of $R = (\mathcal{G}_{k+1})_{\cap[n]} \setminus \left(\mathcal{F}_{k,n} \cup (\mathcal{S}_{k+1})_{\cap[n]} \cup (\mathcal{T}_{k+1})_{\cap[n]}\right)$ can be estimated as

$$\left\lfloor \frac{n}{2k} \right\rfloor 2(k-1) \leq |R| \leq \left\lceil \frac{n}{2k} \right\rceil 2(k-1).$$

Since $\lfloor \frac{n}{2k} \rfloor 2(k-1) \geq (1 - \frac{1}{k})n - 2k + 2$ and $\lceil \frac{n}{2k} \rceil 2(k-1) \leq (1 - \frac{1}{k})n + 2k - 2$, this gives the following inequalities:

$$\left| \mathcal{F}_{k,n} \right| + (3 - \tfrac{1}{k})n - 6k + 1 \leq \left| \mathcal{F}_{k+1,n} \right| \leq \left| \mathcal{F}_{k,n} \right| + (3 - \tfrac{1}{k})n - 2k - 3.$$

Therefore, for $1 = o(k)$ and $k = o(n)$, $\left| \mathcal{F}_{k,n} \right| = 3kn + o(kn)$.

For fixed $k$, these examples do not in general have maximum possible size, not even asymptotically in $n$, as demonstrated by:

**Example 3.8.4** Let

$$\begin{aligned}
\mathcal{F} \quad = \quad & \{[i,j] \mid j - i + 1 \leq 4\} \cup \{[i,l] \mid i \in [1,4]\} \\
& \cup \{ \{i, i+2\} \mid i \in \mathbf{Z}\} \\
& \cup \{ \{6i, 6i+1, 6i+3\} \mid i \in \mathbf{Z}\} \\
& \cup \{ \{6i+2, 6i+4, 6i+5\} \mid i \in \mathbf{Z}\}.
\end{aligned}$$

Then $\mathcal{F}$ is locally 4-wide and $\left| \mathcal{F}_{\cap[n]} \right| \sim 9\frac{1}{3}n$, while $\mathcal{F}_{4,n} \sim (4 + (3 - \frac{1}{2}) + (3 - \frac{1}{3}))n = 9\frac{1}{6}n$.

**Proof of Theorem 3.8.1.** Let $\mathcal{F} = \{U_0 = \emptyset, U_1, U_2, \ldots, U_b\}$ be a locally $k$-wide family of sets on an $n$-set $X$. Let $\mathcal{F}_i = \{U_1, \ldots, U_i\}$. We can assume that $U_i$ is maximal in $\mathcal{F}_i$ for each $i$. Let

$$\mathcal{G}_{w,i} = (\mathrm{SC}_w(\mathcal{F}_1) \cup \ldots \cup \mathrm{SC}_w(\mathcal{F}_i)) \setminus \{\emptyset\},$$

for $1 \leq w \leq k$ and $1 \leq i \leq b$. Let $\mathcal{G}_w = \mathcal{G}_{w,b}$ and $b_w = \left| \mathcal{G}_w \right|$. Note that since $U_i \in \mathrm{SC}_1(\mathcal{F}_i)$, we have $\mathcal{G}_{1,i} = \mathcal{F}_i \setminus \{\emptyset\}$ for each $i$, so that $b = b_1 + 1$. If $i < j$, then $\mathcal{G}_{w,i} \subseteq \mathcal{G}_{w,j}$.

**Lemma 3.8.5** Let $i \leq b$ and $B \subseteq X$. If $\mathrm{w}((\mathcal{F}_i)_{\supseteq B}) \leq w'$, then $\mathrm{w}((\mathcal{G}_{w,i})_{\supseteq B}) \leq \binom{w'}{w}$.

**Proof.** Let $\{A_1, A_2, \ldots, A_r\}$ be an antichain of $(\mathcal{G}_{w,i})_{\supseteq B}$. For $1 \leq j \leq r$, let $f(j)$ be the least $j'$ such that $A_j$ is in $\mathcal{G}_{w,j'}$. Then $A_j \in \mathrm{SC}_w(\mathcal{F}_{f(j)})$ and $f(j) \leq i$ for each $j$. By reordering, we can assume that $f$ is non-decreasing. Let $\mathcal{A}_j$ be the maximal $w$-antichain of $\mathcal{F}_{f(j)}$ with $\bigcap \mathcal{A}_j = A_j$. Since antichains with incomparable intersections are incomparable in the filter order, the antichains $\mathcal{A}_j$ are pairwise incomparable. For $j' \leq j$, $\mathcal{F}_{f(j')} \subseteq \mathcal{F}_{f(j)} \subseteq \mathcal{F}_i$, so that $\mathcal{A}_j$ is a maximal $w$-antichain of $\mathcal{A}_{j'} \cup \mathcal{A}_j$ and $\mathcal{A}_j \subseteq (\mathcal{F}_i)_{\supseteq B}$. It follows that the sequence of antichains $\{\mathcal{A}_1, \ldots, \mathcal{A}_r\}$ satisfies the conditions of Theorem 3.7.11, whence $r \leq \binom{w'}{w}$. ∎

Lemma 3.8.5 implies that $\mathcal{G}_w$ is locally $\binom{k}{w}$-wide. In particular $\mathcal{G}_k$ is locally 1-wide, which gives $b_k \leq 2n - 1$. An amortized counting argument will be used to bound the $b_w$ in terms of the $b_{w'}$ for $w' > w$.

Let $1 \leq w < k$, write

$$\mathcal{G}_w = \{A_1, A_2, \ldots, A_{b_w}\}$$

and let $f(i)$ be the least index $i'$ such that $A_i$ is in $\mathcal{G}_{w,i'}$. Then $A_i \in \mathrm{SC}_w(\mathcal{F}_{f(i)})$ for each $i$. By reordering, we can assume that $f$ is non-decreasing. Let

$$\mathcal{H}_i = \{A_1, A_2, \ldots, A_i\}.$$

**Lemma 3.8.6** *If $A$ is a non-empty member of $\mathrm{SC}_w(\mathcal{F}_i)$, then $A$ is maximal in $\mathcal{G}_{w,i}$. Conversely, if $A$ is maximal in $\mathcal{G}_{w,i}$, then $A$ is in $\mathrm{SC}_{w'}(\mathcal{F}_i)$ for some $w' \geq w$.*

**Proof.**    Suppose that $A$ is a non-empty member of $\mathrm{SC}_w(\mathcal{F}_i)$. By definition, $A$ is in $\mathcal{G}_{w,i}$. Let $\mathcal{A}$ be a maximal $w$-antichain of $\mathcal{F}_i$ such that $\bigcap \mathcal{A} = A$. If $A$ is not maximal in $\mathcal{G}_{w,i}$, then there is a member $B$ of $\mathcal{G}_{w,i}$ with $B \supset A$. Consider such a $B$. There is a $w$-element antichain $\mathcal{B}$ of $\mathcal{F}_i$ such that $\bigcap \mathcal{B} = B$. Since $B \supset A$, $\mathcal{B} \not\leq \mathcal{A}$. Since $\mathcal{A}$ is a maximal $w$-antichain of $\mathcal{F}_i$, $\mathcal{B} \not> \mathcal{A}$. Theorem 3.7.11 applied to the sequence of antichains $\{\mathcal{B}, \mathcal{A}\}$ yields $w' = w((\mathcal{F}_i)_{\supseteq A}) \geq w + 1$. However, this implies that $A = \mathrm{SC}_{\mathcal{F}_i}(A) \in \mathrm{SC}_{w'}(\mathcal{F}_i)$ with $w' > w$, contradicting the fact that $\mathrm{SC}_w(\mathcal{F}_i)$ and $\mathrm{SC}_{w'}(\mathcal{F}_i)$ are disjoint.

Conversely, suppose that $A$ is maximal in $\mathcal{G}_{w,i}$. Let $B = \mathrm{SC}_{\mathcal{F}_i}(A)$. Since $\mathrm{w}((\mathcal{F}_i)_{\supseteq A}) \geq w$, $B$ is in $\mathrm{SC}_{w'}(\mathcal{F}_i)$ for some $w' \geq w$. We show that $A = B$. Since $\mathrm{w}((\mathcal{F}_{j+1})_{\supseteq B}) \leq \mathrm{w}((\mathcal{F}_j)_{\supseteq B}) + 1$ for each $j$, it follows that for some $j \leq i$, $\mathrm{w}((\mathcal{F}_j)_{\supseteq B}) = w$. We have $\mathrm{SC}_{\mathcal{F}_j}(B) \in \mathrm{SC}_w(\mathcal{F}_j) \subseteq \mathcal{G}_{w,i}$. The inclusions $A \subseteq B \subseteq \mathrm{SC}_{\mathcal{F}_j}(B)$ and maximality of $A$ imply that $A = B$.    ∎

The assumption that $f$ is non-decreasing and Lemma 3.8.6 imply that $A_i$ is maximal in $\mathcal{H}_i$ for each $i$.

We consider a procedure which processes each member $A_i$ of $\mathcal{G}_w$ in turn. At each step the procedure will perform some token transactions on three sets of tokens $T_1$, $T_2$ and $T_3$. Initially,

$$
\begin{aligned}
T_1 &= X, \\
T_2 &= \biguplus_{w' > w} \left( \mathcal{G}_{w'} \times \left\{ 1, 2, \ldots, \tbinom{w'}{w} - 1 \right\} \right), \\
T_3 &= \emptyset,
\end{aligned}
$$

where $\biguplus$ denotes disjoint union. The set of tokens $T_2$ has $\binom{w'}{w} - 1$ copies of each set in $\mathcal{G}_{w'}$ for $k \geq w' > w$. After the $i$'th step, $T_1$ and $T_3$ will satisfy:

$(*)$ $\qquad\qquad T_1 = X \setminus \bigcup \mathcal{H}_i$,

$(**)$ $\qquad\quad\; T_3 = \{A \in \mathcal{H}_i \mid A \text{ is maximal in } \mathcal{H}_i\}$.

The $i$'th step will consist of adding $A_i$ to $T_3$ (since $A_i$ is maximal in $\mathcal{H}_i$) and either

(i) removing a total of at least two members from $T_1$, $T_2$ and $T_3$ (including the members $A$ of $T_3$ with $A \subset A_i$), or

(ii) removing at least one member of $T_1$ and no tokens otherwise.

To see how this will bound $b_w$, let $m_w$ be the initial size of the set of tokens $T_2$. Note that $m_w$ is determined by the $b_{w'}$ with $w' > w$:

$$m_w = \sum_{k \geq w' > w} b_{w'} \left( \binom{w'}{w} - 1 \right).$$

The initial size of $T_1$ is $n$. No tokens are ever added to $T_1$ or $T_2$. The total number of tokens added to $T_3$ is $b_w$. After a token has been removed from $T_3$ it is never added again. Since $A_{b_w}$ is added to $T_3$ in the last step, the total number of tokens removed from $T_3$ is at most $b_w - 1$ (or 0 if $b_w = 0$). Let $r_1$ be the number of steps where only one token is removed and let $r_2$ be the number of steps where two or more tokens are removed. Then,

$$b_w = r_1 + r_2.$$

When only one token is removed, this token comes from $T_1$. Therefore,

$$r_1 \leq n.$$

Since the total number of tokens removed is at least $r_1 + 2r_2$,

$$r_1 + 2r_2 \leq n + m_w + b_w - 1.$$

Combining these (in)equalities we obtain

$$2b_w - n \leq 2b_w - r_1 \leq n + m_w + b_w - 1.$$

Solving for $b_w$:

$$b_w \leq 2n - 1 + m_w.$$

We will estimate the right-hand side of this inequality following the description of the steps of the procedure.

To describe the circumstances under which the procedure removes tokens from $T_2$, we prove the following lemma:

**Lemma 3.8.7** *Let $A_j$ be a maximal member of $\mathcal{H}_{i-1}$ such that $A_i$ and $A_j$ are incomparable and $A_i \cap A_j \neq \emptyset$. Then $B = A_i \cap A_j$ is in $\mathcal{G}_{w'}$ where $w' = \mathrm{w}((\mathcal{F}_{f(i)})_{\supseteq B}) > w$.*

**Proof.** To see that $w' > w$, let $\mathcal{A}_i$ be the maximal $w$-antichain of $\mathcal{F}_{f(i)}$ with $\bigcap \mathcal{A}_i = A_i$. Let $\mathcal{A}_j$ be a $w$-element antichain of $\mathcal{F}_{f(j)}$ with $\bigcap \mathcal{A}_j = A_j$. The inclusion $\mathcal{F}_{f(j)} \subseteq \mathcal{F}_{f(i)}$ implies that the width $w'$ of $(\mathcal{F}_{f(i)})_{\supseteq B}$ is at least the width of $\mathcal{A}_i \sqcup \mathcal{A}_j$. Since $A_i$ and $A_j$ are incomparable, the sequence of antichains $\{\mathcal{A}_j, \mathcal{A}_i\}$ satisfies the conditions of Theorem 3.7.11, whence $w' > w$.

Suppose that $A_j$ is maximal in $\mathcal{G}_{w,f(i)}$. Then by Lemma 3.8.6, $A_j \in \mathrm{SC}(\mathcal{F}_{f(i)})$. Since $\mathrm{SC}(\mathcal{F}_{f(i)})$ is intersection-closed and $A_i \in \mathrm{SC}(\mathcal{F}_{f(i)})$, $B = A_i \cap A_j \in \mathrm{SC}(\mathcal{F}_{f(i)})$. By definition of $w'$, $B \in \mathrm{SC}_{w'}(\mathcal{F}_{f(i)}) \subseteq \mathcal{G}_{w'}$, as required.

Suppose that $A_j$ is not maximal in $\mathcal{G}_{w,f(i)}$. Since $A_j$ is maximal in $\mathcal{G}_{w,f(j)}$, we have $f(j) < f(i)$, so that $A_j \in \mathcal{G}_{w,f(i)-1}$. Since $f$ is non-decreasing, $\mathcal{G}_{w,f(i)-1} \subseteq \mathcal{H}_{i-1}$, hence $A_j$ is maximal in $\mathcal{G}_{w,f(i)-1}$ and by Lemma 3.8.6, $A_j \in \mathrm{SC}(\mathcal{F}_{f(i)-1})$. Let $w''$ be the width of $(\mathcal{F}_{f(i)-1})_{\supseteq B}$. Let

$$
\begin{aligned}
B' &= \mathrm{SC}_{\mathcal{F}_{f(i)}}(A_i \cap A_j), \\
B'' &= \mathrm{SC}_{\mathcal{F}_{f(i)-1}}(A_i \cap A_j).
\end{aligned}
$$

Since $A_i \in \mathrm{SC}(\mathcal{F}_{f(i)})$, $\mathrm{SC}_{\mathcal{F}_{f(i)}}(A_i) = A_i$. Since Sperner closure preserves inclusion,

$$B' \subseteq A_i.$$

Similarly,

$$B'' \subseteq A_j.$$

The family $\mathcal{F}_{f(i)}$ is obtained from $\mathcal{F}_{f(i)-1}$ by the addition of the set $U_{f(i)}$. This implies that either $w' = w''$ or $w' = w'' + 1$. Let $\mathcal{B}'$ be the maximal $w'$-antichain of $(\mathcal{F}_{f(i)})_{\supseteq B}$ with $\bigcap \mathcal{B}' = B'$. Let $\mathcal{B}''$ be the maximal $w''$-antichain of $(\mathcal{F}_{f(i)-1})_{\supseteq B}$ with $\bigcap \mathcal{B}'' = B''$.

If $w' > w''$, then $\mathcal{B}' = \mathcal{B}'' \cup \{U_{f(i)}\}$, so that $B'' \subseteq B'$. This, together with $B'' \subseteq A_j$ and $A_i \cap A_j \subseteq B' \subseteq A_i$, implies that $B' = B = A_i \cap A_j$. Therefore $B \in \mathrm{SC}_{w'}(\mathcal{F}_{f(i)}) \subseteq \mathcal{G}_{w'}$.

If $w' = w''$, then $\mathcal{B}'' < \mathcal{B}'$ (since $\mathcal{B}'$ is the maximal Sperner antichain of $(\mathcal{F}_{f(i)})_{\supseteq B}$), so that $B'' \subseteq B'$. This, together with $B' \subseteq A_i$ and $A_i \cap A_j \subseteq B'' \subseteq A_j$, implies that $B'' = B = A_i \cap A_j$. Therefore, $B \in \mathrm{SC}_{w''}(\mathcal{F}_{f(i)-1}) \subseteq \mathcal{G}_{w''} = \mathcal{G}_{w'}$. This completes the proof of the lemma. ∎

Let

$$
\begin{aligned}
\overline{A_i} &= \bigcup \{ A \in \mathcal{H}_{i-1} \mid A \subset A_i \text{ and } A \text{ is maximal in } \mathcal{H}_{i-1} \}, \\
\mathcal{M}_i &= \{ B \mid B \text{ is maximal in } \{ A \cap A_i \mid A \in \mathcal{H}_{i-1},\ A \cap (A_i \setminus \overline{A_i}) \neq \emptyset \} \}.
\end{aligned}
$$

Note that if $A$ is maximal in $\mathcal{H}_{i-1}$ and $A \cap (A_i \setminus \overline{A_i}) \neq \emptyset$, then $A$ and $A_i$ are incomparable. If $B \in \mathcal{M}_i$, then there is a maximal member $A$ of $\mathcal{H}_{i-1}$ with $A \cap A_i = B$, so that by Lemma 3.8.7, $B \in \mathcal{G}_{w'}$ for some $w' > w$.

We now perform the token removals required in the $i$'th step of the procedure. Assume inductively that after the $(i-1)$'th step $(*)$ and $(**)$ are satisfied; and that any tokens removed from $T_2$ in the $j$'th step for $j \leq i - 1$ are (copies of) $B \in \mathcal{M}_j$ taken from (a copy of) $\mathcal{G}_{w'}$, where $w' = \mathrm{w}((\mathcal{F}_{f(i)})_{\supseteq B})$. Note that these conditions are satisfied initially (after the '0'th step). Let $\mathcal{C}$ be the family of members of $T_3$ included in $A_i$. The procedure removes the members of $\mathcal{C}$ from $T_3$ and adds $A_i$ to $T_3$. The procedure also removes the elements of $A_i \setminus \bigcup \mathcal{H}_{i-1}$ from $T_1$. If $A_i \not\subseteq \bigcup \mathcal{H}_{i-1}$, then (ii) is satisfied and the $i$'th step is complete. If $A_i \subseteq \bigcup \mathcal{H}_{i-1}$ but $|\mathcal{C}| \geq 2$, then (i) is satisfied and the $i$'th step is complete. Suppose that $A_i \subseteq \bigcup \mathcal{H}_{i-1}$ and $|\mathcal{C}| \leq 1$. Since $\overline{A}_i = \bigcup \mathcal{C}$, we have $A_i \setminus \overline{A}_i \neq \emptyset$, whence $\mathcal{M}_i \neq \emptyset$. If $|\mathcal{C}| = 0$, then $\overline{A}_i = \emptyset$, so that $|\mathcal{M}_i| \geq 2$ (by maximality of $A_i \in \mathcal{H}_i$). It follows that to satisfy (i), it suffices to remove the appropriate copy of each member of $\mathcal{M}_i$ from $T_2$.

Let $B \in \mathcal{M}_i$. Let $w' = \mathrm{w}((\mathcal{F}_{f(i)})_{\supseteq B})$. Initially, $T_2$ contains $\binom{w'}{w} - 1$ tokens which are copies of $B \in \mathcal{G}_{w'}$. To show that there are such tokens left after the $(i-1)$'th step, let $t - 1$ be the number of copies of $B \in \mathcal{G}_{w'}$ removed from $T_2$ so far, where the $j$'th copy of $B \in \mathcal{G}_{w'}$ was removed at the $g(j)$'th step. We can assume that $g$ is increasing. Let $g(t) = i$. Let $A_{g(0)}$ be the maximal member of $\mathcal{H}_{g(1)-1}$ with $A_{g(0)} \cap A_{g(1)} = B$ (using the third inductive assumption). We show that $\{A_{g(0)}, \ldots, A_{g(t)}\}$ is an antichain of $(\mathcal{G}_{w, f(i)})_{\supseteq B}$. Let $0 \leq j < j' \leq i$. Then $A_{g(j)}$ and $A_{g(j')}$ properly include $B$. Consider the set $A_{g(j)} \cap A_{g(j')}$. We have $A_{g(j)} \cap A_{g(j')} \supseteq B$. Let $A$ be a maximal member of $\mathcal{H}_{g(j')-1}$ with $A \supseteq A_{g(j)}$. Since $B \in \mathcal{M}_{g(j')}$ and by construction of $\mathcal{M}_{g(j')}$, $A \cap A_{g(j')} = B$. This implies that $A_{g(j)} \cap A_{g(j')} = B$, so that $A_{g(j)}$ and $A_{g(j')}$ are incomparable, as desired. By Lemma 3.8.5 $\mathrm{w}((\mathcal{G}_{w, f(i)})_{\supseteq B}) \leq \binom{w'}{w}$, which implies that $t + 1 \leq \binom{w'}{w}$, as required.

To obtain the bound $b \leq (2k)^{k-1}n$, we show by induction on $i$ that for $1 \leq i < k$, $b_{k-i} \leq (2k)^i n - 1$. Using the fact that $b_1 \leq 2n - 1$ observed earlier,

$$m_{k-1} \leq \left(\binom{k}{k-1} - 1\right)(2n - 1) \leq 2kn - 2n.$$

Using $b_{k-1} \leq 2n - 1 + m_{k-1}$,

$$b_{k-1} \leq 2n - 1 + (2kn - 2n) = 2kn - 1.$$

Assume that $b_{k-j} \leq (2k)^j n - 1$ for $1 \leq j < i$. Then

$$
\begin{aligned}
m_{k-i} &= \left(\binom{k}{k-i} - 1\right) b_k + \sum_{j=1}^{i-1} \left(\binom{k-j}{k-i} - 1\right) b_{k-j} \\
&\leq \binom{k}{k-i} 2n - 2n + \sum_{j=1}^{i-1} \binom{k-j}{k-i}(2k)^j n \\
&\leq \frac{k^i}{i!} 2n - 2n + \sum_{j=1}^{i-1} \frac{k^{i-j}}{(i-j)!}(2k)^j n \\
&\leq (2k)^i n \left(\frac{1}{2^{i-1} i!} + \sum_{j=1}^{i-1} \frac{1}{2^j j!}\right) - 2n.
\end{aligned}
$$

The sum in this expression can be approximated using the fact that $j! \geq 2^{j-1}$:

$$\begin{aligned}
m_{k-i} &\leq (2k)^i n \left( \frac{1}{4^{i-1}} + \sum_{j=1}^{i-1} \frac{2}{4^j} \right) - 2n \\
&= (2k)^i n \left( \frac{1}{4^{i-1}} + \frac{1}{2} \frac{1 - \frac{1}{4^{i-1}}}{1 - \frac{1}{4}} \right) - 2n.
\end{aligned}$$

For $i \geq 1$, the expression in parenthesis evaluates to

$$\frac{1}{4^{i-1}} + \frac{2}{3}\left(1 - \frac{1}{4^{i-1}}\right) = \frac{2}{3} + \frac{1}{3 \cdot 4^{i-1}} \leq 1.$$

Therefore,

$$b_{k-i} \leq 2n - 1 + m_{k-i} \leq (2k)^i n - 1.$$

This completes the induction step.

We now have

$$b = b_1 + 1 \leq (2k)^{k-1} n,$$

as desired. ∎

The estimated bound $(2k)^{k-1}n$, obtained in the proof of Theorem 3.8.1 above, can be decreased to yield a slight improvement on the bound for locally $k$-wide families of arcs in Theorem 3.6.1:

**Corollary 3.8.8** *For $n \geq 3$, the maximum size of a locally 2-wide family of sets on an $n$-set is either $4n - 4$ or $4n - 3$.*

The maximum size of an arbitrary family of sets on a 2-set is $4 = 4 \cdot 2 - 4$.

**Proof.** There are examples of such families of size $4n - 4$.

Using notation from the proof of Theorem 3.8.1, we first show that $|\mathcal{G}_2| \leq 2n - 3$. By construction, $\emptyset \notin \mathcal{G}_2$. Let $A \in \mathcal{G}_2$. Then there are incomparable sets $U$ and $V$ in $\mathcal{F}$ such that $U \cap V = A$. Let $X = \bigcup \mathcal{F}$. If $X \supset W \supset X \setminus A$, then $\{U, V, W\}$ is an antichain with non-empty intersection. Since $\mathcal{F}$ is locally 2-wide, $W \notin \mathcal{F}$. This implies that $w(\mathcal{F}_{\supseteq X \setminus A}) \leq 1$, so that no $B \in \mathcal{G}_2$ includes $X \setminus A$. It follows that $X \notin \mathcal{G}_2$ and either $\mathcal{G}_2$ does not contain every singleton of $X$, or $\mathcal{G}_2$ has at least three maximal members. Since $\mathcal{G}_2 \cup \{X\} \cup \{\{x\} \,|\, x \in X\}$ is a tree, Theorem 2.6.1 implies that $|\mathcal{G}_2| \leq 2n - 3$.

By the proof of Theorem 3.8.1,

$$\begin{aligned}
b &= b_1 + 1 \leq 2n + m_1 \\
&\leq 2n + (\binom{2}{1} - 1)(2n - 3) \\
&= 4n - 3
\end{aligned}$$

as desired. ∎

### 3.9  Notes

**Section 3.1.**   The history of the notion of centered families of sets and $k$-pseudotrees is discussed in the introduction.

**Section 3.4.**   The result of this section appears in Knill et al. [27]. The idea for the proof of the bound $(k+2)n - 1$ is due to Ehrenfeucht.

**Section 3.6.**   The proof of the linear bound for locally $k$-wide families of arcs (segments) is based on a generalization of a proof by Ehrenfeucht which shows that 2-pseudotrees of segments are linearly bounded.

**Section 3.7.**   Results on the structure of maximal $r$-antichains were discovered while establishing the linear bound for locally $k$-wide families of sets. The counting techniques used in the proof of Theorem 3.7.3 are similar to the ones used by Greene and Kleitman in [18].

# CHAPTER 4

# POSET DENSITIES IN SEMILATTICES.

## 4.1 Introduction

Let $L$ be a semilattice and $P$ a poset. Recall that $L^P$ denotes the poset of order-preserving maps from $P$ into $L$ where $f \leq g$ iff for every $x \in P$, $f(x) \leq g(x)$.

**Definition.** Let $F$ be a filter of $P$ and $a$ a join-irreducible of $L$. Define

$$\mathrm{T}(L^P, F, a) = \{f \in L^P \mid f(x) \geq a \text{ iff } x \in F\}.$$

An order-preserving map $f : P \to L$ is of *type $(F, a)$* iff $f \in \mathrm{T}(L^P, F, a)$.

**Observation 4.1.1** *For fixed $a \in L$, the families of maps $\mathrm{T}(L^P, F, a)$ partition $L^P$ into convex subfamilies, one for each filter $F$ of $P$.*

We are interested in comparing the cardinalities of the types. One way to do this is by the construction of matchings.

**Definition.** Let $Q_1$ and $Q_2$ be subsets of a poset. A *decreasing (increasing) matching* from $Q_1$ to $Q_2$ is a one-to-one map $\sigma : Q_1 \to Q_2$ such that for every $x \in Q_1$, $\sigma(x) \leq x$ $(\sigma(x) \geq x)$.

**Definition.** The pair $(L, a)$ has the *full downward matching property* for $P$, $(L, a) \in \mathrm{MPf}(P)$ for short, iff $a$ is a join-irreducible of $L$ such that for every pair of filters $F$ and $G$ of $P$ with $F \supseteq G$, there is a decreasing matching $\sigma : \mathrm{T}(L^P, F, a) \to \mathrm{T}(L^P, G, a)$.

**Definition.** The pair $(L, a)$ has the *top downward matching property* for $P$, $(L, a) \in \mathrm{MPt}(P)$ for short, iff $a$ is a join-irreducible of $L$ such that for every filter $F$ of $P$ there is a decreasing matching $\sigma : \mathrm{T}(L^P, P, a) \to \mathrm{T}(L^P, F, a)$.

**Definition.** The pair $(L, a)$ has the *weak downward matching property* for $P$, $(L, a) \in \mathrm{MPw}(P)$ for short, iff $a$ is a join-irreducible of $L$ such that for every filter $F$ of $P$, $\left| \mathrm{T}(L^P, P, a) \right| \leq \left| \mathrm{T}(L^P, F, a) \right|$.

**Definition.** The semilattice $L$ has the *full (top, weak) matching property* for $P$, $L \in \mathrm{MPf}(P)$ ($L \in \mathrm{MPt}(P)$, $L \in \mathrm{MPw}(P)$) for short, iff there exists a join-irreducible $a \in L$ such that $(L, a) \in \mathrm{MPf}(P)$ ( $(L, a) \in \mathrm{MPt}(P)$, $(L, a) \in \mathrm{MPw}(P)$).

**Observation 4.1.2** *If $(L, a) \in \mathrm{MPf}(P)$, then $(L, a) \in \mathrm{MPt}(P)$. If $(L, a) \in \mathrm{MPt}(P)$, then $(L, a) \in \mathrm{MPw}(P)$.*

For $F \neq P$, $\mathrm{T}(L^P, F, \hat{0}) = \emptyset$. This implies:

**Observation 4.1.3** *The pair $(L, \hat{0})$ does not have any of the downward matching properties.*

It follows that the one-element lattice [1] has none of the downward matching properties.

Consider the one-element poset [1]. We have $L^{[1]} \cong L$ so that

$$\mathrm{T}(L^{[1]}, [1], a) \cong [a) \quad \text{and} \quad \mathrm{T}(L^{[1]}, \emptyset, a) \cong L \setminus [a).$$

**Observation 4.1.4** *$(L, a) \in \mathrm{MPf}([1])$ iff $(L, a) \in \mathrm{MPt}([1])$; and $(L, a) \in \mathrm{MPw}([1])$ iff $\big|[a)\big| \leq \frac{1}{2}|L|$.*

A fundamental problem is to determine the (non-trivial) semilattices which have the full (top, weak) matching property for $P$. This problem appears to be very difficult even for the one-element poset. The following (in-)famous conjecture is attributed to P. Frankl:

**Conjecture 4.1.5** *If $L$ is a semilattice with $|L| \geq 2$, then $L \in \mathrm{MPw}([1])$, i.e. there is a join-irreducible $a \in L$ such that $\big|[a)\big| \leq \frac{1}{2}|L|$.*

This conjecture is usually stated in terms of union-closed families of sets (see Section 4.3). Although this problem is well known, little progress has been made toward its solution. See *Notes* at the end of this chapter for a discussion of the history and available literature.

The upward matching properties are defined analogously to the downward ones. In order for these properties to be non-trivial, we have to explicitly consider only the proper join-irreducibles (since otherwise $(L, \hat{0})$ has every upward matching property):

**Definition.** The pair $(L, a)$ has the *full upward matching property* for $P$, $(L, a) \in \mathrm{MPf}^u(P)$ for short, iff $a$ is a proper join-irreducible of $L$ such that for every pair of filters $F$ and $G$ of $P$ with $F \subseteq G$, there is an increasing matching $\sigma : \mathrm{T}(L^P, F, a) \rightarrow \mathrm{T}(L^P, G, a)$.

**Definition.** The pair $(L, a)$ has the *top upward matching property* for $P$, $(L, a) \in \mathrm{MPt}^u(P)$ for short, iff $a$ is a proper join-irreducible of $L$ such that for every filter $F$ of $P$ there is an increasing matching $\sigma : \mathrm{T}(L^P, \emptyset, a) \rightarrow \mathrm{T}(L^P, F, a)$.

**Definition.** The pair $(L, a)$ has the *weak upward matching property* for $P$, $(L, a) \in \mathrm{MPw}^u(P)$ for short, iff $a$ is a proper join-irreducible of $L$ such that for every filter $F$ of $P$, $\big|\mathrm{T}(L^P, \emptyset, a)\big| \leq \big|\mathrm{T}(L^P, F, a)\big|$.

**Definition.** The semilattice $L$ has the *full (top, weak) upward matching property* for $P$, $L \in \mathrm{MPf}^u(P)$ ($L \in \mathrm{MPt}^u(P)$, $L \in \mathrm{MPw}^u(P)$) for short, iff there exists a proper join-irreducible $a \in L$ such that $(L, a) \in \mathrm{MPf}^u(P)$ ($(L, a) \in \mathrm{MPt}^u(P)$, $(L, a) \in \mathrm{MPw}^u(P)$).

Matching properties are downward by default. Except for the next section, we will focus on the downward matching properties.

## 4.2 Examples and Counter-examples

Let $P$ be a poset.

**Example 4.2.1** Consider the one-element lattice [1] (the *trivial* lattice). The only order-preserving map from $P$ to [1] is of type $(P, \hat{0})$. Therefore, [1] has none of the downward matching properties. Since there are no proper join-irreducibles, [1] has none of the upward matching properties.

**Example 4.2.2** Consider the two-element lattice [2]. For every filter $F$ of $P$, $\mathrm{T}([2]^P, F, \hat{1})$ consists of exactly one map $f_F$ defined by $f_F(x) = \hat{1}$ if $x \in F$ and $f_F(x) = \hat{0}$ otherwise. Since $f_F \le f_G$ iff $F \subseteq G$, [2] has all the matching properties.

**Theorem 4.2.3** *If $(L, a)$ has one of the matching properties for $P$ and $M$ is a semilattice, then $(L \times M, \langle a, \hat{0} \rangle)$ has that matching property for $P$.*

**Proof.** Observe that for every join-irreducible $a \in L$,

$$\mathrm{T}((L \times M)^P, F, \langle a, \hat{0} \rangle) \cong \mathrm{T}(L^P, F, a) \times M^P.$$

If $\sigma : \mathrm{T}(L^P, F, a) \to \mathrm{T}(L^P, G, a)$ is a decreasing (increasing) matching, then the map $\sigma'$ defined by

$$\sigma'(f)(\langle u, v \rangle) = \langle \sigma(\pi_1 \circ f)(u), (\pi_2 \circ f)(v) \rangle$$

is a decreasing (increasing) matching from $\mathrm{T}((L \times M)^P, F, \langle a, \hat{0} \rangle)$ to $\mathrm{T}((L \times M)^P, G, \langle a, \hat{0} \rangle)$. ∎

Theorem 4.2.3 will be generalized in Section 4.4.

**Example 4.2.4** Let $n \ge 1$. The Boolean lattice $B_n$ generated by $n$ atoms is (isomorphic to) the $n$-fold cartesian product of the two-element lattice. Example 4.2.2 and Theorem 4.2.3 imply that $(B_n, a) \in \mathrm{MPf}(P)$ and $(B_n, a) \in \mathrm{MPf}^u(P)$ for every atom $a \in B_n$. It follows that for every filter $F$ of $P$, $\left| \mathrm{T}(B_n^P, F, a) \right| = \left| \mathrm{T}(B_n^P, P, a) \right| = \left| B_{n-1}^Q \right|$. Hence, by induction, $\left| B_n^P \right| = \left| [2]^P \right|^n$.

The following theorem together with the Fundamental Theorem of Distributive Lattices (Theorem 2.5.1) shows that every non-trivial distributive lattice has the full downward and upward matching properties.

**Theorem 4.2.5** *Let $Q$ be a poset. Let $L$ be the family of ideals of $Q$. If $x$ is maximal in $Q$, then $(L, (x]) \in \mathrm{MPf}(P)$. If $y$ is minimal in $Q$, then $(L, (y]) \in \mathrm{MPf}^u(P)$.*

**Proof.** If $z \in Q$ and $I$ is an ideal of $Q$ with $I \subset (z]$, then $z \notin I$. It follows that $(z]$ is join-irreducible in $L$ for every $z \in Q$. Let $x$ be a maximal element of $Q$. Suppose that $F$ and $G$ are filters of $P$ such that $F \supseteq G$. Define $\sigma : \mathrm{T}(L^P, F, (x]) \to \mathrm{T}(L^P, G, (x])$ by

$$\sigma(f)(u) = \begin{cases} f(u) \setminus \{x\} & \text{if } u \in F \setminus G, \\ f(u) & \text{otherwise.} \end{cases}$$

Since $x$ is maximal, for every ideal $I$ of $Q$, $I \setminus \{x\}$ is an ideal. Thus $\sigma$ is a decreasing matching from $\mathrm{T}(L^P, F, (x])$ to $\mathrm{T}(L^P, G, (x])$.

Similarly, if $y$ is a minimal element of $Q$, then we obtain an increasing matching $\sigma' : \mathrm{T}(L^P, G, (y]) \to \mathrm{T}(L^P, F, (y])$ by defining

$$\sigma'(f)(u) = \begin{cases} f(u) \cup \{y\} & \text{if } u \in F \setminus G, \\ f(u) & \text{otherwise.} \end{cases}$$

Note that by minimality of $y$, if $I$ is an ideal of $Q$, then so is $I \cup \{y\}$. The theorem follows. ■

Many atomic (semi-)lattices $L$ satisfy $\left| \mathrm{T}(L^P, P, a) \right| < \left| \mathrm{T}(L^P, \emptyset, a) \right|$ for each atom $a$, which implies that $L \notin \mathrm{MPw}^u(P)$. This is illustrated by the following example.

**Example 4.2.6** Let $M_n$ be the semilattice consisting of $\hat{0}$ and $n$ atoms (see Figure 4.1). The only proper join-irreducibles are the atoms, and for each atom $a$, $\left| \mathrm{T}(M_n^P, P, a) \right| = 1$. For $n \geq 2$, $\left| \mathrm{T}(M_n^P, \emptyset, a) \right| \geq 2$, so that $M_n \notin \mathrm{MPw}^u(P)$.

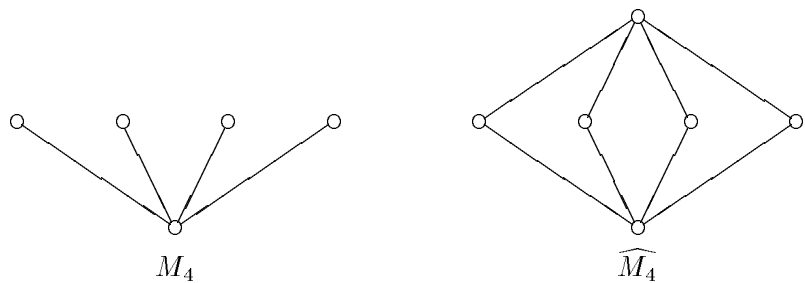

$M_4$        $\widehat{M_4}$

Figure 4.1

Let $\widehat{M_n}$ be $M_n$ with a greatest element $\hat{1}$ adjoined. Then $\widehat{M_n}$ is a (modular) lattice. For $n \geq 3$, $\left|\mathrm{T}(\widehat{M_n}^P, \emptyset, a)\right| > \left|[2]^P\right| = \left|\mathrm{T}(\widehat{M_n}^P, P, a)\right|$ for every atom $a \in \widehat{M_n}$, so that $\widehat{M_n} \notin \mathrm{MPw}^u(P)$.

A weak version of the full downward matching property requires only that for some join-irreducible $a \in L$, the function $f : [2]^P \to \mathbf{N}$ defined by $f(F) = \left|\mathrm{T}(L^P, F, a)\right|$ is order-reversing, i.e. $F \subseteq G$ implies that $f(F) \geq f(G)$. The lattice $\widehat{M_n}$ is an example for which there are posets $P$ such that $\widehat{M_n}$ does not satisfy this weak version of the full downward matching property for $P$:

**Example 4.2.7** Let $P$ be the dual of $M_k$. Then $P$ consists of a $k$-element antichain $A$ covered by $\hat{1}$. The number of filters of $P$ is $2^k + 1$. Let $a$ be an atom of $\widehat{M_n}$. The complement of $[a]$ in $\widehat{M_n}$ is isomorphic to $M_{n-1}$. Thus

$$\left|\mathrm{T}(\widehat{M_n}^P, \emptyset, a)\right| = \left|M_{n-1}^P\right| = (n-1)(\left|[2]^P\right| - 1) + 1 = (n-1)2^k + 1.$$

Let $f \in \mathrm{T}(\widehat{M_n}^P, \{\hat{1}\}, a)$. Either $f(\hat{1}) = a$, in which case for every $x \neq \hat{1}$, $f(x) = \hat{0}$; or $f(\hat{1}) = \hat{1} \in \widehat{M_n}$, in which case the restriction of $f$ to $P \setminus \{\hat{1}\}$ can be an arbitrary map into the complement of $[a]$. This implies that $\left|\mathrm{T}(\widehat{M_n}^P, \{\hat{1}\}, a)\right| = n^k + 1$. For $k \geq 2$ and $n$ sufficiently large, $n^k + 1 > (n-1)2^k + 1$.

Whether all non-trivial (semi)lattices have the weak downward matching property for every poset remains an open question. There are lattices which do not have the top matching property for any poset.

**Example 4.2.8** Let $L$ be the union-closed family of sets generated by the (unordered) pairs of adjacent vertices of the pentagon and the empty set. Since $L$ has a least member, $L$ is a lattice with the join operation given by union. $L$ is atomic, where the atoms are the edges of the pentagon.

Let $\{a, b\}$ be an edge of the pentagon. Let $\{c, a\}$ and $\{b, d\}$ be the edges other than $\{a, b\}$ incident on $a$ or $b$. Consider the ideal $I$ of $L$ generated by $\{c, a, b, d\}$ (see Figure 4.2). Then $\left|I \cap [\{a, b\})\right| = 4$ and $\left|I \setminus [\{a, b\})\right| = 3$. This implies that there is no decreasing matching from $[\{a, b\})$ to $L \setminus [\{a, b\})$, so that $L \notin \mathrm{MPt}([1])$.

More generally, let $P$ be an arbitrary poset. The sublattice $I \cap [\{a, b\})$ of $L$ is isomorphic to the four-element Boolean lattice $B_2$. The subsemilattice $I \setminus [\{a, b\})$ of $L$ is isomorphic to $M_2$ (Example 4.2.6). Since $M_2$ is (isomorphic to) a proper subsemilattice of $B_2$, $\left|\mathrm{T}(I^P, \emptyset, \{a, b\})\right| < \left|\mathrm{T}(I^P, P, \{a, b\})\right|$. Since $I^P$ is an ideal of $L^P$, there is no decreasing matching $\sigma : \mathrm{T}(L^P, P, \{a, b\}) \to \mathrm{T}(L^P, \emptyset, \{a, b\})$.

In Section 4.11 we will show that every lattice which is isomorphic to a union-closed family of sets generated by sets of size at most two, has the weak matching property for [1]. This implies that the lattice $L$ in Example 4.2.8 satisfies $L \in \mathrm{MPw}([1])$.

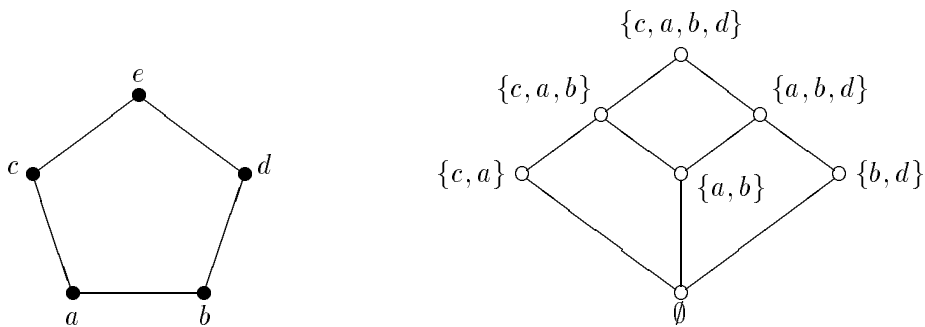

Figure 4.2

## 4.3  Semilattices, Intersection- and Union-closed Families of Sets

For many purposes, the concepts of semilattices and intersection- and union-closed families of sets are interchangeable. Every semilattice has a canonical representation as an intersection-closed family of sets where the meet operation is intersection. An equivalent representation as a union-closed family of sets (where meet is represented by union) is obtained by complementing every member of an intersection-closed representation. If $L$ is a lattice, then the canonical intersection-closed representation of the dual of $L$ is the canonical union-closed representation of $L$ where join is represented by union.

Let $L$ be a semilattice.

**Definition.**  Let $\mathcal{P}(L)$ denote the family of principal ideals of $L$. Let $J(L)$ denote the set of proper join-irreducibles of $L$. For every member $u$ of $L$, let $J(u) = \{a \in J(L) \,|\, a \le u\}$. Note that $J(u) = J(\,(u]\,) = J(L) \cap (u]$.

If $(u]$ and $(v]$ are principal ideals of $L$, then $(u] \cap (v] = (u \wedge v]$. Therefore:

**Observation 4.3.1**  *The family $\mathcal{P}(L)$ is intersection-closed. The map $\sigma : L \to \mathcal{P}(L)$ defined by $\sigma(u) = (u]$ is a semilattice isomorphism.*

**Definition.**  The *canonical intersection-closed representation* of $L$ is given by

$$\mathcal{F}(L) = (\mathcal{P}(L))_{\cap J(L)},$$

the restriction of $\mathcal{P}(L)$ to $J(L)$.

**Theorem 4.3.2** *The family $\mathcal{F}(L)$ is an intersection-closed family of sets with domain $J(L)$. The map $J : L \to \mathcal{F}(L)$ is a semilattice isomorphism with inverse given by $J^{-1}(U) = \bigvee U$ for $U \in \mathcal{F}(L)$.*

**Proof.** Since $\mathcal{F}(L)$ is the restriction of $\mathcal{P}(L)$ to $J(L)$, $\mathcal{F}(L)$ is intersection-closed. Since $J$ is the composition of the isomorphism from $L$ to $\mathcal{P}(L)$ with the restriction map $(u] \mapsto (u] \cap J(L)$, it follows that $J$ is a meet-homomorphism. Every member $u$ of $L$ is the least upper bound of the proper join-irreducibles below $u$, i.e. $u = \bigvee J(u)$. This implies that $J$ is a semilattice isomorphism. ∎

**Theorem 4.3.3** *The intersection-closed family $\mathcal{F} = \mathcal{F}(L)$ has the following property: If $a \in \bigcup \mathcal{F}$ and $U$ is the minimal member of $\mathcal{F}$ with $a \in U$, then $U \setminus \{a\} \in \mathcal{F}$.*

**Proof.** Let $a \in \bigcup \mathcal{F} = J(L)$. The minimal member of $\mathcal{F}_{\supseteq \{a\}}$ is $J(a) = (a] \cap J(L)$. Since $a$ is a proper join-irreducible of $L$, there is a unique element $u$ of $L$ covered by $a$. Every member of $L$ strictly below $a$ is below $u$. Therefore, $J(u) = J(a) \setminus \{a\}$. ∎

**Definition.** The family of sets $\mathcal{F}$ is an *intersection-closed representation* of $L$ iff $\mathcal{F}$ is intersection-closed, and $\mathcal{F}$ and $L$ are isomorphic semilattices. The representation $\mathcal{F}$ of $L$ is *irredundant* iff no restriction of $\mathcal{F}$ to a proper subset of the domain is isomorphic to $\mathcal{F}$.

Theorem 4.3.3 shows that $\mathcal{F}(L)$ is an irredundant representation of $L$. Every intersection-closed family of sets has an irredundant representation obtained by restriction. The following theorem and its proof show how to find such representations.

**Theorem 4.3.4** *Let $L$ be an intersection-closed family of sets. Then there is a one-to-one map $p : J(L) \to \bigcup L$ such that for each $a \in J(L)$, $p(a) \in a$ and the map $U \mapsto \{p(a) \mid a \in U\}$ is an isomorphism from $\mathcal{F}(L)$ onto the restriction of $L$ to the image of $p$.*

**Proof.** For each proper join-irreducible $a \in L$, let $c(a)$ be the unique member of $L$ covered by $a$ and let $p(a)$ be an (arbitrary) element of $a \setminus c(a)$. To show that $p$ is one-to-one, assume that $p(a) = p(b)$. Then $p(a) \in a \cap b$. Since $p(a) \notin c(a)$, $a \cap b = a$. Similarly, $a \cap b = b$, hence $a = b$.

Let $X$ be the image of $p$. For each $a \in J(L)$ and $u \in L$, if $p(a) \in u$, then $a \subseteq u$ (otherwise $p(a) \in u \cap a \subset a$). This implies that for each $u \in L$, $u \cap X = \{p(a) \mid a \in J(u)\}$. The result follows. ∎

Note that all irredundant intersection-closed representations of $L$ have the same number of elements in their domains (the number of proper join-irreducibles of $L$). In fact, if $\mathcal{F}$ and $\mathcal{G}$ are irredundant intersection-closed representations of $L$, then there is a bijection $\sigma : \bigcup \mathcal{F} \to \bigcup \mathcal{G}$ such that $U \in \mathcal{F}$ iff $\sigma(U) \in \mathcal{G}$ (such bijections are called *set-system isomorphisms*).

The correspondence between semilattices and intersection-closed families of sets allows us to translate properties of semilattices to properties of intersection-closed families, where the concept of a proper join-irreducible is replaced by that of an element in the domain. In general this correspondence will not be mentioned explicitly and each member $u$ of a semilattice $L$ will be identified with the set of proper join-irreducibles below $u$. Conversely, each element $x$ in the domain of an intersection-closed family $\mathcal{F}$ will often be identified with the minimal member of $\mathcal{F}_{\supseteq \{x\}}$.

**Definition.** Let $L$ be a join-semilattice. Let $M(L)$ denote the set of meet-irreducibles of $L$. The *canonical union-closed representation* of $L$ is the family $\overline{\mathcal{F}}(L)$ given by

$$\overline{\mathcal{F}}(L) = \{ M(L) \setminus [u) \mid u \in L \}.$$

Since $M(L) = J(L^*)$, $\overline{\mathcal{F}}(L) = \{ M(L) \setminus U \mid U \in \mathcal{F}(L^*) \}$ (where $L^*$ is the dual of $L$). In $\overline{\mathcal{F}}(L)$, the join operation of $L$ is represented by union. The results and definitions for canonical intersection-closed representations of meet-semilattices translate to canonical union-closed representations of join-semilattices.

Every semilattice $L$ can be made into a lattice (the *completion* of $L$) by adjoining a greatest element if necessary. Let $L$ be a lattice. Then $L$ is a join-semilattice.

**Observation 4.3.5** *The join-irreducibles of $L$ are the generators of the canonical union-closed representation of $L$.*

In Sections 4.6 and 4.7 we will use the family of generators of $\overline{\mathcal{F}}(L)$ to define a notion of locality for $L$. In those sections, we will identify each $u \in L$ with the corresponding member $M(L) \setminus [u)$ of $\overline{\mathcal{F}}(L)$.

The canonical intersection- and union-closed representations of (semi)-lattices and their duals can be used to obtain the following equivalent formulations of Conjecture 4.1.5:

**Conjecture 4.3.6** *If $\mathcal{F}$ is an intersection-closed family of sets with $|\mathcal{F}| \geq 2$, then there is an element $x \in \bigcup \mathcal{F}$ such that $\left| \mathcal{F}_{\supseteq \{x\}} \right| \leq \frac{1}{2} |\mathcal{F}|$.*

**Conjecture 4.3.7** (The union-closed sets conjecture) *If $\mathcal{F}$ is a union-closed family of sets with $|\mathcal{F}| \geq 2$, then there is an element $x \in \bigcup \mathcal{F}$ such that $\left| \mathcal{F}_{\supseteq \{x\}} \right| \geq \frac{1}{2} |\mathcal{F}|$.*

**Conjecture 4.3.8** *If $\mathcal{F}$ is a union-closed family of sets with $|\mathcal{F}| \geq 2$ and $\emptyset \in \mathcal{F}$, then there is a generator $U$ of $\mathcal{F}$ such that $\left|\mathcal{F}_{\supseteq U}\right| \leq \frac{1}{2}|\mathcal{F}|$.*

Conjectures 4.3.6 and 4.3.7 are equivalent by complementation. The irredundant intersection-closed representations of semilattices show that 4.3.6 is equivalent to 4.1.5. To see that 4.3.8 is equivalent to 4.1.5, use Observation 4.3.5 and note that if 4.1.5 holds for lattices, then it holds for semilattices. Proof: Assume that 4.1.5 holds for lattices. Let $L$ be a semilattice without a greatest element and let $\widehat{L} = L \cup \{\widehat{1}_L\}$ be the completion of $L$. By assumption, there is a join-irreducible $a \in \widehat{L}$ such that $\left|[a)_{\widehat{L}}\right| \leq \frac{1}{2}|\widehat{L}|$. Then

$$\left|[a)_L\right| = \left|[a)_{\widehat{L}}\right| - 1 < \frac{1}{2}(|\widehat{L}| - 1) = \frac{1}{2}|L|,$$

as required.

## 4.4 Subdirect Products of Semilattices

**Definition.** The semilattice $L$ is a *subdirect product* of the semilattices $L_1, \ldots, L_n$ iff

    (i) $L$ is (isomorphic to) a subsemilattice of the cartesian product $\prod_{i=1}^{n} L_i$,

    (ii) for each $i$ and $u \in L_i$, there is a member $\langle v_1, \ldots, v_n \rangle$ of $L$ such that $v_i = u$.

Condition (ii) implies that the restrictions $\pi_i \restriction L$ of the projections are onto $L_i$. If for some $i$, $\pi_i \restriction L$ is one-to-one, then $L$ is isomorphic to $L_i$. In this case, $L$ is a *trivial* subdirect product of $L_1, \ldots, L_n$.

    Since the family of subsets of an $n$-set is isomorphic to the $n$-fold cartesian product of [2], every intersection-closed family $\mathcal{F}$ with $\emptyset \in \mathcal{F}$ is a subdirect product of $\left|\bigcup \mathcal{F}\right|$ copies of [2]. Since every semilattice has an intersection-closed representation, it follows that every semilattice $L$ with $|L| > 2$ is a non-trivial subdirect product. The following example shows how to express a semilattice as a subdirect product with two components.

    Let $L$ be a semilattice.

**Example 4.4.1** Let $A_1$ and $A_2$ be subsets of $L$ such that $A_1 \cup A_2 \supseteq J(L)$. Let $L_1 = \mathcal{P}(L)_{\cap A_1}$ and $L_2 = \mathcal{P}(L)_{\cap A_2}$. For $u \in L$, let

$$\sigma(u) = \langle (u] \cap A_1, (u] \cap A_2 \rangle \in L_1 \times L_2.$$

Since $A_1 \cup A_2 \supseteq J(L)$, $\sigma$ is a one-to-one meet-homomorphism into $L_1 \times L_2$ (see Theorem 4.3.2). It follows that $L$ is a subdirect product of $L_1$ and $L_2$.

    Corollary 4.4.13 below shows that every subdirect product arises in this fashion.

**Example 4.4.2** Let $P$ be a poset. Then $L^P$ is a subdirect product of $|P|$ copies of $L$. Let $x$ be a non-maximal member of $P$ and $P' = P \setminus \{x\}$. Then every order-preserving map $f : P' \to L$ can be extended to an order-preserving map on $P$ by letting

$$f(x) = \bigwedge \{f(y) \mid y > x\}.$$

(If $L$ is a lattice, this also works for maximal $x \in P$.) This implies that $L^P$ is a subdirect product of $L$ and $L^{P'}$.

The subdirect products of Example 4.4.2 are full subdirect products:

**Definition.** The semilattice $L$ is a *full subdirect product* of the semilattices $L_1, \ldots, L_n$ iff $L$ is a subdirect product of $L_1, \ldots, L_n$ and

(iii) for every $u, v \in L$, if $w = u \vee_L v$, then $w = u \vee_{\Pi_{i=1}^n L_i} v$.

The family $\operatorname{Sub}(L)$ of subsemilattices of $L$ is intersection-closed. Thus $\operatorname{Sub}(L)$ is a semilattice with meets given by intersection. An alternate semilattice structure on $\operatorname{Sub}(L)$ is defined as follows:

**Definition.** Let $\operatorname{Sub}(L)$ denote the family of subsemilattices of $L$. For $L_1, L_2 \in \operatorname{Sub}(L)$, define $L_1 \wedge L_2$ by

$$L_1 \wedge L_2 = \{u \wedge v \mid u \in L_1, v \in L_2\}.$$

Let $L_1$ and $L_2$ be semilattices. We can define subdirect products operationally:

**Definition.** The map $f : L_1 \to \operatorname{Sub}(L_2)$ is *meet-concave* iff for every $u, v \in L_1$, $f(u) \wedge f(v) \subseteq f(u \wedge v)$. If $f : L_1 \to \operatorname{Sub}(L_2)$ is meet-concave, then define

$$L_1 \bowtie_f L_2 = \{\langle u, v \rangle \in L_1 \times L_2 \mid v \in f(u)\}.$$

**Theorem 4.4.3** *If* $f : L_1 \to \operatorname{Sub}(L_2)$ *is meet-concave, then* $L_1 \bowtie_f L_2$ *is a subsemilattice of* $L_1 \times L_2$.

**Proof.** If $\langle u_1, v_1 \rangle$ and $\langle u_2, v_2 \rangle$ are in $L_1 \bowtie_f L_2$, then by meet-concavity of $f$, $v_1 \wedge v_2 \in f(u_1 \wedge u_2)$, so that $\langle u_1 \wedge u_2, v_1 \wedge v_2 \rangle \in L_1 \bowtie_f L_2$, as desired. ∎

**Theorem 4.4.4** *The semilattice* $L$ *is a subdirect product of* $L_1$ *and* $L_2$ *iff there exists a meet-concave map* $f : L_1 \to \operatorname{Sub}(L_2)$ *such that* $L_2 = \bigcup_{u \in L_1} f(u)$ *and* $L \cong L_1 \bowtie_f L_2$.

**Proof.** Suppose that $L$ is a subdirect product of $L_1$ and $L_2$. For $u \in L_1$, let

$$f(u) = \{v \in L_2 \mid \langle u, v \rangle \in L\}.$$

Since $\langle u, v_1 \rangle \wedge \langle u, v_2 \rangle = \langle u, v_1 \wedge v_2 \rangle$, $f(u)$ is a subsemilattice of $L_2$. Since $\langle u_1, v_1 \rangle \wedge \langle u_2, v_2 \rangle = \langle u_1 \wedge u_2, v_1 \wedge v_2 \rangle$, we have $f(u_1) \wedge f(u_2) \subseteq f(u_1 \wedge u_2)$. This implies that $f$ is a meet-concave map, whence $L$ is (isomorphic to) $L_1 \bowtie_f L_2$. Since $\pi_2 \dagger L$ is onto $L_2$, $\bigcup_{u \in L_1} f(u) = L_2$.

Let $f : L_1 \to \mathrm{Sub}(L_2)$ be a meet-concave map with $\bigcup_{u \in L_1} f(u) = L_2$. By Theorem 4.4.3, $L_1 \bowtie_f L_2$ is a subsemilattice of $L_1 \times L_2$. Since every subsemilattice of $L_2$ is non-empty, $\pi_1 \dagger L_1 \bowtie_f L_2$ is onto $L_1$. Since the union of the range of $f$ is $L_2$, $\pi_2 \dagger L_1 \bowtie_f L_2$ is onto $L_2$. It follows that $L_1 \bowtie_f L_2$ is a subdirect product of $L_1$ and $L_2$. ∎

**Definition.** Let $L$ be a subdirect product of $L_1$ and $L_2$. Let $\iota_1 : L_1 \to \mathrm{Sub}(L_2)$ and $\iota_2 : L_2 \to \mathrm{Sub}(L_1)$ be the maps defined by

$$\begin{aligned} \iota_1(u) &= \{ v \in L_2 \mid \langle u, v \rangle \in L \}, \\ \iota_2(v) &= \{ u \in L_1 \mid \langle u, v \rangle \in L \}. \end{aligned}$$

For $i = 1, 2$ let $\underline{\iota_i}(w)$ be the minimal member of $\iota_i(w)$. The maps $\iota_1$ and $\iota_2$ implicitly depend on $L$, $L_1$ and $L_2$. By the proof of Theorem 4.4.4, the maps $\iota_i$ are meet-concave, $L \cong L_1 \bowtie_{\iota_1} L_2$ and $L \cong L_1 {}_{\iota_2}\!\bowtie L_2$.

Let $u \le v$ in $L_i$. Since $\iota_i$ is a meet-concave map, $\underline{\iota_i}(u) \wedge \underline{\iota_i}(v) \in \iota_i(u)$. Hence $\underline{\iota_i}(u) \le \underline{\iota_i}(u) \wedge \underline{\iota_i}(v)$, which yields $\underline{\iota_i}(u) \le \underline{\iota_i}(v)$.

**Observation 4.4.5** *The maps $\underline{\iota_1} : L_1 \to L_2$ and $\underline{\iota_2} : L_2 \to L_1$ are order-preserving.*

**Theorem 4.4.6** *Let $L$ be a subdirect product of $L_1$ and $L_2$. If $\langle u, v \rangle$ is join-irreducible in $L$, then either $u$ is join-irreducible (in $L_1$) and $v = \underline{\iota_1}(u)$ or $v$ is join-irreducible (in $L_2$) and $u = \underline{\iota_2}(v)$. If both $v = \underline{\iota_1}(u)$ and $u = \underline{\iota_2}(v)$, then $u$ and $v$ are both join-irreducible.*

**Proof.** Let $\langle u, v \rangle$ be a join-irreducible of $L$. Since

$$\langle u, v \rangle = \langle u, \underline{\iota_1}(u) \rangle \vee_L \langle \underline{\iota_2}(v), v \rangle,$$

either $\underline{\iota_1}(u) = v$ or $\underline{\iota_2}(v) = u$. Suppose that $\underline{\iota_1}(u) = v$. To show that $u$ is join-irreducible in $L_1$, let $u = u_1 \vee u_2$ for some $u_1, u_2 \in L_1$. Since $\underline{\iota_1}(u_1) \le \underline{\iota_1}(u)$ and $\underline{\iota_1}(u_2) \le \underline{\iota_1}(u)$,

$$\langle u, v \rangle = \langle u, \underline{\iota_1}(u) \rangle = \langle u_1, \underline{\iota_1}(u_1) \rangle \vee_L \langle u_2, \underline{\iota_1}(u_2) \rangle.$$

Since $\langle u, v \rangle$ is join-irreducible, this implies that $u_1 = u$ or $u_2 = u$, as desired. Similarly, if $\underline{\iota_2}(v) = u$, $v$ is join-irreducible in $L_2$. ∎

The converse of Theorem 4.4.6 holds if $L$ is a full subdirect product of $L_1$ and $L_2$:

**Theorem 4.4.7** *Let $L$ be a full subdirect product of $L_1$ and $L_2$. Then $\langle u, v \rangle$ is join-irreducible in $L$ iff either $u$ is join-irreducible in $L_1$ and $v = \underline{\iota}_1(u)$ or $v$ is join-irreducible in $L_2$ and $u = \underline{\iota}_2(v)$.*

**Proof.**    The 'if' direction is given by Theorem 4.4.6.

Assume that $u$ is join-irreducible and $v = \underline{\iota}_1(u)$. Let $\langle u_1, v_1 \rangle \in L$, $\langle u_2, v_2 \rangle \in L$ such that

$$\langle u, v \rangle = \langle u_1, v_1 \rangle \vee_L \langle u_2, v_2 \rangle.$$

Since $L$ is a full subdirect product,

$$\langle u_1, v_1 \rangle \vee_L \langle u_2, v_2 \rangle = \langle u_1 \vee u_2, v_1 \vee v_2 \rangle.$$

Hence $u_1 \vee u_2 = u$ and $v_1 \vee v_2 = v$. Since $u$ is join-irreducible, either $u_1 = u$ or $u_2 = u$. Without loss of generality, assume that $u_1 = u$. Since $v = \underline{\iota}_1(u)$ is minimal in $\iota_1(u)$ and $v_1 \leq v$, it follows that $v_1 = v$, so that $\langle u_1, v_1 \rangle = \langle u, v \rangle$, as desired.

The case where $v$ is join-irreducible and $u = \underline{\iota}_2(v)$ is symmetric.    ∎

Theorem 4.4.7 and inductive use of the representation of $L^P$ as a subdirect product of $L$ and $L^{P \setminus \{x\}}$ (for non-maximal $x \in P$) yields the fact that the join-irreducibles of $L^P$ are the maps $f : P \to L$ such that for some join-irreducible $a \in L$ and element $y \in P$, $f(z) = a$ for $z \in [y)$ and $f(z) = \hat{0}$ otherwise.

A subdirect product $L$ of $L_1$ and $L_2$ can be viewed as an internal subdirect product, where $L_1$ and $L_2$ are subsets of $L$ closed under least upper bounds.

**Definition.**    Let $L' \subseteq L$. Then $L'$ is a *lub-subsemilattice* of $L$ iff for every $U \subseteq L'$ such that $\bigvee U$ exists in $L$, $\bigvee U \in L'$ (i.e. $L'$ is closed under least $\underline{u}$pper $\underline{b}$ounds of $L$).

If $L'$ is a lub-subsemilattice of $L$, then $\bigvee \emptyset = \hat{0} \in L'$; and for every $u, v \in L'$,

$$u \wedge_{L'} v = \bigvee \{ w \in L' \,|\, w \leq u \wedge_L v \},$$

so that $L'$ is a meet-semilattice (the meet of $L'$ need not agree with the meet of $L$).

**Observation 4.4.8** *The subset $L'$ of $L$ is a lub-subsemilattice of $L$ iff $\hat{0}_L \in L'$ and for every $u, v \in L'$, if $u \vee_{L'} v = w \in L$, then $u \vee_{L'} v = w \in L'$.*

**Theorem 4.4.9** *Let $L$ be a subdirect product of $L_1$ and $L_2$. Define $L_1' \subseteq L$ and $L_2' \subseteq L$ by*

$$\begin{aligned} L_1' &= \{ \langle u, \underline{\iota}_1(u) \rangle \,|\, u \in L_1 \}, \\ L_2' &= \{ \langle \underline{\iota}_2(v), v \rangle \,|\, v \in L_2 \}. \end{aligned}$$

*Then, for $i = 1, 2$ $L_i'$ is a lub-subsemilattice of $L$ and $\pi_i \dagger L_i' : L_i' \to L_i$ is a meet-isomorphism.*

**Proof.** By symmetry, it suffices to prove the result for $L_1'$. Observe first that

$$\hat{0}_{L_1 \times L_2} = \langle \hat{0}_{L_1}, \underline{\iota_1}(\hat{0}_{L_1}) \rangle \wedge \langle \underline{\iota_2}(\hat{0}_{L_2}), \hat{0}_{L_2} \rangle \in L,$$

so that $\langle \hat{0}_{L_1}, \underline{\iota_1}(\hat{0}_{L_1}) \rangle = \hat{0}_{L_1 \times L_2} = \hat{0}_L \in L_1'$.

Let $u_1, u_2 \in L_1$ and suppose that $\langle u_1, \underline{\iota_1}(u_1) \rangle \vee_L \langle u_2, \underline{\iota_1}(u_2) \rangle$ exists in $L$. Then $u_1 \vee u_2$ exists in $L_1$, and since $\underline{\iota_1}(u_1 \vee u_2) \geq \underline{\iota_1}(u_1)$ and $\underline{\iota_1}(u_1 \vee u_2) \geq \underline{\iota_1}(u_2)$ (Observation 4.4.5), we have

$$\langle u_1, \underline{\iota_1}(u_1) \rangle \vee_L \langle u_2, \underline{\iota_1}(u_2) \rangle = \langle u_1 \vee u_2, \underline{\iota_1}(u_1 \vee u_2) \rangle.$$

Thus $L_1'$ is closed under pairwise greatest lower bounds of $L$. It follows that $L_1'$ is a lub-subsemilattice of $L$. Since

$$\langle u_1, \underline{\iota_1}(u_1) \rangle \wedge_{L_1'} \langle u_2, \underline{\iota_1}(u_2) \rangle = \langle u_1 \wedge u_2, \underline{\iota_1}(u_1 \wedge u_2) \rangle,$$

$\pi_1 \dagger L_1'$ is a meet-isomorphism onto $L_1$, as required. ∎

Let $L$ be a subdirect product of $L_1$ and $L_2$, and let $L_1'$ and $L_2'$ be the lub-subsemilattices of $L$ defined in Theorem 4.4.9. The isomorphisms $L_1 \cong L_1'$ and $L_2 \cong L_2'$ induce projections $\pi_1'$ and $\pi_2'$ of $L$ onto $L_1'$ and $L_2'$ respectively. For $i = 1, 2$ these projections satisfy

$$\pi_i'(u) = \bigvee_L \{ w \in L_i' \mid w \leq u \} \in L_i'$$

and

$$u = \pi_1'(u) \vee_L \pi_2'(u)$$

for $u \in L$.

**Definition.** Let $A \subseteq L$. The lub-subsemilattice of $L$ generated by $A$ is given by

$$\mathrm{G}_L(A) = \{ \bigvee_L B \mid B \subseteq A \}.$$

The *projection* $\pi_A$ of $L$ onto $\mathrm{G}_L(A)$ is defined by

$$\pi_A(u) = \bigvee_L (\, (u] \cap A)$$

for $u \in L$. This gives $\pi_i' = \pi_{L_i'}$ for $i = 1, 2$.

**Theorem 4.4.10** *The map $\pi_A$ is a meet-homomorphism of $L$ onto $\mathrm{G}_L(A)$.*

**Proof.** Let $L' = G_L(A)$. The proper join-irreducibles of $L'$ are included in $A$. This implies that the family $\mathcal{P}(L')_{\cap A}$ is an intersection-closed representation of $L'$. The map $\sigma : L' \to \mathcal{P}(L')_{\cap A}$ defined by $\sigma(u) = (u] \cap A$ is a meet-isomorphism with inverse given by $\sigma^{-1}(B) = \bigvee_L(B)$ for $B \in \mathcal{P}(L')_{\cap A}$. Let $\rho : L \to 2^A$ be the meet-homomorphism defined by $\rho(u) = (u] \cap A$. Then $\pi_A = \sigma^{-1} \circ \rho$ and the result follows. ∎

**Definition.** The semilattice $L$ is the *internal subdirect product* of $L_1$ and $L_2$ iff $L_1$ and $L_2$ are lub-subsemilattices of $L$ such that $L_1 \cup L_2 \supseteq J(L)$.

**Theorem 4.4.11** *If $L$ is the internal subdirect product of $L_1$ and $L_2$, then $L$ is a subdirect product of $L_1$ and $L_2$.*

**Proof.** Let $\sigma : L \to L_1 \times L_2$ be defined by

$$\sigma(u) = \langle \pi_{L_1}(u), \pi_{L_2}(u) \rangle.$$

By Theorem 4.4.10, $\sigma$ is a meet-homomorphism from $L$ into $L_1 \times L_2$. Define $\rho : L_1 \times L_2 \to L$ by

$$\rho(\langle u, v \rangle) = u \vee v \quad \text{if } u \vee v \text{ exists in } L.$$

Then

$$\rho(\sigma(u)) = \bigvee((u] \cap L_1) \vee \bigvee((u] \cap L_2).$$

Since $J(L) \subseteq L_1 \cup L_2$, $\rho(\sigma(u)) = u$. Therefore, $\sigma$ is one-to-one, whence $L$ is isomorphic to a meet-subsemilattice of $L_1 \times L_2$. Since $\pi_{L_1}$ and $\pi_{L_2}$ are onto $L_1$ and $L_2$ respectively, $L$ is a subdirect product of $L_1$ and $L_2$. ∎

Theorems 4.4.6, 4.4.9 and 4.4.11 yield:

**Corollary 4.4.12** *The semilattice $L$ is a subdirect product of $L_1$ and $L_2$ iff $L$ is the internal subdirect product of $L_1' \cong L_1$ and $L_2' \cong L_2$.*

**Corollary 4.4.13** *The semilattice $L$ is a subdirect product of $L_1$ and $L_2$ iff there are subsets $A_1$ and $A_2$ of $L$ such that $A_1 \cup A_2 \supseteq J(L)$ and $\mathcal{P}(L)_{\cap A_i} \cong L_i$, for $i = 1, 2$.*

## 4.5 Constructions which Preserve Matching Properties

Let $L$ be a semilattice. Let $P$ and $Q$ be posets. Let $F_1$ and $F_2$ be filters of $P$ such that $F_1 \supseteq F_2$, and let $G_1$ and $G_2$ be filters of $Q$ such that $G_1 \supseteq G_2$.

**Theorem 4.5.1** *If there are decreasing matchings*

$$\sigma : \mathrm{T}(L^P, F_1, a) \to \mathrm{T}(L^P, F_2, a) \text{ and } \rho : \mathrm{T}(L^Q, G_1, a) \to \mathrm{T}(L^Q, G_2, a),$$

*then there is a decreasing matching*

$$(\sigma + \rho) : \mathrm{T}(L^{P+Q}, F_1 + G_1, a) \to \mathrm{T}(L^{P+Q}, F_2 + G_2, a).$$

**Proof.**   The matching $\sigma + \rho$ is defined by

$$(\sigma + \rho)(f)(x) = \begin{cases} \sigma(f \dagger P)(x) & \text{if } x \in P, \\ \rho(f \dagger Q)(x) & \text{if } x \in Q, \end{cases}$$

for $f \in \mathrm{T}(L^{P+Q}, F_1 + G_1, a)$, where $f \dagger P$ and $f \dagger Q$ are the restrictions of $f$ to $P$ and $Q$ respectively. ■

The following is a partial converse of Theorem 4.5.1:

**Theorem 4.5.2** *Suppose that $G_1 = G_2$ or that $a$ is an atom of $L$. If there is a decreasing matching*

$$\sigma : \mathrm{T}(L^{P+Q}, F_1 + G_1, a) \to \mathrm{T}(L^{P+Q}, F_2 + G_2, a),$$

*then there is a decreasing matching*

$$\sigma' : \mathrm{T}(L^P, F_1, a) \to \mathrm{T}(L^P, F_2, a).$$

**Proof.**   For $f \in L^P$ and $g \in L^Q$, let $f + g \in L^{P+Q}$ be the map defined by

$$(f + g)(x) = \begin{cases} f(x) & \text{if } x \in P, \\ g(x) & \text{if } x \in Q. \end{cases}$$

Let $h \in \mathrm{T}(L^Q, G_1, a)$ be the map defined by

$$h(x) = \begin{cases} a & \text{if } x \in G_1, \\ \hat{0} & \text{otherwise.} \end{cases}$$

Then $h$ is the minimal member of $\mathrm{T}(L^Q, G_1, a)$.

Given the decreasing matching $\sigma : \mathrm{T}(L^{P+Q}, F_1 + G_1, a) \to \mathrm{T}(L^{P+Q}, F_2 + G_2, a)$, define $\sigma'$ by

$$\sigma'(f) = \sigma(f + h) \dagger P,$$

for $f \in \mathrm{T}(L^P, F_1, a)$. The definition implies that $\sigma'(f) \in \mathrm{T}(L^P, F_2, a)$ and $\sigma'(f) \leq f$ for $f \in \mathrm{T}(L^P, F_1, a)$. It remains to show that $\sigma'$ is one-to-one. If $G_1 = G_2$, the only map $h' \in \mathrm{T}(L^Q, G_2, a)$ with $h' \leq h$ is $h' = h$. If $a$ is an atom, the only map $h' \in \mathrm{T}(L^Q, G_2, a)$ with $h' \leq h$ is the map $h'$ defined by $h'(x) = a$ if $x \in G_2$ and $h'(x) = \hat{0}$ otherwise. This implies that $\sigma(f + h) \dagger Q = h'$ for every $f \in \mathrm{T}(L^P, F_1, a)$. Since $\sigma$ is one-to-one, $\sigma'$ is one-to-one. ■

Theorems 4.5.1 and 4.5.2 have the following corollaries:

**Corollary 4.5.3** $(L, a) \in \mathrm{MPf}(P+Q)$ *iff* $(L, a) \in \mathrm{MPf}(P)$ *and* $(L, a) \in \mathrm{MPf}(Q)$.

**Corollary 4.5.4** *If* $(L, a) \in \mathrm{MPt}(P)$ *and* $(L, a) \in \mathrm{MPt}(Q)$, *then* $(L, a) \in \mathrm{MPt}(P + Q)$.

**Corollary 4.5.5** *Let a be an atom of L. Then $(L, a) \in \mathrm{MPt}(P + Q)$ iff $(L, a) \in \mathrm{MPt}(P)$ and $(L, a) \in \mathrm{MPt}(Q)$.*

**Theorem 4.5.6** *Let I be an ideal of L. If there is a decreasing matching*

$$\sigma : \mathrm{T}(L^P, F_1, a) \to \mathrm{T}(L^P, F_2, a),$$

*then there is a decreasing matching*

$$\sigma \dagger I^P : \mathrm{T}(I^P, F_1, a) \to \mathrm{T}(I^P, F_2, a).$$

**Proof.**   Let $\sigma : \mathrm{T}(L^P, F_1, a) \to \mathrm{T}(L^P, F_2, a)$ be a decreasing matching. The decreasing matching $\sigma \dagger I^P : \mathrm{T}(I^P, F_1, a) \to \mathrm{T}(I^P, F_2, a)$ is obtained by restricting the domain of $\sigma$ to $\mathrm{T}(I^P, F_1, a)$. Since $I$ is an ideal and $\sigma$ is decreasing, the range of the restriction is included in $\mathrm{T}(I^P, F_2, a)$. ∎

**Corollary 4.5.7** *Let I be an ideal of L. If $(L, a) \in \mathrm{MPf}(P)$, then $(I, a) \in \mathrm{MPf}(P)$. If $(L, a) \in \mathrm{MPt}(P)$, then $(I, a) \in \mathrm{MPt}(P)$.*

Let $L_1$ and $L_2$ be semilattices. The following theorem generalizes part of Theorem 4.2.3.

**Theorem 4.5.8** *Let $f : L_1 \to \mathrm{Sub}(L_2)$ be a meet-concave map such that for every $u, v \in L_1$ with $u \leq v$, $f(u) \supseteq f(v)$. Let a be a join-irreducible of $L_1$. If there is a decreasing matching $\sigma : \mathrm{T}(L_1^P, F_1, a) \to \mathrm{T}(L_1^P, F_2, a)$, then there is a decreasing matching*

$$\sigma' : \mathrm{T}(\, (L_1 \bowtie_f L_2)^P, F_1, \langle a, \underline{\iota}_1(a) \rangle \,) \to \mathrm{T}(\, (L_1 \bowtie_f L_2)^P, F_2, \langle a, \underline{\iota}_1(a) \rangle \,).$$

**Proof.**   Let $L = L_1 \bowtie_f L_2$. To check that $\langle a, \underline{\iota}_1(a) \rangle$ is join-irreducible in $L$, suppose

$$\langle a, \underline{\iota}_1(a) \rangle = \langle u_1, v_1 \rangle \vee_L \langle u_2, v_2 \rangle.$$

Then $u_1 \vee u_2 \leq a$. Using $\underline{\iota}_1(a) \in f(a) \subseteq f(u_1 \vee u_2)$, we obtain $\langle u_1 \vee u_2, \underline{\iota}_1(a) \rangle \in L$. Since $\underline{\iota}_1(a) \geq v_1 \vee v_2$,

$$\langle u_1 \vee u_2, \underline{\iota}_1(a) \rangle \geq \langle u_1, v_1 \rangle \vee_L \langle u_2, v_2 \rangle = \langle a, \underline{\iota}_1(a) \rangle,$$

which yields $u_1 \vee u_2 = a$. Since $a$ is join-irreducible, either $u_1 = a$ or $u_2 = a$. Without loss of generality, assume $u_1 = a$. Then, by definition of $\underline{\iota}_1$, $v_1 \geq \underline{\iota}_1(a)$, so that $\langle u_1, v_1 \rangle = \langle a, \underline{\iota}_1(a) \rangle$, as required.

Let $\sigma : \mathrm{T}(L_1^P, F_1, a) \to \mathrm{T}(L_1^P, F_2, a)$ be a decreasing matching. For $g \in \mathrm{T}(L^P, F_1, \langle a, \underline{\iota}_1(a) \rangle)$, define $\sigma'(g)$ by

$$\sigma'(g)(x) = \langle\, \sigma(\pi_1 \circ g)(x), \pi_2(g(x)) \,\rangle.$$

Since $\sigma$ is decreasing, $\sigma(\pi_1 \circ g)(x) \leq \pi_1(g(x))$. Using $\langle \pi_1(g(x)), \pi_2(g(x)) \rangle \in L$ and the assumption on $f$, we get $\langle \sigma(\pi_1 \circ g)(x), \pi_2(g(x)) \rangle \in L$. Since $\sigma(\pi_1 \circ g)$ and $\pi_2 \circ g$ are order-preserving functions, so is $\sigma'(g)$. Since $\sigma(\pi_1 \circ g)(x) \geq a$ iff $x$ is in $F_2$, it follows that $\sigma'$ maps $\mathrm{T}(L^P, F_1, \langle a, \iota_1(a) \rangle)$ into $\mathrm{T}(L^P, F_2, \langle a, \iota_1(a) \rangle)$. Since $\sigma$ is decreasing, so is $\sigma'$. Using the assumption that $\sigma$ is one-to-one, the map $g$ can be recovered from $\sigma'(g)$ by

$$g = \langle \sigma^{-1}(\pi_1 \circ \sigma'(g)), \pi_2 \circ \sigma'(g) \rangle.$$

Hence $\sigma'$ is one-to-one and it follows that $\sigma'$ is the required matching. ∎

**Corollary 4.5.9** *Let $f : L_1 \to \mathrm{Sub}(L_2)$ be a meet-concave map such that for every $u, v \in L_1$ with $u \leq v$, $f(u) \supseteq f(v)$. If $L_1 \in \mathrm{MPf}(P)$, then $L_1 \bowtie_f L_2 \in \mathrm{MPf}(P)$. If $L_1 \in \mathrm{MPt}(P)$, then $L_1 \bowtie_f L_2 \in \mathrm{MPt}(P)$.*

**Example 4.5.10** Let $x$ be a minimal member of $Q$ and $Q' = Q \setminus \{x\}$. By Example 4.4.2, $L^Q$ is a subdirect product of $L$ and $L^{Q'}$. The associated function $\iota_1$ satisfies that $\iota_1(u) \supseteq \iota_1(v)$ if $u \leq v$.

**Corollary 4.5.11** *Let $Q$ be a poset. If $L \in \mathrm{MPf}(P)$ then $L^Q \in \mathrm{MPf}(P)$. If $L \in \mathrm{MPt}(P)$ then $L^Q \in \mathrm{MPt}(P)$.*

This generalizes the part of Theorem 4.2.5 which shows that every distributive lattice, i.e. every lattice of the form $2^Q$, has the full matching property.

## 4.6 Neighborhoods in Lattices

Though we now focus on lattices, much of what will be shown can be applied to semilattices by completion and use of Theorem 4.5.6.

Let $L$ be a lattice. By replacing $L$ with the canonical union-closed representation $\overline{\mathcal{F}}(L)$ of $L$ if necessary (Section 4.3), we can assume that $L$ is a union-closed family of sets such that $\emptyset \in L$ and $L$ is generated by the join-irreducibles of $L$. We then have $\hat{0}_L = \emptyset$ and the domain of $L$ is given by $\hat{1}_L = \bigcup L$.

**Definition.** Let $U \subseteq \bigcup L$. The *lattice neighborhood* in $L$ of $U$, denoted by $\mathcal{N}_L(U)$, is the join-subsemilattice of $L$ generated by $\emptyset = \hat{0}$ and the generators $a$ of $L$ with $a \cap U \neq \emptyset$.

If $U \subseteq \bigcup L$, $u \in L$ and $x \in u \cap U$, then there is a generator $a$ of $L$ such that $x \in a \subseteq u$. By the definition of lattice neighborhood, $a \in \mathcal{N}_L(U)$. This yields the following observation:

**Observation 4.6.1** *If $u \in L$ and $U \subseteq \bigcup L$, then $\pi_{\mathcal{N}_L(U)}(u) \cap U = u \cap U$.*

The map $\pi_{\mathcal{N}_L(U)}$ is the projection of $L$ onto $\mathcal{N}_L(U)$ defined in Section 4.4.

A sequence of lattice neighborhoods of $U$ which includes increasingly more distant generators of $L$ is obtained by iterating $\mathcal{N}_L$:

**Definition.** Let $\mathcal{N}_L^1(U) = \mathcal{N}_L(U)$. For $n \geq 1$, define

$$\mathcal{N}_L^{2n+1}(U) = \mathcal{N}_L(\bigcup \mathcal{N}_L^{2n-1}(U)).$$

Intermediate lattice neighborhoods are given by the ideals generated by $\mathcal{N}_L^{2n+1}(U)$. For $n \geq 1$, define

$$\mathcal{N}_L^{2n}(U) = \left(\mathcal{N}_L^{2n-1}(U)\right] = \{V \in L \mid V \subseteq \bigcup \mathcal{N}_L^{2n-1}(U)\}.$$

Let $\mathcal{N}_L^0(U) = \{u \in L \mid u \subseteq U\}$ (the ideal of $L$ generated by $U$).

Lattice neighborhoods can be used to characterize subdirect products satisfying the conditions of Theorem 4.5.8.

**Theorem 4.6.2** *Let $L$ and $L_1$ be lattices. The following are equivalent:*

*(i) There exists a semilattice $L_2$ and a meet-concave map $f : L_1 \to \mathrm{Sub}(L_2)$ such that $L \cong L_1 \bowtie_f L_2$ and for every $u, v \in L$ with $u \subseteq v$, $f(u) \supseteq f(v)$.*

*(ii) There is a subset $U$ of the domain of $L$ such that $L_1 \cong \mathcal{N}_L(U)$ and for every $u, v \in \mathcal{N}_L(U)$, $u \cap U = v \cap U$ iff $u = v$.*

**Proof.** Assume (i). By Corollary 4.4.12 we can assume that $L$ is the internal subdirect product of $L_1$ and $L_2$. Let $\pi_1 = \pi_{L_1}$ and $\pi_2 = \pi_{L_2}$ be the projections. The map $f$ is given by

$$f(u) = \{\pi_2(w) \mid w \in L \text{ and } \pi_1(w) = u\}.$$

Let $U = \bigcup L \setminus \hat{1}_{L_2}$. Suppose that $u, v \in L_1$ and $u \cap U = v \cap U$. Let $u' = \pi_1(u \cup \hat{1}_{L_2}) \in L_1$. Then $u \subseteq u'$. Since $\pi_2(u \cup \hat{1}_{L_2}) = \hat{1}_{L_2}$, the assumption on $f$ implies that there exists $w \in L$ such that $\pi_1(w) = u$ and $\pi_2(w) = \hat{1}_{L_2}$. This yields $w = \pi_1(w) \vee \pi_2(w) = u \vee \hat{1}_{L_2}$, hence $u = \pi_1(w) = \pi_1(u \vee \hat{1}_{L_2}) = u'$. Similarly, using $u \cup \hat{1}_{L_2} = v \cup \hat{1}_{L_2}$, we obtain $v = u'$, so that $u = v$. Note that this implies that $(L_1)_{\cap U} \cong L_1$. To show that $L_1 = \mathcal{N}_L(U)$, suppose that $a$ is a join-irreducible of $L$ with $a \cap U \neq \emptyset$. Since $\pi_2(a) \subseteq \hat{1}_{L_2}$, $\pi_2(a) \subset a$. Since $\pi_1(a) \vee \pi_2(a) = a \in J(L)$, $\pi_1(a) = a$, so that $a \in L_1$. By arbitrariness of $a$, $L_1 \supseteq \mathcal{N}_L(U)$. Let $\sigma : \mathcal{N}_L(U) \to (L_1)_{\cap U}$ be the map defined by $\sigma(u) = u \cap U$. Since $\pi_{\mathcal{N}_L(U)}(u) \cap U = u \cap U$ for every $u \in L$, it follows that $\sigma \circ (\pi_{\mathcal{N}_L} \dagger L_1)$ is onto $(L_1)_{\cap U}$. This implies that $\sigma$ is onto $(L_1)_{\cap U}$, so that $\left|\mathcal{N}_L(U)\right| \geq \left|(L_1)_{\cap U}\right| = \left|L_1\right|$. The result follows.

For the converse, assume (ii). Let $L_2$ be the lub-subsemilattice of $L$ generated by the join-irreducibles of $L$ disjoint from $U$. Then $L$ is the internal subdirect product of $\mathcal{N}_L(U)$ and $L_2$. Let $\pi_1 = \pi_{\mathcal{N}_L(U)}$ and $\pi_2 = \pi_{L_2}$. We have to show that if $u, u' \in \mathcal{N}_L(U)$, $w \in L$, $\pi_1(w) = u$, $\pi_2(w) = v$ and $u' \subseteq u$, then there exists a $w' \in L$ such that $\pi_1(w') = u'$ and $\pi_2(w') = v$. Such a $w'$ is given by $w' = (w \setminus U) \cup u'$. Since $w' \setminus U = w \setminus U$, $\pi_2(w') = \pi_2(w) = v$. We have $\pi_1(w') \cap U = w' \cap U = u' \cap U$, so that by assumption on $\mathcal{N}_L(U)$, $\pi_1(w') = u'$, as required. ∎

**Observation 4.6.3** *Let $U \subseteq \bigcup L$. The following are equivalent:*

   *(i) For every $u, v \in \mathcal{N}_L(U)$, $u \cap U = v \cap U$ iff $u = v$.*

   *(ii) $\mathcal{N}_L(U)_{\cap U}$ is a union-closed representation of $\mathcal{N}_L(U)$.*

   *(iii) There exists $V \subseteq U$ such that $\mathcal{N}_L(V)_{\cap V}$ is an irredundant union-closed representation of $\mathcal{N}_L(V)$ and $\mathcal{N}_L(V) = \mathcal{N}_L(U)$.*

This implies that every subdirect product of lattices $L_1$ and $L_2$ satisfying condition (i) of Theorem 4.6.2 can be obtained by the following construction. Start with the canonical union-closed representation $\mathcal{F} = \overline{\mathcal{F}}(L_1)$ of $L_1$ with domain $U = M(L_1)$ and an arbitrary union-closed representation $\mathcal{H}$ of $L_2$ on a disjoint domain $V$. Let $\mathcal{G}$ be the (join-irreducible) generators of $\mathcal{F}$. Modify $\mathcal{G}$ by adjoining to each non-empty member $u$ of $\mathcal{G}$ a subset $V(u)$ of $V$:

$$\mathcal{G}' = \{u \cup V(u) \mid u \in \mathcal{G}\}.$$

Let $\mathcal{F}'$ be the union-closed family of sets generated by $\mathcal{G}'$. The union-closed family of sets $\mathcal{F}' \vee \mathcal{H}$ is then a subdirect product of the desired kind.

     The conditions of Observation 4.6.3 are satisfied if for some element $x$ of the domain of $L$, $x$ is contained in exactly one generator $a$ of $L$. Then the lattice neighborhood of $\{x\}$ is the two-element chain $\{\emptyset, a\}$. Theorems 4.5.8 and 4.6.2 and the fact that $[2] \in \mathrm{MPf}(P)$ imply that $L \in \mathrm{MPf}(P)$ for every poset $P$. This can be used to prove the following theorem.

**Theorem 4.6.4** *Let $h$ be the height of $L$. If $h = \big| J(L) \big|$, then $L \in \mathrm{MPf}(P)$ for every poset $P$.*

**Proof.** Let $u_0 = \hat{0} < u_1 <, \ldots, < u_h = \hat{1}$ be a saturated chain of $L$. For each $i > 0$, there is a (non-empty) generator $a_i \subseteq u_i$ such that $a_i \nsubseteq u_{i-1}$. Since $h = \big| J(L) \big|$, $\{a_1, \ldots, a_h\} = J(L)$. By construction, there exists $x \in a_h \setminus (\bigcup_{i=1}^{h-1} a_i)$. Thus the only generator of $L$ which contains $x$ is $a_h$. The proof is completed by the discussion preceding the statement of the theorem. ∎

     Note that the proof of Theorem 4.6.4 shows that for every lattice $L$ of height $h$, $h \leq \big| J(L) \big|$.

## 4.7  Lattices with Special Lattice Neighborhoods
     Let $L$ be a lattice and $P$ a poset.

**Definition.** The decreasing matching $\sigma : \mathrm{T}(L^P, F, a) \to \mathrm{T}(L^P, G, a)$ is $a$-*invertible* iff for every function $f \in \mathrm{T}(L^P, F, a)$,

$$f(x) = \begin{cases} \sigma(f)(x) \vee a & \text{if } x \in F, \\ \sigma(f)(x) & \text{otherwise.} \end{cases}$$

**Observation 4.7.1** *If* $\sigma : \mathrm{T}(L^P, F, a) \to \mathrm{T}(L^P, G, a)$ *and* $\rho : \mathrm{T}(L^P, G, a) \to \mathrm{T}(L^P, H, a)$ *are $a$-invertible decreasing matchings, then so is* $\rho \circ \sigma$.

Many of the commonly studied lattices have matchings which are $a$-invertible. (see Theorems 4.7.4 and 4.7.9 below). The existence of such matchings is completely determined by a neighborhood of $a$.

**Theorem 4.7.2** *Let* $a \in J(L)$ *and let* $\sigma : \mathrm{T}(L^P, F, a) \to \mathrm{T}(L^P, G, a)$ *be an $a$-invertible decreasing matching. Then the restriction of $\sigma$ to* $\mathrm{T}((\mathcal{N}_L^2(a))^P, F, a)$ *is an $a$-invertible decreasing matching into* $\mathrm{T}((\mathcal{N}_L^2(a))^P, G, a)$.

**Proof.** Let $I = \mathcal{N}_L^2(a)$. Then $I$ is the ideal of $L$ generated by the union of the generators of $L$ which intersect $a$. As in Theorem 4.5.6, the matching $\sigma$ restricts to a decreasing matching $\sigma' : \mathrm{T}(I^P, F, a) \to \mathrm{T}(I^P, G, a)$. If $\sigma$ is $a$-invertible, then so is $\sigma'$. ∎

The converse also holds:

**Theorem 4.7.3** *Let* $a \in J(L)$ *and let* $L'$ *be a join-subsemilattice of $L$ such that* $L' \supseteq \mathcal{N}_L(a)$. *If* $\sigma' : \mathrm{T}(L'^P, F, a) \to \mathrm{T}(L'^P, G, a)$ *is an $a$-invertible decreasing matching, then there is an $a$-invertible decreasing matching* $\sigma : \mathrm{T}(L^P, F, a) \to \mathrm{T}(L^P, G, a)$.

**Proof.** Let $B = \{b \in J(L) \mid b \cap a = \emptyset\}$. Define $\sigma$ by

$$\sigma(f)(x) = \pi_B(f(x)) \cup \sigma'(\pi_{L'} \circ f)(x)$$

for $f \in \mathrm{T}(L^P, F, a)$ and $x \in P$.

Let $f \in \mathrm{T}(L^P, F, a)$. Since $a \in \mathcal{N}_L(a) \subseteq L'$, $(\pi_{L'} \circ f)(x) \supseteq a$ iff $f(x) \supseteq a$. This implies that $\pi_{L'} \circ f \in \mathrm{T}(L'^P, F, a)$, so that $\sigma(f)$ is well-defined. Since $f$ is order-preserving, so are $\pi_B \circ f$ and $\sigma'(\pi_{L'} \circ f)$, whence $\sigma(f)$ is order-preserving. Since $\pi_B(f(x)) \cap a = \emptyset$, $\sigma(f)(x) \supseteq a$ iff $\sigma'(\pi_{L'} \circ f)(x) \supseteq a$, which implies that $\sigma(f) \in \mathrm{T}(L^P, G, a)$. Since $\sigma'(\pi_{L'} \circ f)(x) \subseteq (\pi_{L'} \circ f)(x)$, we have $\sigma(f)(x) \subseteq f(x)$ for each $x \in P$, so that $\sigma$ is a decreasing map from $\mathrm{T}(L^P, F, a)$ into $\mathrm{T}(L^P, G, a)$. It remains to show that $\sigma$ is $a$-invertible. Let $x \in F$. Since $\sigma'$ is $a$-invertible, $\sigma'(\pi_{L'} \circ f)(x) \cup a = (\pi_{L'} \circ f)(x)$. Therefore, using the fact that $L$ is the internal subdirect product of $\mathrm{G}_L(B)$ and $L'$,

$$\sigma(f)(x) \cup a = \pi_B(f(x)) \cup \pi_{L'}(f(x)) = f(x)$$

as required. ∎

**Definition.** Let $\mathcal{C}$ be a class of lattices and $U \subseteq \bigcup L$. Then $L'$ is a *$\mathcal{C}$-neighborhood* of $U$ iff $L'$ is a join-subsemilattice of $L$ such that $L' \supseteq \mathcal{N}_L(u)$ and $L' \in \mathcal{C}$. If $\mathcal{C}$ is the class of lower semimodular (geometric, etc.) lattices, a $\mathcal{C}$-neighborhood is called a lower semimodular (geometric, etc.) neighborhood.

**Definition.** The member $u$ of $L$ is *lower semimodular* iff for every $v, w \in L$ such that $v$ covers $u$, either $v \wedge w$ covers $u \wedge w$ or $v \wedge w = u \wedge w$ (see Figure 4.3).

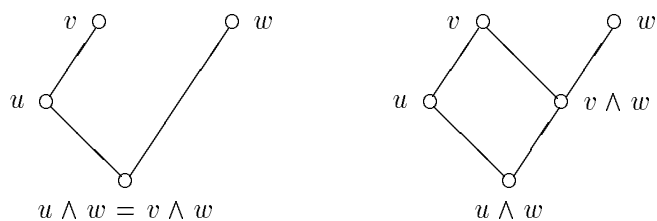

Figure 4.3

Note that $L$ is lower semimodular iff every $u \in L$ is lower semimodular.

**Theorem 4.7.4** *If $L$ has a lower semimodular coatom $c$, then for every $a \in J(L)$ with $a \not\leq c$, $(L, a) \in \mathrm{MPt}(P)$ where the matchings are $a$-invertible.*

Observe that since $\hat{1} = \bigvee J(L)$, there is at least one join-irreducible $a \in L$ such that $a \not\leq c$.

**Proof.** Let $\eta : [a) \to L \setminus [a)$ be the map given by $\eta(u) = u \wedge c$ for $u \in [a)$.

**Lemma 4.7.5** *The map $\eta$ is an order-preserving, decreasing matching and for every $u \in [a)$, $\eta(u) \vee a = u$.*

**Proof.** The fact that $\eta$ is decreasing and order-preserving follows by definition. The element $\hat{1}$ covers $c$, hence, by lower semimodularity of $c$ and since $u \wedge c \neq u$, $u = u \wedge \hat{1}$ covers $\eta(u) = u \wedge c$. We have $\eta(u) \leq \eta(u) \vee a \leq u$ and since $\eta(u) \notin [a)$, $\eta(u) \vee a \neq \eta(u)$. Therefore $\eta(u) \vee a = u$. This also shows that $\eta$ is one-to-one. ∎

Let $F$ be a filter of $P$. Then $\eta$ induces a decreasing matching $\sigma : \mathrm{T}(L^P, P, a) \to \mathrm{T}(L^P, F, a)$ defined by

$$\sigma(f)(x) = \begin{cases} f(x) & \text{if } x \in F, \\ \eta(f(x)) & \text{otherwise,} \end{cases}$$

for $f \in \mathrm{T}(L^P, F, a)$ and $x \in P$. The properties of $\eta$ imply that if $f$ is order-preserving, then so is $\sigma(f)$. The identity $\eta(u) \vee a = u$ implies that $\sigma$ is $a$-invertible. ∎

**Corollary 4.7.6** *If $L$ is a non-trivial lower semimodular lattice, then $L \in \mathrm{MPt}(P)$.*

**Proof.** If $c$ is a coatom of $L$, then $c$ is lower semimodular and there exists a proper join-irreducible $a \not\leq c$. ∎

**Theorem 4.7.7** *If the dual of $L$ is a geometric lattice and $a \in J(L)$, then $(L, a) \in \mathrm{MPf}([n])$ where the matchings are $a$-invertible.*

**Proof.** The dual of $L$ is geometric iff $L$ is coatomic and lower semimodular. Every $u \in L$ is the meet of the coatoms above $u$. This implies that for every $u \notin [a)$, there exists a coatom $\gamma(u) \geq u$ with $\gamma(u) \notin [a)$.

To prove the theorem, it suffices to construct $a$-invertible decreasing matchings $\sigma_k : \mathrm{T}(L^{[n]}, [k, n], a) \longrightarrow \mathrm{T}(L^{[n]}, [k+1, n], a)$ for $1 \leq k \leq n$. The other matchings are obtained by composition.

Let $f \in \mathrm{T}(L^{[n]}, [k, n], a)$. Let $f(0) = \hat{0}$. Define $\sigma_k(f)$ by

$$\sigma_k(f)(i) = \begin{cases} f(k) \wedge \gamma(f(k-1)) & \text{if } i = k, \\ f(i) & \text{otherwise,} \end{cases}$$

for $i \in [n]$. Then $\sigma_k(f) \in \mathrm{T}(L^{[n]}, [k+1, n], a)$ and $\sigma_k(f) \leq f$ (in $L^{[n]}$). By Lemma 4.7.5 and since $\sigma_k(f)(k-1) = f(k-1)$, $f$ can be recovered from $\sigma_k(f)$ by

$$f(i) = \begin{cases} \sigma_k(f)(k) \vee \gamma(\sigma_k(f)(k-1)) & \text{if } i = k, \\ \sigma_k(f)(i) & \text{otherwise,} \end{cases}$$

for $i \in [n]$. Thus $\sigma_k$ is as required. ∎

Theorems 4.7.3 and 4.7.7 yield:

**Corollary 4.7.8** *If $a \in J(L)$ has a dually geometric neighborhood, then $(L, a) \in \mathrm{MPf}([n])$.*

**Theorem 4.7.9** *If $L$ is a geometric lattice and $a$ is an atom of $L$, then $(L, a) \in \mathrm{MPf}([n])$ where the matchings are $a$-invertible.*

**Proof.** As in the proof of Theorem 4.7.7, it suffices to construct $a$-invertible decreasing matchings $\sigma_k : \mathrm{T}(L^{[n]}, [k, n], a) \to \mathrm{T}(L^{[n]}, [k+1, n], a)$ for $1 \leq k \leq n$.

Let $f \in \mathrm{T}(L^{[n]}, [k, n], a)$. Let $f(0) = \hat{0}$. Since $f(k-1) \not\geq a$, $f(k-1) < f(k)$. Let $A$ be an independent set of atoms such that $\bigvee A = f(k-1)$ (see Theorem 2.5.3). Then $A \cup \{a\}$ is independent. Let $B$ be an independent set of atoms such that $B \supseteq A \cup \{a\}$ and $\bigvee B = f(k)$. Let $B' = B \setminus \{a\}$. Then $\bigvee B' \geq f(k)$ and $\bigvee B' \not\geq a$. Furthermore, $\bigvee B' \vee a = f(k)$. Define $\sigma_k(f)$ by

$$\sigma_k(f)(i) = \begin{cases} \bigvee B' & \text{if } i = k, \\ f(i) & \text{otherwise,} \end{cases}$$

for $i \in [n]$. Then $\sigma_k(f) \in \mathrm{T}(L^{[n]}, [k+1, n], a)$ and $\sigma_k(f) \leq f$ (in $L^{[n]}$). The map $f$ can be recovered from $\sigma_k(f)$ by

$$f(i) = \begin{cases} \sigma_k(f)(k) \vee a & \text{if } i = k, \\ \sigma_k(f)(i) & \text{otherwise,} \end{cases}$$

for $i \in [n]$. Therefore, $\sigma_k$ is as required. ∎

Theorems 4.7.3 and 4.7.7 yield:

**Corollary 4.7.10** *If $a \in J(L)$ has a geometric neighborhood, then $(L, a) \in$ MPf([n]).*

## 4.8 The Density Perspective

Let $P$ and $Q$ be posets.

**Definition.** The *P-size* of $Q$ is the cardinality of $Q^P$.

The notion of $P$-size is a generalization of cardinality. The [1]-size of $Q$ is the cardinality of $Q$ and the $[n]$-size of $Q$ is the number of multichains of length $n - 1$ of $Q$.

If $S$ is a subset of $T$, then the *density* in $T$ of $S$ is the ratio $\frac{|S|}{|T|}$ of the cardinalities of $S$ and $T$. More generally, if we have a size measure $w$ for the members of a class of structures and $S$ is a substructure of the structure $T$ in the class, then the *w-density* in $T$ of $S$ is the ratio $\frac{w(S)}{w(T)}$. We formalize this for posets and $P$-sizes.

**Definition.** If $R \subseteq Q$, then the *P-density* in $Q$ of $R$ is the ratio $\frac{|R^P|}{|Q^P|}$ of the $P$-sizes of $R$ and $Q$. If $x \in Q$, then the $P$-density in $Q$ of $x$ is the $P$-density of the principal filter $[x]$ generated by $x$. If $P$ is the one-element poset [1], the prefix '$P$-' is dropped.

Let $L$ be a semilattice. Let $p = \left| [2]^P \right|$. Thus $p$ is the number of filters of $P$.

If $L$ has one of the (downward) matching properties for $P$, then there is a join-irreducible $a \in L$ such that the $P$-size of $[a]$ is at most $\frac{1}{p}$ times the $P$-size of $L$. In other words, there is a join-irreducible $a$ with $P$-density in $L$ at most $\frac{1}{p}$.

If the join-irreducible $a$ of $L$ has density in $L$ at most $\frac{1}{2}$, then $(L, a) \in$ MPw([1]). The corresponding statement for posets other than [1] may not hold.

**Definition.** The pair $(L, a)$ has the *P-density property* iff $a$ is a join-irreducible of $L$ with $P$-density in $L$ at most $\frac{1}{p}$. The semilattice $L$ has the *P-density property* iff there is a join-irreducible $a \in L$ such that $(L, a)$ has the $P$-density property.

Whether all non-trivial (semi-)lattices have the $P$-density property is unkown.

The problem of finding bounds for the minimum $P$-density of the filters generated by join-irreducibles can also be pursued from the following perspective. Let $\mathcal{C}$ be a class of semilattices. Define $\min(\emptyset) = \infty$ and $\max(\emptyset) = 0$.

**Definition.** For $n \geq 0$, define

$$\underline{h}(\mathcal{C}, P, n) = \min\{ \left| L^P \right| \mid L \in \mathcal{C} \text{ and for every } a \in J(L) \cup \{\hat{0}\}, \left| [a)^P \right| \geq n \}.$$

Thus $\underline{h}(\mathcal{C}, P, n)$ is the minimum $P$-size of a semilattice in $\mathcal{C}$ with the property that the $P$-size of each principal filter generated by a join-irreducible is at least $n$. For $n \geq 0$, define

$$\overline{h}(\mathcal{C}, P, n) = \max\{ \left| L^P \right| \mid L \in \mathcal{C} \text{ and for every } a \in J(L) \cup \{\hat{0}\}, \left| L^P \setminus [a)^P \right| \leq n \}.$$

Thus $\overline{h}(\mathcal{C}, P, n)$ is the maximum $P$-size of a semilattice $L$ in $\mathcal{C}$ with the property that for every join-irreducible $a \in L$, the number of order-preserving maps $f : P \to L$ of type other than $(P, a)$ is at most $n$.

The functions $\underline{h}$ and $\overline{h}$ are related by

$$\begin{array}{llllll}
\text{if} & \underline{h}(n) & = & n + m, & \text{then} & \overline{h}(m) \geq n + m, \\
\text{if} & \overline{h}(m) & = & n + m, & \text{then} & \underline{h}(n) \leq n + m,
\end{array}$$

where the arguments $\mathcal{C}$ and $P$ have been omitted. Since $\underline{h}$ and $\overline{h}$ are increasing in the last argument, this implies that:

$$\begin{array}{lllll}
\max\{\underline{h}(n) \mid \underline{h}(n) \leq n + m\} & \leq & \overline{h}(m) & \leq & m + \max\{n \mid \underline{h}(n) \leq n + m\}, \\
n + \min\{m \mid \overline{h}(m) \geq n + m\} & \leq & \underline{h}(n) & \leq & \min\{\overline{h}(m) \mid \overline{h}(m) \geq n + m\}.
\end{array}$$

This does not completely determine one function in terms of the other.

**Observation 4.8.1** *The following are equivalent:*

*(i) Every non-trivial semilattice in $\mathcal{C}$ has the $P$-density property.*

*(ii) For all $n \geq 2$, $\underline{h}(\mathcal{C}, P, n) \geq pn$.*

*(iii) For all $n \geq 1$, $\overline{h}(\mathcal{C}, P, n) \leq (1 + \frac{1}{p-1})n$.*

To prove that (ii) implies (i), observe that if $L$ is a non-trivial semilattice then $\left| L^P \right| \geq p$. This implies that every non-trivial semilattice $L$ with a join-irreducible $a$ such that $\left| [a)^P \right| = 1$ (i.e. $a$ is maximal) has the $P$-density property.

Let $\mathcal{L}$ be the class of all semilattices. Since every semilattice has at least one element, $\underline{h}(\mathcal{L}, P, 0) = \underline{h}(\mathcal{L}, P, 1) = 1$. If $L$ is a non-trivial semilattice and $a \in J(L)$, then $\left| L^P \setminus [a)^P \right| \geq p - 1$ (count the order-preserving maps $f : P \to \{\hat{0}, a\}$). This implies that $\overline{h}(\mathcal{L}, P, n) = 1$ for $0 \leq n < p - 1$.

If $L$ is a semilattice with $|L| \geq 3$, then $L$ has either at least two atoms or a non-atomic join-irreducible. In either case, there is a join-irreducible $a$ with an element $b \in L \setminus [a)$ such that $b \neq \hat{0}$. The set $L^P \setminus [a)^P$ includes $p - 2$ maps $f : P \to L$ onto $\{\hat{0}, a\}$ and $p$ maps $f : P \to L$ into $\{\hat{0}, b\}$. Thus $\left| L^P \setminus [a)^P \right| \geq 2p - 2$. Using the two-element lattice as an example, we obtain $\overline{h}(\mathcal{L}, P, n) = p$ for $p - 1 \leq n < 2p - 2$.

If for every join-irreducible $a \in L$, $\left| [a)^P \right| \geq 2$, then no maximal member of $L$ is join-irreducible, so that $L$ contains a subset (order-)isomorphic to the four-element Boolean lattice $B_2$. This implies that $\left| L^P \right| \geq p^2$. Using $B_2$ as an example, we obtain $\underline{h}(\mathcal{L}, P, n) = p^2$ for $2 \leq n \leq p$.

In summary:

**Observation 4.8.2** *For $0 \leq n < p - 1$,*

$$\underline{h}(\mathcal{L}, P, 0) = \underline{h}(\mathcal{L}, P, 1) = \overline{h}(\mathcal{L}, P, n) = 1.$$

*For $2 \leq n \leq p$,*

$$\underline{h}(\mathcal{L}, P, n) = p^2$$

*and for $p - 1 \leq n < 2p - 2$,*

$$\overline{h}(\mathcal{L}, P, n) = p.$$

The best known general bounds on $\overline{h}$ are given in Section 4.10.

Before continuing the study of the $P$-density property and the functions $\underline{h}$ and $\overline{h}$, consider the function $\underline{g}$ defined as follows.

**Definition.** For $n \geq 0$, define

$$\underline{g}(\mathcal{C}, P, n) = \min\{\left| L^P \right| \mid L \in \mathcal{C} \text{ and for every } a \in J(L), \left| \mathrm{T}(L^P, \emptyset, a) \right| \geq n\}.$$

Thus $\underline{g}(\mathcal{C}, P, n)$ is the minimum $P$-size of a semilattice $L$ in $\mathcal{C}$ with the property that for every proper join-irreducible $a \in L$, the number of maps of type $(\emptyset, a)$ in $L^P$ is at least $n$.

The values of $\underline{g}$ for the class of semilattices are more easily obtained than the values of $\underline{h}$.

**Theorem 4.8.3** *Let $n \geq 1$. If $P$ has a greatest element $\hat{1}$, then*

$$\underline{g}(\mathcal{L}, P, n) = \left( \left\lceil \frac{n-1}{p-1} \right\rceil + 1 \right)(p-1) + 1.$$

**Proof.** Let $m = \lceil \frac{n-1}{p-1} \rceil$. To show that $\underline{g}(\mathcal{L}, P, n) \le (m+1)(p-1)+1$, consider the semilattices $M_k$ of Example 4.2.6. Since $P$ has a greatest member, the $P$-size of $M_k$ is $k(p-1)+1$. Let $a$ be an atom of $M_{m+1}$. Then the complement of $[a)$ is isomorphic to $M_m$. Therefore, the number of maps $f \in M_{m+1}^P$ of type $(\emptyset, a)$ is $m(p-1)+1 \ge n$. The $P$-size of $M_{m+1}$ is $(m+1)(p-1)+1$, which is the right-hand side of the equality in the statement of the theorem. Observe that if $k < m+1$, then $\left| T(M_k^P, \emptyset, a) \right| \le (m-1)(p-1)+1 < n$ for every atom $a$ of $M_k$.

To show that $\underline{g}(\mathcal{L}, P, n) \ge (m+1)(p-1)+1$, suppose that the semilattice $L$ has at least $n$ maps $f \in L^P$ of type $(\emptyset, a)$ for every $a \in J(L)$. Suppose that $L$ has a non-maximal atom $a$ and let $b > a$. There are $p-1$ order-preserving maps $f : P \to \{\hat{0}, b\}$ such that $b$ is in the range of $f$ and $p-1$ order-preserving maps $f : P \to \{\hat{0}, a\}$ such that $a$ is in the range of $f$. Therefore the number of maps $f \in L^P$ not of type $(\emptyset, a)$ is at least $2(p-1)$. This implies that unless $L$ is isomorphic to $M_k$ for some $k$, the $P$-size of $L$ is at least $n + 2(p-1)$. Since $n + 2(p-1) > (m+1)(p-1)+1$, the theorem follows. ∎

## 4.9 The $[n]$-density Property for Large $n$

Let $L$ be a non-trivial lattice.

**Theorem 4.9.1** *There exists $m \in \mathbf{N}$ such that for all $n \ge m$, $L$ has the $[n]$-density property.*

**Proof.** Let $h(P)$ denote the height of the poset $P$ and let $c_i(P)$ denote the number of (proper) chains of length $i$ of $P$. Let $h = h(L)$.

The proof of the theorem requires the following result from the combinatorics of posets:

**Lemma 4.9.2** *Let $P$ be a poset. Then*

$$\left| P^{[n+1]} \right| = \sum_{i=0}^{h(P)} c_i(P) \binom{n}{i}.$$

The sum on the right-hand side is known as the Zeta polynomial of $P$ (see Stanley [41]).

**Proof.** The chains of length $i$ of $P$ are the subsets of $P$ isomorphic to $[i+1]$. Since the number of order-preserving maps from $[n+1]$ onto $[i+1]$ is given by $\binom{n}{i}$ and since the image of every order-preserving map $f : [n+1] \to P$ is a chain, the result follows. ∎

By rewriting $\binom{n}{i}$ as $\frac{(n)_i}{i!}$, Lemma 4.9.2 shows that $\left| P^{[n+1]} \right|$ is a polynomial in $n$ of degree $h(P)$, where the coefficient of $n^{h(P)}$ is given by $\frac{c_{h(P)}}{h(P)!}$.

By Theorem 4.6.4 we can assume that $h < \left| J(L) \right|$ (otherwise $L \in \mathrm{MPf}(P)$ for every poset $P$, and we are done).

First suppose that there is a join-irreducible $a \in L$ such that $h([a)) \leq h - 2$. Since the degree of $\left|[a)^{[n+1]}\right|$ as a polynomial in $n$ is $h([a)) \leq h - 2$ while the degree of $\left|L^{[n+1]}\right|$ is $h$,

$$\lim_{n \to \infty} (n+2) \frac{\left|[a)^{[n+1]}\right|}{\left|L^{[n+1]}\right|} = 0.$$

Since the number of filters of $[n+1]$ is $n+2$, this implies that $(L, a)$ has the $[n+1]$-density property for all sufficiently large $n$.

Now suppose that for every join-irreducible $a \in L$, $h([a)) \geq h - 1$. Then $L$ is atomic, so by assumption on $L$, $L$ has at least $h+1$ atoms. Every chain $C$ of length $h$ of $L$ consists of $\hat{0}$ followed by a chain $C' = C \setminus \{\hat{0}\}$ of length $h - 1$ of $[a)$ for some atom $a$ of $L$. This shows that

$$c_h(L) = \sum_{a \in J(L)} c_{h-1}([a)),$$

so that for some atom $a \in L$, $c_{h-1}([a)) \leq \frac{1}{h+1} c_h(L)$. For this atom, $\left|[a)^{[n+1]}\right|$ is a polynomial of degree $h - 1$ in $n$, so that

$$\lim_{n \to \infty} (n+2) \frac{\left|[a)^{[n+1]}\right|}{\left|L^{[n+1]}\right|} = \frac{c_{h-1}([a))}{c_h(L)/h} \leq \frac{h}{h+1} < 1.$$

Thus $(L, a)$ has the $n$-density property for all sufficiently large $n$. ∎

## 4.10   A Bound on the Minimum $P$-density of Join-irreducibles

**Definition.** A *multi-(semi)lattice* $(L, \alpha)$ consists of a (semi)lattice $L$ and a *multiplicity function* $\alpha : L \to \mathbf{P}$. For $u \in L$, $\alpha(u)$ is the *multiplicity* of $u$ in $(L, \alpha)$. For $M \subseteq L$, the *($\alpha$-)cardinality* of $M$ is given by

$$|M|_\alpha = \sum_{u \in M} \alpha(u).$$

Note that every semilattice $L$ can be considered as a multi-semilattice $(L, \alpha)$ where $\alpha(u) = 1$ for $u \in L$.

Let $P$ be a poset and $L$ a lattice. Let $p$ be the number of filters of $P$.

**Theorem 4.10.1** *Let $(L^P, \alpha)$ be a multi-lattice with $\alpha(\hat{1}_{L^P}) = 1$. Let $n \geq 1$. If for every join-irreducible $a \in L$, $\left|L^P \setminus [a)^P\right|_\alpha \leq n$, then $\left|L^P\right|_\alpha \leq M(n)$, where $M(n)$ is defined by*

$$M(n) = \begin{cases} \max\{kn - p^{k-1}(k(p-1) - p) \\ \qquad \mid k \geq 1 \text{ and } p^{k-1}(p-1) \leq n\} & \text{if } n \geq p - 1, \\ \\ 1 & \text{otherwise.} \end{cases}$$

*This bound is best possible.*

Asymptotically, $M(n) \sim n \log_p(n)$ (Theorem 4.10.5).

**Proof.**  Suppose that for every join-irreducible $a \in L$, $\left| L^P \setminus [a]^P \right|_\alpha \leq n$. It will be shown that $\left| L^P \right|_\alpha \leq M(n)$.

Consider first the case $n \leq p - 2$. If $|L| \geq 2$, then for every proper join-irreducible $a \in L$, $\left| L^P \setminus [a]^P \right|_\alpha \geq \left| L^P \setminus [a]^P \right| \geq p - 1$ (count the number of maps $f : P \to \{\hat{0}, a\}$ with $\hat{0}$ in the range of $f$). Therefore, given that $n \leq p - 2$, $L$ must be the one-element lattice, and since $\alpha(\hat{1}_{L^P}) = 1$, $\left| L^P \right|_\alpha = 1 = M(n)$.

Now assume that $n \geq p - 1$. Then $M(n) \geq n + 1$. Hence we can assume that $|L| \geq 2$. We reduce the problem to the case where $L$ is a Boolean lattice.

Let $M$ be a minimal subset of $J(L)$ with $\bigvee M = \hat{1}$. Since $L$ is non-trivial, $M \neq \emptyset$. Let $B = 2^M$ and $k = |M|$. Then $B$ is a Boolean lattice generated by $k$ atoms. Let $\pi : L \to B$ be the restriction map defined by

$$\pi(u) = M \cap (u].$$

**Lemma 4.10.2** *The map $\pi$ is onto $B$.*

**Proof.**  Let $N \subseteq M$. We have $\pi(\bigvee N) \supseteq N$ and $\bigvee \pi(\bigvee N) = \bigvee N$. Consider $\bigvee(N \cup (M \setminus \pi(\bigvee N)))$:

$$
\begin{aligned}
\bigvee(N \cup (M \setminus \pi(\bigvee N))) &= (\bigvee N) \vee (\bigvee(M \setminus \pi(\bigvee N))) \\
&= (\bigvee \pi(\bigvee N)) \vee (\bigvee(M \setminus \pi(\bigvee N))) \\
&= \bigvee M \\
&= \hat{1}.
\end{aligned}
$$

Minimality of $M$ implies that $N \cup (M \setminus \pi(\bigvee N)) = M$, hence $\pi(\bigvee N) = N$, as desired. ∎

Since $\pi : L \to B$ is order-preserving, $\pi$ induces a map $\overline{\pi} : L^P \to B^P$ defined by

$$\overline{\pi}(f) = \pi \circ f,$$

for $f \in L^P$. The map $\overline{\pi}$ is onto $B^P$. To see this, let $\sigma : B \to L$ be the order-preserving map defined by $\sigma(N) = \bigvee N$ for $N \subseteq M$. If $g \in B^P$, then $\sigma \circ g \in L^P$. By the proof of Lemma 4.10.2, $\pi(\sigma(N)) = N$, so that $\overline{\pi}(\sigma \circ g) = g$ for every $g \in B^P$.

Let $(B^P, \beta)$ be the multi-lattice defined by

$$\beta(g) = \left| \{f \mid \overline{\pi}(f) = g\} \right|_\alpha.$$

Since $\overline{\pi}$ is onto, $\beta(g) > 0$ for every $g \in B^P$. Since $\bigvee M = \hat{1}_L$, we have $\pi(u) = \hat{1}_B$ iff $u = \hat{1}_L$, which gives $\overline{\pi}(f) = \hat{1}_{B^P}$ iff $f = \hat{1}_{L^P}$. Thus $\beta(\hat{1}_{B^P}) = \alpha(\hat{1}_{L^P}) = 1$. Let $a \in M$, so that $\{a\}$ is an atom of $B$. If $h \in L^P$ and $\overline{\pi}(h) \in B^P \setminus [\{a\}]^P$, then $h \in L^P \setminus [a]^P$. This implies that for every atom $\{a\}$ of $B$, $\left|B^P \setminus [\{a\}]^P\right|_\beta \leq n$. Since $\left|B^P\right|_\beta = \left|L^P\right|_\alpha$, the reduction is complete.

It is now shown that we can assume every $g \in B^P$ with $\beta(g) > 1$ is a coatom of $B^P$. The coatoms of $B^P$ are given by the functions $g_{\bar{a},x}$ where $\bar{a}$ is the complement of an atom $a \in B$, $x$ is a minimal member of $P$ and

$$g_{\bar{a},x}(y) = \begin{cases} \bar{a} & \text{if } y = x, \\ \hat{1} & \text{otherwise,} \end{cases}$$

for $y \in P$.

Let $C$ be the set of coatoms of $B^P$. Suppose that $g \in B^P \setminus C$ and $\beta(g) > 1$. Let $g_{\bar{a},x}$ be a coatom with $g_{\bar{a},x} \geq g$. Let $\beta$ be modified to the multiplicity function $\beta'$ defined by

$$\beta'(f) = \begin{cases} 1 & \text{if } f = g, \\ \beta(f) + \beta(g) - 1 & \text{if } f = g_{\bar{a},x}, \\ \beta(f) & \text{otherwise,} \end{cases}$$

for $f \in B^P$. Then $\left|B^P\right|_{\beta'} = \left|B^P\right|_\beta$, $\beta'(\hat{1}_{B^P}) = \beta(\hat{1}_{B^P}) = 1$, and since $g_{\bar{a},x} \geq g$, for every atom $b \in B$, $\left|B^P \setminus [b)^P\right|_{\beta'} \leq \left|B^P \setminus [b)^P\right|_\beta \leq n$.

If we successively perform this modification for every $g \in B^P \setminus C$, then the resulting multiplicity function has the desired property.

Let $m$ be the number of minimal members of $P$. Then $|C| = mk$. We have $\left|B^P\right| = p^k$ (Example 4.2.4). Let $r$ be given by

$$\begin{aligned} r &= \left|B^P\right|_\beta - \left|B^P\right| \\ &= \left|B^P\right|_\beta - p^k, \end{aligned}$$

so that $r$ is the total *excess* multiplicity of $B^P$. Since every map $f \in B^P$ with $\beta(f) > 1$ is a coatom,

$$r = |C|_\beta - |C|.$$

Let $a$ be an atom of $B$. Since $[a)$ is a Boolean algebra generated by $k - 1$ atoms, $\left|[a)^P\right| = p^{k-1}$. If $g \in C \setminus [a)^P$, then $g = g_{\bar{a},x}$ for some minimal $x \in P$. This implies that the family $\{C \setminus [b)^P \mid b \text{ is an atom of } B\}$ partitions $C$ into $k$ sets of cardinality $m$. Let

$$s(a) = \left|C \setminus [a)^P\right|_\beta - m,$$

so that $s(a)$ is the total excess multiplicity of $C \setminus [a)^P$.

Let $a$ be an atom of $B$ such that $s(a)$ is maximal. Then $r \leq ks(a)$. We have

$$\left|B^P \setminus [a)^P\right|_\beta = p^k - p^{k-1} + s(a) \leq n,$$

and

$$\left| B^P \right|_\beta = p^k + r \le p^k + ks(a).$$

The first inequality yields $ks(a) \le k(n - p^{k-1}(p-1))$. Using the second inequality, we obtain

$$(*) \qquad \begin{aligned} \left| B^P \right|_\beta &\le p^k + k(n - p^{k-1}(p-1)) \\ &= kn - p^{k-1}(k(p-1) - p). \end{aligned}$$

Since $kn - p^{k-1}(k(p-1) - p) \le M(n)$, this completes the proof of the bound.

To see that the bound is best possible, observe that if $s(a) = n - (p^{k-1}(p-1))$ for every atom $a \in B$, then equality holds in $(*)$ and for every atom $a \in B$, $\left| B^P \setminus [a]^P \right|_\beta = n$. Thus $\left| B_k \right|_\beta = M(n)$ for some $k$ and $\beta$. ∎

Let $B_k^\circ = B_k \setminus \{\hat{1}\}$, where $B_k$ is the Boolean lattice generated by $k$ atoms. Let $\delta_k = \left| B_k^P \setminus (B_k^\circ)^P \right|$, so that $\left| (B_k^\circ)^P \right| = p^k - \delta_k$.

**Theorem 4.10.3** *If $P$ has a greatest member or $k = 1$, then $\delta_k = (p-1)^k$. In general, $\delta_k \le (p-1)^k$.*

**Proof.** Suppose that $P$ has a greatest member $\hat{1}_P$. Then $f \in B_k^P \setminus (B_k^\circ)^P$ iff $f(\hat{1}_P) = \hat{1}_{B_k}$. Thus $\left| B_k^P \setminus (B_k^\circ)^P \right| = \left| B_k^{P^\circ} \right|$ where $P^\circ = P \setminus \{\hat{1}_P\}$. Since the number of filters of $P^\circ$ is $p - 1$, $\delta_k = (p-1)^k$.

Suppose that $P$ does not have a greatest member. Consider first the case $k = 1$. Since $B_1 \cong [2]$, the maps in $B_1^P \setminus (B_1^\circ)^P$ are determined by the non-empty filters of $P$. Thus $\delta_1 = (p-1)$. For $k \ge 2$, let $\pi_i$ denote the projection of $B_k$ onto the $i$'th component of the representation $B_k \cong B_1^k \cong [2]^k$. Consider the family $F$ of maps $f \in B_k^P$ such that for each $i \in 1, \dots, k$, there is an $x \in P$ with $\pi_i(f(x)) = \hat{1}$. Then $F \supseteq B_k^P \setminus (B_k^\circ)^P$. Since $F \cong (B_1^P \setminus (B_1^\circ)^P)^k$, $\delta_k \le |F| = (p-1)^k$. ∎

**Theorem 4.10.4** *Let $(L^P, \alpha)$ be a multi-semilattice where $L$ is a semilattice with no greatest member. Let $n \ge 1$. If for every join-irreducible $a \in L$, $\left| L^P \setminus [a]^P \right|_\alpha \le n$, then $\left| L^P \right|_\alpha \le M^\circ(n)$ where $M^\circ(n)$ is defined by*

$$\begin{aligned} M^\circ(n) = \max\{ &kn - p^{k-1}(k(p-1) - p) + (k-1)\delta_k - k\delta_{k-1} \\ &\mid k \ge 2 \text{ and } p^{k-1}(p-1) - \delta_k + \delta_{k-1} \le n \}. \end{aligned}$$

*This bound is best possible.*

The requirement that $k \ge 2$ in the expression being maximized implies that for $n < p(p-1) - \delta_2 + \delta_1$, $M^\circ(n) = 0$ (using the convention that $\max \emptyset = 0$). This reflects the fact that if $L$ is a semilattice with no greatest member, then $L$ has two incomparable join-irreducibles $a$ and $b$, so that

$$\left| L^P \setminus [a]^P \right| \ge \left| \{a, b, \hat{0}\}^P \setminus [a]^P \right| = (p^2 - \delta_2) - (p - \delta_1).$$

Asymptotically, $M^o(n) \sim M(n) \sim n \log_p(n)$ (Theorem 4.10.5).

**Proof.** Let $\widehat{L} = L \cup \{\widehat{1}_L\}$ be the completion of $L$, where $\widehat{1}_L$ is a new greatest element for $L$. We have $\widehat{L} \geq 4$. We can perform the reduction to the Boolean lattice $B$ as in the proof of Theorem 4.10.1, using $\widehat{L}$ instead of $L$. We obtain the multi-semilattice $((B^o)^P, \beta)$ where $B^o = B \setminus \widehat{1}_B$. Let $k$ be the number of atoms of $B$. Since $\widehat{1}_L$ is not join-irreducible, $k \geq 2$.

Let $C$ be the set of maximal members of $(B^o)^P$. Then $C$ consists of the maps $g_{\bar{a}} : P \to B^o$ where $\bar{a}$ is the complement of an atom $a \in B$ and $g_{\bar{a}}(x) = \bar{a}$ for every $x \in P$. We can assume that if $\beta(g) > 1$ then $g \in C$ (this follows as in the proof of Theorem 4.10.1, using the maps $g_{\bar{a}}$ instead of the maps $g_{\bar{a},x}$).

Since $B^o = B^o_k$, $\left|(B^o)^P\right| = p^k - \delta_k$. Thus the total excess multiplicity $r$ is given by:

$$r = \left|(B^o)^P\right|_\beta - \left|(B^o)^P\right| = \left|(B^o)^P\right|_\beta - (p^k - \delta_k).$$

Since every map $f \in (B^o)^P$ with $\beta(f) > 1$ is maximal,

$$r = |C|_\beta - |C|.$$

Let $s(a) = \beta(g_{\bar{a}}) - 1$. Choose an atom $a \in B$ such that $s(a)$ is maximal. Since $C \setminus [a]^P = \{g_{\bar{a}}\}$, we have

$$\left|(B^o)^P \setminus [a]^P\right|_\beta = (p^k - \delta_k) - (p^{k-1} - \delta_{k-1}) + s(a) \leq n.$$

Since $|C| = k$, we have $r \leq ks(a)$, which yields

$$\left|(B^o)^P\right|_\beta = p^k - \delta_k + r \leq p^k - \delta_k + ks(a).$$

The first inequality gives $ks(a) \leq k(n - p^{k-1}(p-1) + (\delta_k - \delta_{k-1}))$. Using the second inequality, we obtain

$$(*) \qquad \begin{aligned} \left|(B^o)^P\right|_\beta &\leq p^k - \delta_k + k(n - p^{k-1}(p-1) + (\delta_k - \delta_{k-1})) \\ &= kn - p^{k-1}(k(p-1) - p) + (k-1)\delta_k - k\delta_{k-1}. \end{aligned}$$

Thus $\left|(B^o)^P\right|_\beta \leq M^o(n)$.

If $s(a) = n - (p^{k-1}(p-1)) + (\delta_k - \delta_{k-1})$ for every atom $a \in B$, then equality holds in $(*)$. This implies that the bound $M(n)$ is attained by $(B_k, \beta)$ for some $k$ and $\beta$. ∎

**Theorem 4.10.5**

$$\begin{aligned} M(n) &= n \log_p(n)(1 + o(1)), \\ M^o(n) &= n \log_p(n)(1 + o(1)). \end{aligned}$$

**Proof.**   Define the functions $f$, $f^\circ$, $c$ and $c^\circ$ by

$$
\begin{aligned}
f(k) &= kn - p^{k-1}(k(p-1) - p), \\
f^\circ(k) &= kn - p^{k-1}(k(p-1) - p) + (k-1)\delta_k - k\delta_{k-1}, \\
c(k) &= p^{k-1}(p-1), \\
c^\circ(k) &= p^{k-1}(p-1) - \delta_k + \delta_{k-1}.
\end{aligned}
$$

Then

$$
\begin{aligned}
M(n) &= \max\{f(k) \mid k \geq 1 \text{ and } c(k) \leq n\}, \\
M^\circ(n) &= \max\{f^\circ(k) \mid k \geq 2 \text{ and } c^\circ(k) \leq n\}.
\end{aligned}
$$

Since $\delta_k \leq (p-1)^k$ (Theorem 4.10.3), either one of $c(k) \leq n$ and $c^\circ(k) \leq n$ implies that

$$(*) \qquad\qquad p^{k-1}(p-1) - (p-1)^{k-1}p \leq n.$$

We show that $(*)$ implies $k < \log_p(n) + 2$ for $n$ sufficiently large. Suppose that $k \geq \log_p(n) + 2$. Then

$$
\begin{aligned}
p^{k-1}(p-1) - (p-1)^{k-1}p &= p(p-1)(p^{k-2} - (p-1)^{k-2}) \\
&\geq p(p-1)(n - n^{\log_p(p-1)}),
\end{aligned}
$$

where we used the fact that $p^{k-2} - (p-1)^{k-2}$ is increasing in $k$. Since $p \geq 2$ and $n^{\log_p(p-1)} = o(n)$, it follows that $p(p-1)(n - n^{\log_p(p-1)}) > n$ for $n$ sufficiently large.

Let $2 \leq k < \log_p(n) + 2$. Then $p^{k-1}(k(p-1) - p) \geq 0$. Thus we can estimate $f(k)$ by

$$f(k) \leq n\log_p(n) + 2n = n\log_p(n)(1 + o(1)).$$

For $k = 1$, $f(k) = n + 1$, hence

$$M(n) \leq n\log_p(n)(1 + o(1)).$$

Using $\delta_k \leq (p-1)^k$, we get

$$
\begin{aligned}
f^\circ(k) &\leq n\log_p(n) + 2n + (\log_p(n) + 1)(p-1)^2 n^{\log_p(p-1)} \\
&= n\log_p(n)(1 + o(1)),
\end{aligned}
$$

which yields

$$M^\circ(n) \leq n\log_p(n)(1 + o(1)).$$

To show that the bounds on $M(n)$ and $M^\circ(n)$ are asymptotically optimal, let $k_n = \log_p(n) - \log_p(t_n)$, where $t_n$ is any sequence such that $\lim_{n\to\infty} t_n = \infty$ and $t_n = o(n)$. Then

$$k_n n = n\log_p(n)(1 + o(1))$$

and

$$k_n p^{k_n} = \frac{n \log_p(n)}{t_n}(1 + o(1)) = o(n \log_p(n)).$$

This implies

$$f(k_n) \geq n \log_p(n) - o(n \log_p(n))$$

and

$$f^o(k_n) \geq n \log_p(n) - o(n \log_p(n)).$$

Since $p^{k_n} = o(n)$, the constraints $c(k_n) \leq n$ and $c^o(k_n) \leq n$ are satisfied for sufficiently large $n$. The result follows. ∎

Let $\mathcal{L}$ be the class of all lattices and $\mathcal{S}$ the class of all semilattices. Since $\mathcal{L}$ and $\mathcal{S}$ can be considered as subclasses of the class of multi-lattices and the class of multi-semilattices respectively, we have the following corollary of Theorems 4.10.1 and 4.10.4:

**Corollary 4.10.6**

$$\begin{aligned} \overline{h}(\mathcal{L}, P, n) &\leq M(n), \\ \overline{h}(\mathcal{S}, P, n) &\leq \max(M(n), M^o(n)). \end{aligned}$$

Neither bound is exact. For example, suppose that $3 \leq p \leq n < 2p - 2$. By Observation 4.8.2, $\overline{h}(\mathcal{S}, P, n) = \overline{h}(\mathcal{L}, P, n) = p$. Since

$$2p - 2 = p(p - 1) - (p - 1)^2 + (p - 1) \leq p(p - 1) - \delta_2 + \delta_1,$$

it follows that $M^o(n) = 0$. If $k \geq 1$ and $p^{k-1}(p - 1) \leq n$, then $k = 1$, so that $M(n) = n + 1 > p$.

## 4.11 Union-closed Families of Sets Generated by Graphs

**Definition.** The family of sets $\mathcal{G}$ is a *graph* iff for every $U \in \mathcal{G}$, $|U| \leq 2$. The members of $\mathcal{G}$ are referred to as *edges*. The graph $\mathcal{G}$ is *simple* iff for every $U \in \mathcal{G}$, $|U| = 2$. For $x \in \bigcup \mathcal{G}$, the *degree* in $\mathcal{G}$ of $x$ is

$$\mathrm{d}_{\mathcal{G}}(x) = \left| \mathcal{G}_{\supseteq \{x\}} \right|.$$

If $U$ is an edge of $G$, the *degree* in $\mathcal{G}$ of $U$ is

$$\mathrm{d}_{\mathcal{G}}(U) = \left| \{ V \in \mathcal{G} \setminus \{U\} \mid V \cap U \neq \emptyset \} \right|.$$

**Definition.** Let $\mathcal{F}$ be a union-closed family of sets. Let $G(\mathcal{F})$ denote the family of generators of $\mathcal{F}$ and define $J(\mathcal{F}) = G(\mathcal{F}) \setminus \{\emptyset\}$. The *closure* in $\mathcal{F}$ of the set $X$ is given by

$$\pi_{\mathcal{F}}(X) = \bigcup \mathcal{F}_{\subseteq X} = \bigcup \{ V \in J(\mathcal{F}) \mid V \subseteq X \}.$$

(This agrees with the notation for projections onto lub-subsemilattices introduced in Section 4.4.) The set of *isolated* elements of $X$ is

$$\overline{\pi}_{\mathcal{F}}(X) = X \setminus \pi_{\mathcal{F}}(X).$$

In terms of hypergraphs, $\overline{\pi}_{\mathcal{F}}(X)$ is the set of isolated points of the hypergraph $(\mathcal{F}_{\subseteq X}, X)$.

**Observation 4.11.1** *Let $\mathcal{F}$ be a union-closed family of sets. The following are equivalent:*

> *(i)* $U \in \mathcal{F}$,
>
> *(ii)* $\pi_{\mathcal{F}}(U) = U$,
>
> *(iii)* $\overline{\pi}_{\mathcal{F}}(U) = \emptyset$.

In this section we will prove:

**Theorem 4.11.2** *If $\mathcal{F}$ is a union-closed family of sets such that $\emptyset \in \mathcal{F}$ and $J(\mathcal{F})$ is a non-empty graph, then $\mathcal{F}$ has the density property.*

Let $\mathcal{F}$ be a union-closed family of sets such that $\emptyset \in \mathcal{F}$ and let $U \subseteq \bigcup \mathcal{F}$. Recall that $\mathcal{N}_{\mathcal{F}}(U)$ (the lattice neighborhood in $\mathcal{F}$ of $U$) is the union-closed family generated by the empty set and the members $V$ of $J(\mathcal{F})$ with $V \cap U \neq \emptyset$. The third lattice neighborhood is given by $\mathcal{N}_{\mathcal{F}}^3(U) = \mathcal{N}_{\mathcal{F}}(\bigcup \mathcal{N}_{\mathcal{F}}(U))$.

To prove Theorem 4.11.2, we develop a technique for estimating the density of $\mathcal{F}_{\supseteq U}$ in $\mathcal{F}$. This estimate depends only on the third lattice neighborhood $\mathcal{N}_{\mathcal{F}}^3(U)$. The estimate will be used to show that if $J(\mathcal{F})$ is a graph and $U \in \mathcal{G} = \{V \in J(\mathcal{F}) \,|\, |V| = 2\}$ has minimal degree in $\mathcal{G}$, then $(\mathcal{F}, U)$ has the density property (Theorem 4.11.28).

**Definition.** Let $\mathcal{H}$ be a family of sets. The *density* in $\mathcal{H}$ of the set $X$ is

$$\rho_{\mathcal{H}}(X) = \frac{\left|\mathcal{H}_{\supseteq X}\right|}{|\mathcal{H}|}.$$

Thus $\rho_{\mathcal{H}}(X)$ is the density of $\mathcal{H}_{\supseteq X}$ in $\mathcal{H}$. The reciprocal of $\rho_{\mathcal{H}}(X)$ is denoted by $(\frac{1}{\rho})_{\mathcal{H}}(X)$. Let $\rho = \rho_{\mathcal{F}}(U)$.

**Definition.** The union-closed family $\mathcal{F}'$ is a *(conservative) extension* of $(\mathcal{F}, U)$ iff there is a non-empty union-closed family of sets $\mathcal{H}$ such that $(\bigcup \mathcal{H}) \cap U = \emptyset$ and $\mathcal{F}' = \mathcal{F} \vee \mathcal{H}$.

Associativity of $\vee$ for families of sets yields:

**Observation 4.11.3** *The extension relation is transitive; i.e. if $\mathcal{F}_1$ is an extension of $(\mathcal{F}, U)$ and $\mathcal{F}_2$ is an extension of $(\mathcal{F}_1, U)$, then $\mathcal{F}_2$ is an extension of $(\mathcal{F}, U)$.*

Since $\mathcal{F} = \mathcal{F} \vee \{\emptyset\}$:

**Observation 4.11.4** *$\mathcal{F}$ is an extension of $(\mathcal{F}, U)$.*

**Definition.** Let

$$\frac{1}{\bar{\rho}} = \inf\{(\tfrac{1}{\rho})_{\mathcal{F}'}(U) \mid \mathcal{F}' \text{ is an extension of } (\mathcal{F}, U)\}.$$

**Observation 4.11.5** $\frac{1}{\bar{\rho}} \leq \frac{1}{\rho}$.

**Observation 4.11.6** *If $U \in J(\mathcal{F})$ and $\frac{1}{\bar{\rho}} \geq 2$, then $(\mathcal{F}, U)$ has the density property.*

The goal is to find lower bounds on $\frac{1}{\bar{\rho}}$ in terms of the local properties of $\mathcal{F}$ at $U$. To this end, let $\mathcal{F}'$ be an (arbitrary) extension of $(\mathcal{F}, U)$ and consider

$$(\tfrac{1}{\rho})_{\mathcal{F}'}(U) = \frac{|\mathcal{F}'|}{|\mathcal{F}'_{\supseteq U}|}.$$

**Definition.** For $X \cap U = \emptyset$, let

$$T_{\mathcal{F}'}(X) = \{Y \subseteq U \mid X \cup Y \in \mathcal{F}'\}.$$

We have

$$|\mathcal{F}'| = \sum_{X \in \mathcal{F}'_{\setminus U}} |T_{\mathcal{F}'}(X)|$$

and

$$|\mathcal{F}'_{\supseteq U}| \leq |\mathcal{F}'_{\setminus U}|,$$

where equality holds if $U \in \mathcal{F}$ (since $\mathcal{F}'$ is closed under union with members of $\mathcal{F}$). This gives

$$(*) \qquad (\tfrac{1}{\rho})_{\mathcal{F}'}(U) \;\; \geq \;\; \frac{\sum_{X \in \mathcal{F}'_{\setminus U}} |T_{\mathcal{F}'}(X)|}{|\mathcal{F}'_{\setminus U}|}.$$

Note that if $U \in \mathcal{F}$, then this is an identity, so that $\frac{1}{\bar{\rho}} \leq 2^{|U|}$. We will determine a lower bound for $|T_{\mathcal{F}'}(X)|$ which is independent of $\mathcal{F}'$.

**Lemma 4.11.7** *Let $X \in \mathcal{F}'_{\setminus U}$. If $Y \subseteq U$ satisfies*

*(i) $\pi_{\mathcal{F}}(X \cup Y) \cap U = Y$,*

*(ii) $\pi_{\mathcal{F}}(X \cup Y) \supseteq \pi_{\mathcal{F}}(X \cup U) \setminus U$,*

*then $X \cup Y \in \mathcal{F}'$.*

Conditions (i) and (ii) are independent of $\mathcal{F}'$. Observe that (i) is equivalent to

$$\overline{\pi}_{\mathcal{F}}(X \cup Y) \cap U = \emptyset$$

and (ii) is equivalent to

$$\overline{\pi}_{\mathcal{F}}(X \cup Y) \setminus U \subseteq \overline{\pi}_{\mathcal{F}}(X \cup U).$$

Lemma 4.11.10 below gives conditions equivalent to (i) and (ii) in terms of the family of generators of $\mathcal{F}$.

**Proof.**  Since $X \in \mathcal{F}'_{\setminus U}$, there exists $Z \subseteq U$ such that $X \cup Z \in \mathcal{F}'$. Since $\mathcal{F}'$ is an extension of $(\mathcal{F}, U)$, there is a union-closed (non-empty) family $\mathcal{H}$ such that $(\bigcup \mathcal{H}) \cap U = \emptyset$ and $\mathcal{F}' = \mathcal{F} \vee \mathcal{H}$. We have $X \cup Z = A \cup B$ for some $A \in \mathcal{F}$ and $B \in \mathcal{H}$. We can assume that $A = \pi_{\mathcal{F}}(X \cup Z)$. Let $A' = \pi_{\mathcal{F}}(X \cup Y)$. By (i), $A' \cap U = Y$. If we can show that $X \setminus A' \subseteq B$, then $X \cup Y = A' \cup B \in \mathcal{F}'$ and we are done.

Let $x \in X \setminus A'$. Then $x \in \overline{\pi}_{\mathcal{F}}(X \cup Y) \setminus U$. By (ii), $x \in \overline{\pi}_{\mathcal{F}}(X \cup U) \setminus U$. Since $Z \subseteq U$, $\overline{\pi}_{\mathcal{F}}(X \cup U) \subseteq \overline{\pi}_{\mathcal{F}}(X \cup Z)$, which yields $x \in \overline{\pi}_{\mathcal{F}}(X \cup Z) \setminus U$. Since $\overline{\pi}_{\mathcal{F}}(X \cup Z) \setminus U = X \setminus A \subseteq B$, we have $x \in B$, as desired.  ∎

**Definition.**  For every set $X$ with $X \cap U = \emptyset$, let $E_{\mathcal{F},U}(X)$ consist of the subsets $Y$ of $U$ satisfying the conditions of Lemma 4.11.7:

$$E_{\mathcal{F},U}(X) = \{Y \subseteq U \mid \pi_{\mathcal{F}}(X \cup Y) \cap U = Y \text{ and} \\ \pi_{\mathcal{F}}(X \cup Y) \supseteq \pi_{\mathcal{F}}(X \cup U) \setminus U\}.$$

Let $E(X) = E_{\mathcal{F},U}(X)$.

**Example 4.11.8**  Suppose that $\mathcal{F}$ is the union-closed family of sets generated by the edges of the graph depicted in Figure 4.4 and the empty set. Then

$$J(\mathcal{F}) = \Big\{\{a, b\}, \{x_1, a\}, \{x_2, a\}, \{x_3, a\}, \{x_3, b\}, \{x_4, b\}, \{x_4, x_5\}\Big\}.$$

Let $U = \{a, b\}$. We have

$$E(\{x_4, x_5\}) = \Big\{\emptyset, \{b\}, \{a, b\}\Big\},$$
$$E(\{x_3\}) = \Big\{\{a\}, \{b\}, \{a, b\}\Big\},$$
$$E(\{x_1, x_4\}) = \Big\{\{a, b\}\Big\}.$$

By Lemma 4.11.7, $E(X) \subseteq T_{\mathcal{F}'}(X)$ for every $X \in \mathcal{F}'_{\setminus U}$. Using inequality $(*)$, we get

$$(**) \qquad (\tfrac{1}{\rho})_{\mathcal{F}'}(U) \geq \frac{\sum_{X \in \mathcal{F}'_{\setminus U}} |E(X)|}{|\mathcal{F}'_{\setminus U}|}.$$

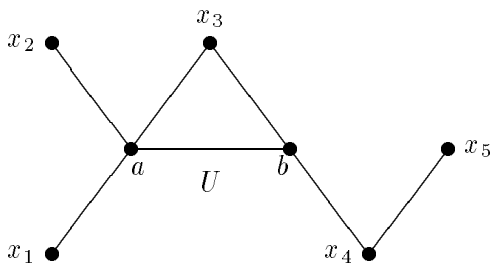

Figure 4.4

**Definition.** Let

$$\mu_{\mathcal{F},\mathcal{F}'}(U) = \frac{\sum_{X \in \mathcal{F}'_{\setminus U}} |E(X)|}{|\mathcal{F}'_{\setminus U}|}.$$

Let $\mu_{\mathcal{F}'} = \mu_{\mathcal{F},\mathcal{F}'}(U)$.

By inequality (**):

**Observation 4.11.9** $\mu_{\mathcal{F}'} \leq (\frac{1}{\rho})_{\mathcal{F}'}(U)$.

We can define $E(X)$ in terms of the set of generators of $\mathcal{F}$.

**Lemma 4.11.10** *The set $Y$ is in $E(X)$ iff*

(i)' *for every $x \in Y$ there exists $V \in J(\mathcal{F})$ with $x \in V \subseteq X \cup Y$,*

(ii)'' *for every $V \in J(\mathcal{F})$, if $V \setminus U \subseteq X$, then for every $x \in V \setminus U$ there exists $V' \in J(\mathcal{F})$ with $x \in V' \subseteq (X \cup Y)$.*

**Proof.** In fact, condition (i) of Lemma 4.11.7 is equivalent to (i)' and condition (ii) of Lemma 4.11.7 is equivalent to (ii)''.

If $x \in W \in \mathcal{F}$, then there is a generator $V \in J(\mathcal{F})$ such that $x \in V \subseteq W$. This implies that $\pi_{\mathcal{F}}(X \cup Y) \cap U = Y$ iff (i)' holds, hence (i) iff (i)'.

Suppose that $X$ and $Y$ satisfy (ii). Let $V \in J(\mathcal{F})$ and $V \setminus U \subseteq X$. Then $V \subseteq \pi_{\mathcal{F}}(X \cup U)$. By (ii), $\pi_{\mathcal{F}}(X \cup Y) \supseteq \pi_{\mathcal{F}}(X \cup U) \setminus U$. Hence, if $x \in V \setminus U$, then $x \in \pi_{\mathcal{F}}(X \cup Y)$, which implies that there is a generator $V'$ of $\mathcal{F}$ such that $x \in V' \subseteq X \cup Y$. Thus (ii)'' holds.

Conversely, suppose that $X$ and $Y$ satisfy (ii)''. We show that $\pi_{\mathcal{F}}(X \cup Y) \supseteq \pi_{\mathcal{F}}(X \cup U) \setminus U$. Let $x \in \pi_{\mathcal{F}}(X \cup U) \setminus U$. Then there exists $V \in J(\mathcal{F})$ such that $x \in V \subseteq X \cup U$. We have $V \setminus U \subseteq X$, so by (ii)'', there exists $V' \in J(\mathcal{F})$ with $x \in V' \subseteq X \cup Y$. This implies that $x \in \pi_{\mathcal{F}}(X \cup Y)$, as desired. ∎

**Definition.** The *neighborhood* in $\mathcal{F}$ of a set $X$ is given by

$$N_{\mathcal{F}}(X) = X \cup \bigcup \{V \in J(\mathcal{F}) \mid V \cap X \neq \emptyset\}.$$

Let $N = N_{\mathcal{F}}(U)$ and $N^2 = N_{\mathcal{F}}(N_{\mathcal{F}}(U))$. Observe that if $U \in \mathcal{F}$, then $N = \bigcup \mathcal{N}_{\mathcal{F}}(U)$ and $N^2 = \bigcup \mathcal{N}_{\mathcal{F}}^3(U)$.

**Lemma 4.11.11** *The family $\mathcal{F}$ is an extension of $(\mathcal{N}_{\mathcal{F}}^3(U), U)$.*

**Proof.** Let $\mathcal{H}$ be the union-closed family of sets generated by the empty set and the generators $V$ of $\mathcal{F}$ with $V \notin \mathcal{N}_{\mathcal{F}}^3(U)$. Then $\mathcal{F} = \mathcal{N}_{\mathcal{F}}^3(U) \vee \mathcal{H}$ and $(\bigcup \mathcal{H}) \cap U = \emptyset$. (This expresses $\mathcal{F}$ as the internal subdirect product of $\mathcal{N}_{\mathcal{F}}^3(U)$ and $\mathcal{H}$, see Section 4.4.) ∎

**Theorem 4.11.12** *Let $\mathcal{F}'$ be an extension of $(\mathcal{F}, U)$. Then*

$$\mu_{\mathcal{F}, \mathcal{F}'}(U) = \mu_{\mathcal{N}_{\mathcal{F}}^3(U), \mathcal{F}'}(U).$$

**Proof.** By Lemma 4.11.11 and Observation 4.11.3, $\mathcal{F}'$ is an extension of $(\mathcal{N}_{\mathcal{F}}^3(U), U)$, so that $\mu_{\mathcal{N}_{\mathcal{F}}^3(U), \mathcal{F}'}(U)$ is well-defined. By Lemma 4.11.10 and by definition of $\mathcal{N}_{\mathcal{F}}^3(U)$, whether a given subset of $U$ is in $E(X)$ depends only on the generators $V \in J(\mathcal{F})$ with $V \in \mathcal{N}_{\mathcal{F}}^3(U)$. Since $J(\mathcal{N}_{\mathcal{F}}^3(U)) = J(\mathcal{F}) \cap \mathcal{N}_{\mathcal{F}}^3(U)$, Lemma 4.11.10 implies that

$$E_{\mathcal{F}, U}(X) = E_{\mathcal{N}_{\mathcal{F}}^3, U}(X)$$

for every $X$ disjoint from $U$. The result now follows by definition of $\mu(U)$. ∎

**Definition.** Let $D \supseteq \bigcup \mathcal{F}$. An extension $\mathcal{F}'$ of $(\mathcal{F}, U)$ *minimizes* $\mu$ *in* $D$ iff $\bigcup \mathcal{F}' \subseteq D$ and for every extension $\mathcal{F}''$ of $(\mathcal{F}, U)$ with $\bigcup \mathcal{F}'' \subseteq D$, $\mu_{\mathcal{F}'} \leq \mu_{\mathcal{F}''}$.

**Theorem 4.11.13** *Let $D \supseteq \bigcup \mathcal{F}$. There exists an extension $\mathcal{F}' = \mathcal{F} \vee \mathcal{H}$ of $(\mathcal{F}, U)$ such that*

*(i) $\mathcal{H}$ is a filter of $2^{D \setminus U}$,*

*(ii) $\mathcal{F}'_{\setminus U} = \mathcal{H}$,*

*(iii) $\mathcal{F}'$ minimizes $\mu$ in $D$.*

**Proof.** Note that (i) implies (ii): Let $\mathcal{H}$ be a filter of $2^{D \setminus U}$. Since $\emptyset \in \mathcal{F}$, $(\mathcal{F} \vee \mathcal{H})_{\setminus U} \supseteq \mathcal{H}$. For the reverse inclusion, let $V \in \mathcal{F} \vee \mathcal{H}$. Then there is a $W \in \mathcal{H}$ with $W \subseteq V \setminus U$, hence $V \setminus U \in \mathcal{H}$.

Let $\mathcal{F} \vee \mathcal{G}$ be an extension of $(\mathcal{F}, U)$ such that $\bigcup \mathcal{G} \subseteq D \setminus U$ and $\mathcal{F} \vee \mathcal{G}$ minimizes $\mu$ in $D$. Since $(\mathcal{F} \vee (\mathcal{F} \vee \mathcal{G})_{\setminus U})_{\setminus U} = (\mathcal{F} \vee \mathcal{G})_{\setminus U}$, we have $\mu_{\mathcal{F} \vee (\mathcal{F} \vee \mathcal{G})_{\setminus U}} = \mu_{\mathcal{F} \vee \mathcal{G}}$. This implies that we can assume $\mathcal{G} = (\mathcal{F} \vee \mathcal{G})_{\setminus U}$.

Let $\mathcal{H}$ be the filter of $2^{D \setminus U}$ generated by $\mathcal{G}$ and let $\mathcal{F}' = \mathcal{F} \vee \mathcal{H}$. We show that $\mu_{\mathcal{F}'} \leq \mu_{\mathcal{F} \vee \mathcal{G}}$ which implies that $\mathcal{F}'$ is as desired.

For $X \in \mathcal{G}$, let

$$P(X) = \{Y \in \mathcal{H} \mid \pi_{\mathcal{G}}(Y) = X\}.$$

**Lemma 4.11.14** *Let $X, Y \in \mathcal{G}$ with $X \subseteq Y$. Then $\big|P(X)\big| \geq \big|P(Y)\big|$.*

**Proof.**    Define

$$\sigma(Z) = (Z \setminus Y) \cup X$$

for $Z \in P(Y)$. Since every $Z \in P(Y)$ includes $Y$, $\sigma$ is a one-to-one map. To show that $\sigma$ maps $P(Y)$ into $P(X)$, let $Z \in P(Y)$. Let $X' = \pi_{\mathcal{G}}(\sigma(Z))$ and suppose that $X' \neq X$. Then $X' \supset X$. By definition of $\sigma(Z)$, $X' \not\subseteq Y$. Let $Y' = Y \cup X'$. Then $Y' \in \mathcal{G}$ and $Z \supseteq Y' \supset Y$, contradicting $Z \in P(Y)$.    ■

**Lemma 4.11.15** *Let $X \in \mathcal{G}$. If $Y \in P(X)$, then $E(Y) = E(X)$.*

**Proof.**    We show that for every $Z \subseteq U$, $\pi_{\mathcal{F}}(X \cup Z) = \pi_{\mathcal{F}}(Y \cup Z)$. Let $Z \subseteq U$. Let $X' = \pi_{\mathcal{F}}(X \cup Z)$ and $Y' = \pi_{\mathcal{F}}(Y \cup Z)$. The inclusion $X \subseteq Y$ implies $X' \subseteq Y'$. Since $Y' \in \mathcal{F}$ and $X \in \mathcal{G}$, we have $Y' \cup X \in \mathcal{F} \vee \mathcal{G}$, so that $(Y' \cup X) \setminus U = (Y' \setminus U) \cup X \in \mathcal{G}$. Since $Y \in P(X)$ and $Y \supseteq (Y' \setminus U) \cup X$, it follows that $(Y' \setminus U) \cup X \subseteq \pi_{\mathcal{G}}(Y) = X$. Using $Y' \cap U \subseteq Z$, we get $Y' \subseteq X \cup Z$, hence $Y' \subseteq X'$.

The identities $\pi_{\mathcal{F}}(X \cup Z) = \pi_{\mathcal{F}}(Y \cup Z)$ and $\pi_{\mathcal{F}}(X \cup U) = \pi_{\mathcal{F}}(Y \cup U)$ imply that $Z \in E(X)$ iff $Z \in E(Y)$, as required.    ■

For $n \geq 0$, let

$$\mathcal{G}_{\leq n} = \{X \in \mathcal{G} \mid \big|P(X)\big| \leq n\}.$$

By Lemma 4.11.14, the $\mathcal{G}_{\leq n}$ are filters of $\mathcal{G}$. Since $(\mathcal{F} \vee \mathcal{G})_{\setminus U} = \mathcal{G}$, this implies that if $\mathcal{G}_{\leq n} \neq \emptyset$, then $\mathcal{F} \vee \mathcal{G}_{\leq n}$ is an extension of $(\mathcal{F}, U)$ such that $(\mathcal{F} \vee \mathcal{G}_{\leq n})_{\setminus U} = \mathcal{G}_{\leq n}$. Let $N$ be the maximum value of $\big|P(X)\big|$. Then $\mathcal{G}_{\leq N} = \mathcal{G}$. We use Lemmas 4.11.14 and 4.11.15 and the fact that the family $\{P(X) \mid X \in \mathcal{G}\}$ is a partition of $\mathcal{H}$ to compute $\mu_{\mathcal{F}'}$:

$$
\begin{aligned}
\mu_{\mathcal{F}'} &= \frac{\sum_{X \in \mathcal{F}'_{\setminus U}} \big|E(X)\big|}{\big|\mathcal{F}'_{\setminus U}\big|} \\
&= \frac{\sum_{X \in \mathcal{H}} \big|E(X)\big|}{\sum_{X \in \mathcal{H}} 1} \\
&= \frac{\sum_{X \in \mathcal{G}} \big|P(X)\big| \cdot \big|E(X)\big|}{\sum_{X \in \mathcal{G}} \big|P(X)\big|}
\end{aligned}
$$

$$= \frac{\sum_{n=0}^{N-1} \sum_{X \in \mathcal{G} \setminus \mathcal{G}_{\leq n}} \left| E(X) \right|}{\sum_{n=0}^{N-1} \sum_{X \in \mathcal{G} \setminus \mathcal{G}_{\leq n}} 1}$$

$$= \frac{\sum_{n=0}^{N-1} \left( \sum_{X \in \mathcal{G}} \left| E(X) \right| - \sum_{X \in \mathcal{G}_{\leq n}} \left| E(X) \right| \right)}{\sum_{n=0}^{N-1} |\mathcal{G}| - |\mathcal{G}_{\leq n}|}.$$

Define $\mu_{\mathcal{F} \vee \emptyset} = \mu_{\mathcal{F} \vee \mathcal{G}}$. We have $\mu_{\mathcal{F} \vee \mathcal{G}} = \left( \sum_{X \in \mathcal{G}} \left| E(X) \right| \right) / |\mathcal{G}|$ and if $\mathcal{G}_{\leq n} \neq \emptyset$, then $\mu_{\mathcal{F} \vee \mathcal{G}_{\leq n}} = \left( \sum_{X \in \mathcal{G}_{\leq n}} \left| E(X) \right| \right) / |\mathcal{G}_{\leq n}|$ and $\mu_{\mathcal{F} \vee \mathcal{G}} \leq \mu_{\mathcal{F} \vee \mathcal{G}_{\leq n}}$. This yields

$$\mu_{\mathcal{F}'} = \frac{\sum_{n=0}^{N-1} \left( \mu_{\mathcal{F} \vee \mathcal{G}} |\mathcal{G}| - \mu_{\mathcal{F} \vee \mathcal{G}_{\leq n}} |\mathcal{G}_{\leq n}| \right)}{\sum_{n=0}^{N-1} \left( |\mathcal{G}| - |\mathcal{G}_{\leq n}| \right)}$$

$$\leq \frac{\sum_{n=0}^{N-1} \left( \mu_{\mathcal{F} \vee \mathcal{G}} |\mathcal{G}| - \mu_{\mathcal{F} \vee \mathcal{G}} |\mathcal{G}_{\leq n}| \right)}{\sum_{n=0}^{N-1} \left( |\mathcal{G}| - |\mathcal{G}_{\leq n}| \right)}$$

$$= \mu_{\mathcal{F} \vee \mathcal{G}} \frac{\sum_{n=0}^{N-1} \left( |\mathcal{G}| - |\mathcal{G}_{\leq n}| \right)}{\sum_{n=0}^{N-1} \left( |\mathcal{G}| - |\mathcal{G}_{\leq n}| \right)}$$

$$= \mu_{\mathcal{F} \vee \mathcal{G}},$$

which implies that $\mathcal{F}'$ satisfies (iii), as required. ∎

**Lemma 4.11.16** *Let $X, Y$ be sets disjoint from $U$. If $Y \cap N^2 = X \cap N^2$, then $E(X) = E(Y)$.*

**Proof.** By Lemma 4.11.10, whether a given subset of $U$ is in $E(Z)$ depends only on the generators of $\mathcal{F}$ included in $N^2$. This implies that $E(Z)$ depends only on $Z \cap N^2$. ∎

We strengthen Theorem 4.11.13.

**Theorem 4.11.17** *Let $D \supseteq \bigcup \mathcal{F}$. There is an extension $\mathcal{F}' = \mathcal{F} \vee \mathcal{H}$ of $(\mathcal{F}, U)$ such that*

(i) *$\mathcal{H}$ is a filter of $2^{D \setminus U}$,*

(ii) *$\mathcal{F}'_{\setminus U} = \mathcal{H}$,*

(iii) *$\mathcal{F}'$ minimizes $\mu$ in $D$,*

(iv) *$\mathcal{F}'_{\setminus N^2} = \{D \setminus N^2\}$.*

**Proof.** By Theorem 4.11.13 there is an extension $\mathcal{F} \vee \mathcal{G}$ of $\mathcal{F}$ which satisfies (i), (ii) and (iii). Let $\mathcal{H} = \mathcal{G} \vee \{D \setminus N^2\}$. We show that $\mathcal{F}' = \mathcal{F} \vee \mathcal{H}$ is as required. By construction, $\mathcal{F}'$ satisfies (i), (ii) and (iv). For $X \in \mathcal{H}$, let

$$P(X) = \{Y \in \mathcal{G} \mid Y \cap N^2 = X \cap N^2\}.$$

Then $\{P(X) \mid X \in \mathcal{H}\}$ is a partition of $\mathcal{G}$. If $Y \in P(X)$, then by Lemma 4.11.16, $E(Y) = E(Y \cap N^2) = E(X \cap N^2) = E(X)$. Since $\mathcal{G}$ is a filter of $2^{D \setminus U}$, if $X \subseteq Y$, then $P(X)_{\setminus N^2} \subseteq P(Y)_{\setminus N^2}$. Since $P(X)_{\cap N^2} = \{X \cap N^2\}$ for each $X \in \mathcal{H}$, this implies that if $X \subseteq Y$, then $\big|P(X)\big| \le \big|P(Y)\big|$. Let

$$\mathcal{H}_{\ge n} = \{X \in \mathcal{H} \mid \big|P(X)\big| \ge n\}.$$

Then $\mathcal{H}_{\ge n}$ is a filter included in $\mathcal{H}$ for each $n$ and $\mathcal{H}_{\ge 1} = \mathcal{H}$. Let $N$ be the maximum value of $\big|P(X)\big|$. By assumption on $\mathcal{F} \vee \mathcal{G}$, $\mu_{\mathcal{F} \vee \mathcal{H}_{\ge n}} \ge \mu_{\mathcal{F} \vee \mathcal{G}}$. Suppose that for some $n$, $\mu_{\mathcal{F} \vee \mathcal{H}_{\ge n}} > \mu_{\mathcal{F} \vee \mathcal{G}}$. Then

$$
\begin{aligned}
\mu_{\mathcal{F} \vee \mathcal{G}} &= \frac{\sum_{X \in \mathcal{G}} \big|E(X)\big|}{\big|\mathcal{G}\big|} \\
&= \frac{\sum_{X \in \mathcal{H}} \big|P(X)\big| \cdot \big|E(X)\big|}{\sum_{X \in \mathcal{H}} \big|P(X)\big|} \\
&= \frac{\sum_{n=1}^{N} \sum_{X \in \mathcal{H}_{\ge n}} \big|E(X)\big|}{\sum_{n=1}^{N} \sum_{X \in \mathcal{H}_{\ge n}} 1} \\
&= \frac{\sum_{n=1}^{N} \left(\big|\mathcal{H}_{\ge n}\big| \mu_{\mathcal{F} \vee \mathcal{H}_{\ge n}}\right)}{\sum_{n=1}^{N} \big|\mathcal{H}_{\ge n}\big|} \\
&> \frac{\sum_{n=1}^{N} \left(\big|\mathcal{H}_{\ge n}\big| \mu_{\mathcal{F} \vee \mathcal{G}}\right)}{\sum_{n=1}^{N} \big|\mathcal{H}_{\ge n}\big|} \\
&= \mu_{\mathcal{F} \vee \mathcal{G}}(U) \cdot \frac{\sum_{n=1}^{N} \big|\mathcal{H}_{\ge n}\big|}{\sum_{n=1}^{N} \big|\mathcal{H}_{\ge n}\big|} \\
&= \mu_{\mathcal{F} \vee \mathcal{G}},
\end{aligned}
$$

which is impossible. Hence, for each $n$, $\mathcal{F} \vee \mathcal{H}_{\ge n}$ also minimizes $\mu$ in $D$. In particular, $\mathcal{F}'$ minimizes $\mu$ in $D$, as required. ∎

**Corollary 4.11.18** *There is an extension $\mathcal{F}' = \mathcal{N}_{\mathcal{F}}^3(U) \vee \mathcal{H}$ of $(\mathcal{N}_{\mathcal{F}}^3(U), U)$ such that*

*(i) $\mathcal{H}$ is a filter of $2^{N^2 \setminus U}$,*

*(ii) $\mathcal{F}'_{\setminus U} = \mathcal{H}$,*

*(iii) if $\mathcal{F}''$ is an extension of $(\mathcal{F}, U)$, then $\mu_{\mathcal{F}, \mathcal{F}''}(U) \ge \mu_{\mathcal{N}_{\mathcal{F}}^3(U), \mathcal{F}'}(U)$.*

**Proof.** Let $D \supseteq \bigcup \mathcal{F}$ and let $\mathcal{F} \vee \mathcal{G}$ be an extension of $(\mathcal{F}, U)$ satisfying the conditions given in Theorem 4.11.17. By Theorem 4.11.12, $\mathcal{F} \vee \mathcal{G}$ is an extension of $(\mathcal{N}_{\mathcal{F}}^3(U), U)$ and

$$\mu_{\mathcal{F}, \mathcal{F} \vee \mathcal{G}}(U) = \mu_{\mathcal{N}_{\mathcal{F}}^3(U), \mathcal{F} \vee \mathcal{G}}(U).$$

By condition (iv) of Theorem 4.11.17, $(\mathcal{F} \vee \mathcal{G})_{\cap N^2}$ is isomorphic to $\mathcal{F} \vee \mathcal{G}$. This and Lemma 4.11.16 imply that

$$\mu_{\mathcal{N}_{\mathcal{F}}^3(U), \mathcal{F} \vee \mathcal{G}}(U) = \mu_{\mathcal{N}_{\mathcal{F}}^3(U), (\mathcal{F} \vee \mathcal{G})_{\cap N^2}}(U).$$

Since $\mathcal{F}$ is an extension of $(\mathcal{N}_{\mathcal{F}}^3(U), U)$, $\mathcal{F} = \mathcal{N}_{\mathcal{F}}^3(U) \vee \mathcal{G}'$ for some (non-empty) union-closed family $\mathcal{G}'$ with $\bigcup \mathcal{G}' \cap U = \emptyset$. The family $\mathcal{H}' = \mathcal{G}' \vee \mathcal{G}$ is a filter of $2^{D \setminus U}$ with $\mathcal{H}'_{\setminus N^2} = \{D \setminus N^2\}$. Let $\mathcal{H} = \mathcal{H}'_{\cap N^2}$. Then

$$\begin{aligned}
(\mathcal{F} \vee \mathcal{G})_{\cap N^2} &= (\mathcal{N}_{\mathcal{F}}^3(U) \vee \mathcal{G}' \vee \mathcal{G})_{\cap N^2} \\
&= \mathcal{N}_{\mathcal{F}}^3(U) \vee \mathcal{H},
\end{aligned}$$

so that $\mathcal{F}' = \mathcal{N}_{\mathcal{F}}^3(U) \vee \mathcal{H}$ is the required extension. ∎

Corollary 4.11.18 implies that to determine the minimum possible value of $\mu_{\mathcal{F}, \mathcal{F}'}(U)$, $\mathcal{F}$ can be replaced by $\mathcal{N}_{\mathcal{F}}^3(U)$; thus assume that $\mathcal{F} = \mathcal{N}_{\mathcal{F}}^3(U)$.

Henceforth we assume that $U \in J(\mathcal{F})$. This implies that for every $X$ disjoint from $U$, $U \in E(X)$, so that $|E(X)| \geq 1$. To show that $(\mathcal{F}, U)$ has the density property, it suffices to show that for every extension $\mathcal{F}'$ of $(\mathcal{F}, U)$, $\mu_{\mathcal{F}'} \geq 2$.

**Theorem 4.11.19** *The following are equivalent:*

(i) *There exists an extension $\mathcal{F}'$ of $(\mathcal{F}, U)$ such that $\mu_{\mathcal{F}'} < 2$.*

(ii) *There is a filter $\mathcal{H}$ of $2^{N^2 \setminus U}$ such that $\mu_{\mathcal{F} \vee \mathcal{H}} < 2$ and for every minimal member $X$ of $\mathcal{H}$, $|E(X)| = 1$.*

**Proof.** Assertion (ii) is a special case of (i). Suppose that (i) holds. By Corollary 4.11.18, there is an extension $\mathcal{F} \vee \mathcal{H}$ of $(\mathcal{F}, U)$ such that $\mathcal{H}$ is a filter of $2^{N^2 \setminus U}$ and $\mu_{\mathcal{F} \vee \mathcal{H}} < 2$. Let $\mathcal{H}$ be a minimal filter of $2^{N^2 \setminus U}$ such that $\mu_{\mathcal{F} \vee \mathcal{H}} < 2$. Suppose that there is a minimal member $X$ of $\mathcal{H}$ such that $E(X) \geq 2$. Let $\mathcal{H}' = \mathcal{H} \setminus \{X\}$. Note that the assumption on $X$ and $\mu_{\mathcal{F} \vee \mathcal{H}} < 2$ imply that $\mathcal{H}' \neq \emptyset$. The family $\mathcal{H}'$ is a filter of $2^{N^2 \setminus U}$ and

$$\begin{aligned}
\mu_{\mathcal{F} \vee \mathcal{H}'} &= \frac{\sum_{Y \in \mathcal{H}'} |E(Y)|}{|\mathcal{H}'|} \\
&= \frac{\left(\sum_{Y \in \mathcal{H}} |E(Y)|\right) - |E(X)|}{|\mathcal{H}'| - 1} \\
&< \frac{\sum_{Y \in \mathcal{H}} |E(Y)|}{|\mathcal{H}|} \\
&= \mu_{\mathcal{F} \vee \mathcal{H}} \\
&< 2,
\end{aligned}$$

where we used the fact that if $a, b > 0$, $c > 1$ and $a/c < b$, then $(a-b)/(c-1) < a/c$. This contradicts the minimality assumption on $\mathcal{H}$, so that $\mathcal{H}$ is as desired. ∎

If there is no filter $\mathcal{H}$ satisfying the conditions in assertion (ii) of Theorem 4.11.19, then $\mu_{\mathcal{F}'}(U) \geq 2$ for every extension $\mathcal{F}'$ of $(\mathcal{F}, U)$, so that $(\mathcal{F}, U)$ has the density property. This yields:

**Observation 4.11.20** *If every extension $\mathcal{F} \vee \mathcal{H}$ of $(\mathcal{F}, U)$ such that*

(i) $\mathcal{H}$ *is a filter of* $2^{N^2 \setminus U}$,

(ii) *for every minimal* $X \in \mathcal{H}$, $\left| E(X) \right| = 1$

*satisfies $\mu_{\mathcal{F} \vee \mathcal{H}} \geq 2$, then $(\mathcal{F}, U)$ has the density property.*

Let $\mathcal{H}$ be an arbitrary filter satisfying (i) and (ii) of Observation 4.11.20 and let $\mathcal{F}' = \mathcal{F} \vee \mathcal{H}$. We now assume that $J(\mathcal{F})$ is a graph. Thus, for some $a, b \in \bigcup \mathcal{F}$, $U = \{a, b\}$. Assume that $a \neq b$ and let

$$
\begin{aligned}
N_a &= \{x \in N \setminus U \mid \{x, a\} \in J(\mathcal{F}) \text{ and } \{x, b\} \notin J(\mathcal{F})\}, \\
N_b &= \{x \in N \setminus U \mid \{x, a\} \notin J(\mathcal{F}) \text{ and } \{x, b\} \in J(\mathcal{F})\}, \\
N_{ab} &= \{x \in N \setminus U \mid \{x, a\} \in J(\mathcal{F}) \text{ and } \{x, b\} \in J(\mathcal{F})\}.
\end{aligned}
$$

Then $N \setminus U = N_a \cup N_b \cup N_{ab}$ (see Figure 4.5).

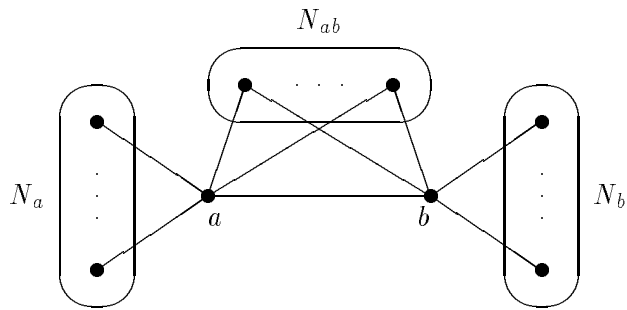

Figure 4.5

**Lemma 4.11.21** *If $X \in \mathcal{H}$, then $X \cap N_a \neq \emptyset$ and $X \cap N_b \neq \emptyset$.*

**Proof.** Let $X \in \mathcal{H}$. Let $X'$ be a minimal member of $\mathcal{H}$ such that $X' \subseteq X$. Suppose $X' \cap N_a = \emptyset$. Then either $\{b\} \in E(X')$ (if $X' \cap (N_{ab} \cup N_b) \neq \emptyset$) or $\emptyset \in E(X')$ (if $X' \cap (N_{ab} \cup N_b) = \emptyset$). Since $\{a, b\} \in E(X')$, this contradicts $\left| E(X') \right| = 1$. Thus $X' \cap N_a \neq \emptyset$. By symmetry, $X' \cap N_b \neq \emptyset$ and we are done. ∎

**Corollary 4.11.22** *If $X \in \mathcal{H}$ and $Y \subseteq U$, then $\pi_{\mathcal{F}}(X \cup Y) \cap U = Y$.*

**Proof.** Let $X \in \mathcal{H}$ and $Y \subseteq U$. Suppose that $a \in Y$. By Lemma 4.11.21, there exists $x \in N_a \cap X$, so that $\{x, a\} \in \mathcal{F}$ and $\{x, a\} \subseteq X \cup Y$, which implies that $a \in \pi_{\mathcal{F}}(X \cup Y)$. Similarly, if $b \in Y$, then $b \in \pi_{\mathcal{F}}(X \cup Y)$, as required. ∎

**Definition.** Let $Y \subseteq U$ and $x \in N \setminus U$. Let $\mathcal{E}(Y, x)$ consist of the subsets $X$ of $N^2 \setminus U$ such that there is an edge $\{x, y\} \in J(\mathcal{F})$ with $y \in X \cup Y$ or $y = x$. Define

$$\mathcal{E}(Y) = \bigcap_{x \in N \setminus U} \mathcal{E}(Y, x).$$

**Observation 4.11.23** *For $Y \subseteq U$ and $x \in N \setminus U$, $\mathcal{E}(Y, x)$ is a filter of $2^{N^2 \setminus U}$.*

**Observation 4.11.24** *If $X \in \mathcal{E}(Y)$, then $X$ and $Y$ satisfy condition (ii)″ of Lemma 4.11.10.*

**Lemma 4.11.25** *If $X \in \mathcal{E}(Y) \cap \mathcal{H}$, then $Y \in E(X)$.*

**Proof.** Let $X \in \mathcal{E}(Y) \cap \mathcal{H}$. By Observation 4.11.24, $X$ and $Y$ satisfy (ii)″ of Lemma 4.11.10. By Corollary 4.11.22, $X$ and $Y$ satisfy (i) of Lemma 4.11.7. The result follows by the proof of Lemma 4.11.10. ∎

**Definition.** Let $\mathcal{U}$ be a filter of $2^S$. Define

$$\nu(\mathcal{U}) = \frac{|\mathcal{U}|}{2^{|S|}}.$$

Thus $\nu(\mathcal{U})$ is the density of $\mathcal{U}$ in $2^S$.

**Theorem 4.11.26**

$$\mu_{\mathcal{F}, \mathcal{F}'}(U) \geq 1 + \sum_{Y \subset U} \prod_{x \in N \setminus U} \nu\left(\mathcal{E}(Y, x)\right).$$

**Proof.** By Kleitman's Lemma 2.6.3, for every $Y \subseteq U$,

$$\nu(\mathcal{H} \cap \mathcal{E}(Y)) \geq \nu(\mathcal{H})\nu(\mathcal{E}(Y)).$$

Lemma 4.11.25 and the fact that $U \in E(X)$ for every $X \in \mathcal{H}$ yield

$$
\begin{aligned}
\sum_{X \in \mathcal{H}} \left| E(X) \right| &= \sum_{Y \subseteq U} \left| \{ X \in \mathcal{H} \mid Y \in E(X) \} \right| \\
&= \left| \{ X \in \mathcal{H} \mid U \in E(X) \} \right| + \sum_{Y \subset U} \left| \{ X \in \mathcal{H} \mid Y \in E(X) \} \right| \\
&\geq \left| \mathcal{H} \right| + \sum_{Y \subset U} \left| \mathcal{H} \cap \mathcal{E}(Y) \right|.
\end{aligned}
$$

Using $\mathcal{H} = \mathcal{F}'_{\setminus U}$ we get

$$
\begin{aligned}
\mu_{\mathcal{F}'} &= \frac{\sum_{X \in \mathcal{H}} \left| E(X) \right|}{|\mathcal{H}|} \\
&\geq \frac{|\mathcal{H}| + \sum_{Y \subset U} \left| \mathcal{H} \cap \mathcal{E}(Y) \right|}{|\mathcal{H}|} \\
&= 1 + \sum_{Y \subset U} \frac{\nu(\mathcal{H} \cap \mathcal{E}(Y))}{\nu(\mathcal{H})} \\
&\geq 1 + \sum_{Y \subset U} \nu(\mathcal{E}(Y)).
\end{aligned}
$$

Multiple applications of Kleitman's Lemma yield

$$
\nu(\mathcal{E}(Y)) \geq \prod_{x \in N \setminus U} \nu(\mathcal{E}(Y, x))
$$

and the result follows. ∎

Define
$$
J'(\mathcal{F}) = \{ V \in J(\mathcal{F}) \mid |V| = 2 \}.
$$

For $x \in N \setminus U$, define

$$
A(x) = \left\{ \{x, y\} \mid y \notin U \text{ and } \{x, y\} \in J'(\mathcal{F}) \right\}.
$$

Let $n(x) = \left| A(x) \right|$.

**Lemma 4.11.27** *Let $x \in N \setminus U$. Then*

    *(i)* $\nu(\mathcal{E}(\emptyset, x)) \geq 1 - (\frac{1}{2})^{n(x)}$.

    *(ii) If $x \in N_a$, then $\nu(\mathcal{E}(\{a\}, x)) = 1$ and $\nu(\mathcal{E}(\{b\}, x)) \geq 1 - (\frac{1}{2})^{n(x)}$.*

    *(iii) If $x \in N_b$, then $\nu(\mathcal{E}(\{a\}, x)) \geq 1 - (\frac{1}{2})^{n(x)}$, and $\nu(\mathcal{E}(\{b\}, x)) = 1$.*

    *(iv) If $x \in N_{ab}$, then $\nu(\mathcal{E}(\{a\}, x)) = \nu(\mathcal{E}(\{b\}, x)) = 1$.*

**Proof.** Consider $\nu(\mathcal{E}(\emptyset, x))$. If $\{x\} \in J(\mathcal{F})$, then $\nu(\mathcal{E}(\emptyset, x)) = 1$. If $\{x\} \notin J(\mathcal{F})$, then

$$
2^{N^2 \setminus U} \setminus \mathcal{E}(\emptyset, x) = 2^{N^2 \setminus (U \cup A(x))}.
$$

We have
$$
\left| 2^{N^2 \setminus (U \cup A(x))} \right| = 2^{\left| N^2 \setminus U \right|} 2^{-\left| A(x) \right|},
$$

which implies that

$$
\nu(\mathcal{E}(\emptyset, x)) \geq 1 - (\tfrac{1}{2})^{n(x)}.
$$

The remaining cases are proved similarly. ∎

Recall that $\mathrm{d}_{\mathcal{G}}(U)$ is the degree in the graph $\mathcal{G}$ of $U$. If $U$ has minimal degree in $J'(\mathcal{F})$, then $(\mathcal{F}, U)$ has the density property:

**Theorem 4.11.28** *If for every edge $V \in J'(\mathcal{F})$ with $V \cap U \neq \emptyset$, $\mathrm{d}_{J'(\mathcal{F})}(V) \geq \mathrm{d}_{J'(\mathcal{F})}(U)$, then for every extension $\mathcal{F}'$ of $(\mathcal{F}, U)$, $\mu_{\mathcal{F}, \mathcal{F}'}(U) \geq 2$.*

**Proof.** We can assume that $\mathcal{F}' = \mathcal{F} \vee \mathcal{H}$ where $\mathcal{H}$ is a filter satisfying (i) and (ii) of Observation 4.11.20. Let $n = \mathrm{d}_{J'(\mathcal{F})}(U)$, $n_a = \lfloor N_a \rfloor$, $n_b = \lfloor N_b \rfloor$ and $n_{ab} = \lfloor N_{ab} \rfloor$. Then

$$n = n_a + n_b + 2n_{ab}.$$

Let $x \in N_a$. We have $\mathrm{d}_{J'(\mathcal{F})}(\{a, x\}) \geq n$ and $\mathrm{d}_{J'(\mathcal{F})}(\{a, x\}) = n_a + n_{ab} + n(x)$. This gives $n(x) \geq n - n_a - n_{ab} = n_b + n_{ab}$. Similarly, if $x \in N_b$, then $n(x) \geq n_a + n_{ab}$. Lemma 4.11.21 and the assumptions on $\mathcal{H}$ imply that $n_a \geq 1$ and $n_b \geq 1$. By Theorem 4.11.26 and Lemma 4.11.27,

$$
\begin{aligned}
\mu_{\mathcal{F}'} &\geq 1 + \prod_{x \in N \setminus U} \nu(\mathcal{E}(\{a\}, x)) + \prod_{x \in N \setminus U} \nu(\mathcal{E}(\{b\}, x)) \\
&\geq 1 + \prod_{x \in N_b} \left(1 - (\tfrac{1}{2})^{n(x)}\right) + \prod_{x \in N_a} \left(1 - (\tfrac{1}{2})^{n(x)}\right) \\
&\geq 1 + \left(1 - (\tfrac{1}{2})^{n_a}\right)^{n_b} + \left(1 - (\tfrac{1}{2})^{n_b}\right)^{n_a}.
\end{aligned}
$$

Let $c_1 = \left(1 - (\tfrac{1}{2})^{n_a}\right)^{n_b}$ and $c_2 = \left(1 - (\tfrac{1}{2})^{n_b}\right)^{n_a}$. It remains to show that $1 + c_1 + c_2 \geq 2$. Without loss of generality, assume that $n_a \leq n_b$. If $n_a = 1$, then

$$1 + c_1 + c_2 = 1 + (\tfrac{1}{2})^{n_b} + 1 - (\tfrac{1}{2})^{n_b} = 2.$$

Suppose that $n_a \geq 2$. Observe that $(\tfrac{3}{2})^m \geq m$ for $m \geq 1$. Using the fact that $(1 - a)^b \geq 1 - ba$ for $0 \leq a \leq 1$ and $1 \leq b$ (Lemma 4.11.29), we obtain

$$
\begin{aligned}
1 + c_1 + c_2 &\geq 2 + (\tfrac{3}{4})^{n_b} - n_a (\tfrac{1}{2})^{n_b} \\
&= 2 + (\tfrac{1}{2})^{n_b} \left((\tfrac{3}{2})^{n_b} - n_a\right) \\
&\geq 2 + (\tfrac{1}{2})^{n_b} \left((\tfrac{3}{2})^{n_a} - n_a\right) \\
&\geq 2,
\end{aligned}
$$

as required. $\blacksquare$

**Lemma 4.11.29** *If $0 \leq a \leq 1$ and $1 \leq b$, then $(1 - a)^b \geq 1 - ba$.*

**Proof.** The result is true for $a = 0$. Differentiating both expressions relative to $a$ yields

$$\frac{d}{da}(1 - a)^b = -b(1 - a)^{b-1} \geq -b = \frac{d}{da}(1 - ba).$$

The result follows. $\blacksquare$

If $\mathcal{J}'(\mathcal{F}) = \emptyset$, then $\mathcal{F}$ is generated by one-element sets, so that $\mathcal{F}$ is a Boolean lattice. Thus by Observation 4.11.20, Theorem 4.11.2 has been proved. In fact we have the following stronger result:

**Theorem 4.11.30** *Let $\mathcal{F}$ be a union-closed family of sets and $U \in J(\mathcal{F})$. Let $\mathcal{G} = \{V \in J(\mathcal{F}) \mid V \cap U \neq \emptyset\}$. If*

*(i)* $|U| = 2$,

*(ii)* $\mathcal{G}$ *is a graph,*

*(iii) there is a simple graph $\mathcal{G}' \subseteq J(\mathcal{F})$ such that for every $V \in \mathcal{G}$ with $|V| = 2$, $V \in \mathcal{G}'$ and $d_{\mathcal{G}'}(V) \geq d_{\mathcal{G}'}(U)$,*

*then $(\mathcal{F}, U)$ has the density property.*

**Proof.** Suppose that $\mathcal{N}$ is a union-closed family such that $\mathcal{F} \supseteq \mathcal{N} \supseteq \mathcal{N}_{\mathcal{F}}(U)$. Then $\mathcal{F}$ is an extension of $(\mathcal{N}, U)$ (see the proof of Lemma 4.11.11). Lemma 4.11.10 shows that $E_{\mathcal{F},U}(X) \supseteq E_{\mathcal{N},U}(X)$ for every $X$ with $X \cap U = \emptyset$. Therefore

$$\mu_{\mathcal{F},\mathcal{F}'}(U) \geq \mu_{\mathcal{N},\mathcal{F}'}(U)$$

for every extension $\mathcal{F}'$ of $(\mathcal{F}, U)$. Let $\mathcal{N}$ be the union-closed family generated by $\mathcal{G} \cup \mathcal{G}' \cup \{\emptyset\}$. Then $\mathcal{F} \supseteq \mathcal{N} \supseteq \mathcal{N}_{\mathcal{F}}(U)$. Theorems 4.11.19 and 4.11.28 show that $\mu_{\mathcal{N},\mathcal{F}'}(U) \geq 2$ for every extension $\mathcal{F}'$ of $(\mathcal{F}, U)$. It follows that $(\mathcal{F}, U)$ has the density property. ∎

## 4.12 Notes

**Section 4.1.** The theory of matchings is a well developed area of combinatorics. One of the earliest and best known results is the theorem of König [28] and Hall [21]. It characterizes the families of sets $\mathcal{F}$ for which there is an injective map $f : \mathcal{F} \to \bigcup \mathcal{F}$ such that $f(U) \in U$ for each $U \in \mathcal{F}$ (Theorem 2.2.1 in Anderson [3]).

Decreasing and increasing matchings in lattices are discussed by Duffus [13] and Kung [29], [30]. Kung [29] shows that if $L$ is a modular lattice, then there is a decreasing matching $\sigma : M(L) \to J(L)$, thereby answering a question posed by Rival.

The union-closed sets conjecture (Conjecture 4.1.5) is stated in Duffus [14] and appears as an (open) problem in Stanley [41] (Exercise 39 pg. 161), where it is called "diabolical". The problem is well known. Winkler [45] writes that it is "one of the most embarrassing gaps in combinatorial knowledge".

Since the problem was posed, virtually no progress has been made toward its solution. It is known to be true for modular lattices (R. McKenzie personal communication) and geometric lattices (as mentioned in Duffus [14]). Neither of these results has been published (they follow from Corollary 4.7.6 and Corollary 4.7.10). It is also known to be true for small union-closed families of sets ($|\mathcal{F}| \lesssim 20$ and $|\bigcup \mathcal{F}| \lesssim 8$) and for some other special cases (Sarvate and Renaud [38], [39] and other unpublished work).

**Section 4.2.** Theorems 4.2.3 and 4.2.5 generalize some well-known elementary results related to the union-closed sets conjecture (A. Ehrenfeucht, M. Main personal communications).

**Section 4.3.** The intersection-closed representations of lattices appear (implicitly) in most textbooks on lattice theory.

Conjectures 4.3.7 and 4.3.6 are the usual way of formulating the union-closed sets conjecture.

**Section 4.4.** The subdirect product is a useful generalization of the direct product. Usually, trivial subdirect products are excluded from the definition. The concept is discussed in Grätzer [17] from a lattice-theoretic persective.

**Section 4.6.** Locality in lattices is usually defined in terms of intervals. For example, a lattice $L$ is *locally distributive* iff for every $x, y \in L$ such that $y$ is the join of the atoms of $[x)$, the interval $[x, y]$ is a Boolean lattice (see Greene and Kleitman [18]).

**Section 4.7.** Most of the special classes of lattices for which the union-closed sets conjecture is known to be true appear in this section. The lattice-theoretic proof techniques used here are standard.

**Section 4.8.** There are many unsolved density related problems. Some examples follow.

**Problem 4.12.1** (Colbourn and Rival, see Sands [37]) *Does there exist $r > 2$ such that every finite distributive lattice $L$ contains a join-irreducible $a$ with $\frac{1}{r} \leq \frac{|[a)|}{|L|} \leq 1 - \frac{1}{r}$?*

Observe that Theorem 4.2.5 shows that in a distributive lattice, every minimal join-irreducible has density at least $\frac{1}{2}$ and every maximal join-irreducible has density at most $\frac{1}{2}$.

**Problem 4.12.2** (Daykin and Frankl [9]) *Is it true that for every convex subfamily $\mathcal{C}$ of $2^X$, $\frac{w(\mathcal{C})}{|\mathcal{C}|} \geq \frac{w(2^X)}{|2^X|}$?*

Note that $\frac{w(P)}{|P|}$ is the density in $P$ of a Sperner antichain of $P$. Daykin [9] gives a powerful strengthening of Problem 4.12.1.

**Problem 4.12.3** (Daykin [9]) *Which posets $P$ have the property that for every $Q \subseteq P$, the density of $Q$ in $P$ is less than or equal to the ratio of the number of maximal chains of $P$ intersecting $Q$ to the number of maximal chains of $P$?*

**Problem 4.12.4** (Knill) *Is it true that for every $\epsilon > 0$, there exists $w > 0$ such that if $L$ is a lattice with $\mathrm{w}(L) \geq w$, then there is a join-irreducible $a \in L$ such that $\frac{\mathrm{w}([a))}{\mathrm{w}(L)} \leq \frac{1}{2} + \epsilon$? Is this true for distributive lattices?*

(The ratio $\frac{\mathrm{w}([a))}{\mathrm{w}(L)}$ is the width-density of $[a)$ in $L$.)

The computations of Observation 4.8.2 can be continued for some other small values of $n$. For $P = [1]$, J.-C. Renaud has computed values up to $n = 9$ (unpublished).

**Section 4.9.** Bounds on the minimum value of $m$ such that $L$ has the $[n]$-density property for all $n \geq m$, are implicit in the proof of Theorem 4.9.1 but are quite large. Lemma 4.9.2 appears in Stanley [41]. The Zeta polynomial finds many applications in the study of the combinatorics of posets.

**Section 4.10.** Corollary 4.10.6 and Theorem 4.10.5 generalize the known (but unpublished) bound $\overline{h}(\mathcal{L}, [1], n) \leq Cn \log_2(n)$.

**Section 4.11.** The class of lattices shown to have the density property in Theorem 4.11.2 is a subclass of the class of 2-distributive lattices. The 1-distributive identity is given by

$$(u_1 \vee u_2) \wedge v = (u_1 \wedge v) \vee (u_2 \wedge v).$$

The 2-distributive identity is given by

$$(u_1 \vee u_2 \vee u_3) \wedge v = ((u_1 \vee u_2) \wedge v) \vee ((u_1 \vee u_3) \wedge v) \vee ((u_2 \vee u_3) \wedge v).$$

(This is the dual of the (non-equivalent) identity given in Grätzer [17], pg 219.) Graph-generated union-closed families $\mathcal{F}$ with $\emptyset \in \mathcal{F}$ satisfy the 2-distributive identity. Proof: If $a \in J(\mathcal{F})$, then $|a| \leq 2$, which implies that $a \subseteq u_1 \cup u_2 \cup u_3$ iff $a \subseteq u_1 \cup u_2$ or $a \subseteq u_1 \cup u_3$ or $a \subseteq u_2 \cup u_3$. Since $u \wedge v = \bigcup \{a \in J(\mathcal{F}) \mid a \subseteq u \cap v\}$, the result follows.

Huhn [23] discusses $n$-distributive modular lattices.

The estimates of $\mu_{\mathcal{F}'}$ can be generalized for $P$-densities with $P \neq [1]$ by using the generalization of Kleitman's Lemma to distributive lattices given in Anderson [3] pp 91–94. However, the resulting estimates are in general not strong enough to prove the $P$-density property for graph generated union-closed families of sets.

Computations similar to the ones in this section suggest that $(\mathcal{F}, U)$ has the density property if for each $V \in J(\mathcal{F})$ with $V \cap U \neq \emptyset$, $\mathrm{d}_{J(\mathcal{F})}(V)$ is large enough compared to $\mathrm{d}_{J(\mathcal{F})}(U)$. In some cases $\mathrm{d}_{J(\mathcal{F})}(V) \gtrsim r2^r \log \mathrm{d}_{J(\mathcal{F})}(U)$ suffices, where $r = \max\{|V| \mid V \in J(\mathcal{F})\}$.